\documentclass[11pt,a4paper]{article}%

\usepackage{amsmath}
\usepackage{amsfonts}
\usepackage{amssymb}
\usepackage{graphicx}%
\usepackage{float}
\usepackage{enumitem,hyperref}
\hypersetup{
    colorlinks=true,
    linkcolor=blue,
    filecolor=magenta,      
    urlcolor=cyan,
    citecolor=blue,
}
\usepackage{xcolor}
\providecommand{\U}[1]{\protect\rule{.1in}{.1in}}
\newtheorem{theorem}{Theorem}[section]

\newtheorem{corollary}[theorem]{Corollary}
\newtheorem{definition}[theorem]{Definition}

\newtheorem{lemma}[theorem]{Lemma}

\newtheorem{remark}[theorem]{Remark}

\usepackage[left=3cm,right=3cm,top=3cm,bottom=3cm]{geometry}

\makeatletter
\let\orgdescriptionlabel\descriptionlabel
\renewcommand*{\descriptionlabel}[1]{%
  \let\orglabel\label
  \let\label\@gobble
  \phantomsection
  \edef\@currentlabel{#1\unskip}%
  \let\label\orglabel
  \orgdescriptionlabel{#1}%
}
\makeatother

\title{\textbf{Estimation and testing on independent not identically distributed observations
based on R\'{e}nyi's pseudodistances}}
\author{Elena Castilla$^{1}$, Maria Jaenada$^{1}$ and Leandro Pardo$^{1}$\\$^{1}${\small Department of Statistics and O.R., Complutense University of
Madrid,  Spain}}
\date{}
\begin{document}

\maketitle

\begin{abstract}
In real life we often deal with independent but not identically distributed observations (i.n.i.d.o), for which the most well-known statistical model is the multiple linear regression model (MLRM) without random covariates. While the classical methods are based on the maximum likelihood estimator (MLE), it is well known its lack of robustness to small deviations from the assumed conditions. In this paper, and based on the R\'{e}nyi's pseudodistance (RP), we introduce a new family of estimators in case our information about the unknown parameter is given for i.n.i.d.o.. This family of estimators, let say minimum RP estimators (as they are obtained by minimizing the RP between the assumed distribution and the empirical distribution of the data), contains the MLE as a particular case and can be applied, among others, to the MLRM without random covariates. Based on these estimators, we introduce Wald-type tests for testing simple and composite null hypotheses, as an extension of the classical MLE-based Wald test. Influence functions for the estimators and Wald-type tests are also obtained and analysed.  Finally, a simulation study is developed in order to asses the performance of the proposed methods and some real-life data are analysed for illustrative purpose.
\end{abstract}

\bigskip

%\noindent\underline{\textbf{AMS 2001 Subject Classification}}\textbf{: }


\section{Introduction\label{sec1}}

In parametric estimation the role of divergence measures is very intuitive: minimizing a suitable divergence measure between the data and the assumed model  in order to estimate the unknown parameters. These estimators are called ``minimum divergence estimators'' (MDEs). There is a growing body of literature that recognizes the importance of MDEs because their robustness, without a significant loss of efficiency, in relation to the maximum likelihood estimator (MLE). See, for instance, Beran \cite{Beran},  Tamura and Boos \cite{Tamura},  Simpson \cite{Simpson1,Simpson2}, Lindsay \cite{Lindsay}, Pardo \cite{Pardo}  and Basu et al. \cite{Basu2011}. In the case of continuous models is convenient to consider families of divergence measures for which non-parametric estimators of the unknown density function are not needed.  From this perspective,  the density power divergence (DPD) family, leading to the minimum DPD estimators, is the most important family of divergence measures. For more details see \cite{Basu2011}. However, there is another important family of divergence measures which neither needs non-parametric estimators, the R\'{e}nyi's pseudodistances (RP).

Let $X_{1},...,X_{n}$ be a random sample of size $n$ from a population having true and unknown density function $g,$ modelled by a parametric family of densities $f_{\boldsymbol{\theta}}$ with $\boldsymbol{\theta}\in\Theta\subset\mathbb{R}^{p}.$ The RP between the densities $f_{\boldsymbol{\theta}}$ and $g$ is given , for $\alpha>0,$ by
\begin{align}
R_{\alpha}\left(  f_{\boldsymbol{\theta}},g\right)  =&\frac{1}{\alpha+1}%
\log\left(  \int f_{\boldsymbol{\theta}}(x)^{\alpha+1}dx\right) \notag \\
&  +\frac
{1}{\alpha\left(  \alpha+1\right)  }\log\left(  \int g(x)^{\alpha+1}dx\right)
-\frac{1}{\alpha}\log\left(  \int f_{\boldsymbol{\theta}}(x)^{\alpha
}g(x)dx\right).  \label{1.1}%
\end{align}
The RP can be defined for $\alpha=0$ taking continuous limits, yielding the expression
\[
R_{0}\left(  f_{\boldsymbol{\theta}},g\right)  =\lim_{\alpha\downarrow
0}R_{\alpha}\left(  f_{\boldsymbol{\theta}},g\right)  =\int g(x)\log
\frac{g(x)}{f_{\theta}(x)}dx,
\]
i.e., the RP coincides with the Kullback-Leibler divergence  between $g$ and $f_{\boldsymbol{\theta}}$, $D_{Kullback}(g,f_{\boldsymbol{\theta}})$, at $\alpha=0$ (see \cite{Pardo}).  The RP was considered for the first time in Jones et al. \cite{Jones}; later Broniatowski et al. \cite{Broniatowski} established that the RP is positive for any two densities and for all values of the parameter $\alpha,$  $R_{\alpha}\left(f_{\boldsymbol{\theta}},g\right)  \geq0$ and further $R_{\alpha}\left(f_{\boldsymbol{\theta}},g\right)  =0$ if and only if $f_{\boldsymbol{\theta}}=g.$ This property embolden the definition of the minimum RP estimator as the minimizer of the RP between the assumed distribution and the empiric distribution of the data. Therefore, the minimum RP estimator based on the random sample $X_{1},...,X_{n}$ for the unknown parameter $\boldsymbol{\theta}$  is given, for $\alpha>0$, by
\begin{equation}
	\widehat{\boldsymbol{\theta}}_{\alpha}^{\ast}=\arg\sup_{\boldsymbol{\theta}%
		\in\Theta}%
	%TCIMACRO{\tsum \limits_{i=1}^{n}}%
	%BeginExpansion
	{\textstyle\sum\limits_{i=1}^{n}}
	%EndExpansion
	\frac{f_{\boldsymbol{\theta}}(X_{i})^{\alpha}}{C_{\alpha}(\boldsymbol{\theta
			)}}, \label{1.2}%
\end{equation}
where%
\[
C_{\alpha}(\boldsymbol{\theta)=}\left(
%TCIMACRO{\tint }%
%BeginExpansion
{\textstyle\int}
%EndExpansion
f_{\boldsymbol{\theta}}(x)^{\alpha+1}dx\right)  ^{\frac{\alpha}{\alpha+1}}.%
\]
Note that the value $\alpha=0$ was defined as the KL divergence and hence, the minimum RP estimator coincides with the MLE at $\alpha=0.$ Besides, \cite{Broniatowski} studied the asymptotic properties and robustness of the minimum RP estimators and presented an application to the multiple linear regression model (MLRM) with random covariates. In the same vein, \cite{Castilla1} introduced Wald-type tests based on the minimum RP estimators for the MLRM and \cite{Castilla2} studied the minimum RP estimator for the linear regression model in the ultra-high dimensional set-up. Moreover, $\widehat{\boldsymbol{\theta}}_{\alpha}^{\ast}$ is a M-estimator and thus it asymptotic distribution and influence function (IF) can be obtained based on the asymptotic theory of the M-estimators. 

However, far too little attention has been paid to the case of independent but non identically distributed observations (i.n.i.d.o.), for which the most well-known statistical model is the MLRM  without random covariates.  The nicest study for i.n.i.d.o. based on divergence measures, until now, is the paper of Ghosh and Basu \cite{Ghosh2013} based on DPD measures. Some extensions are given in \cite{Ghosh2018} and \cite{Basu2018}. The main aim of this paper is to introduce and study the minimum RP estimator for i.n.i.d.o. We study their asymptotic properties as well as we obtain its influence function in order to study its robustness. In Section \ref{sec:RP} we introduce the minimum RP estimator for i.n.i.d.o. The consistency and asymptotic distribution is presented in Section \ref{sec:CONS}. Section \ref{sec:Wald} is devoted to introduce and study Wald-type tests based on minimum RP estimators. The robustness of these estimators and Wald-type test is studied through its influence functions in Section \ref{sec:IF}. The special case of MLRM is considered in Section \ref{sec:MRM}. Finally, an extensively study and two numerical examples  of the MLRM   are presented in Section \ref{sec:simulation} and \ref{sec:data}, respectively.

\section{The minimum RP estimator for independent but not identically distributed observations \label{sec:RP}}

Let $Y_{1},...,Y_{n}$ be i.n.i.d.o. being $g_{1},...,g_{n}$ the corresponding density functions with respect to some common dominating measure. We are interested in modeling $g_{i}$ by the density function $f_{i}(y,\boldsymbol{\theta}),$ $i=1,...,n,$ being $\boldsymbol{\theta}$ common for
all the density functions $f_{i}(y,\boldsymbol{\theta})$. For each observation $i,$ the ``R\'{e}nyi's pseudodistance''  between
$f_{i}(y,\boldsymbol{\theta})$ and $\widehat{g}_{i},$ can be defined for positives values of $\alpha$ as
\begin{equation}
R_{\alpha}\left(  f_{i}(y,\boldsymbol{\theta}),\widehat{g}_{i}\right)
=\frac{1}{\alpha+1}\log\left(  \int f_{i}(y,\boldsymbol{\theta})^{\alpha
+1}dy\right)  -\frac{1}{\alpha}\log\left(  \int f_{i}(y,\boldsymbol{\theta
})^{\alpha}\widehat{g}_{i}(y)dy\right)  +k, \label{2.1}%
\end{equation}
where
\[
k=\frac{1}{\alpha\left(  \alpha+1\right)  }\log\left(  \int\widehat{g}%
_{i}(y)^{\alpha+1}dy\right)
\]
does not depend on $\boldsymbol{\theta.}$ As we only have one observation of $Y_{i}$ the best way to estimate $\widehat{g}_{i}$ is assuming that the distribution is degenerate in $y_{i}.$
Therefore, (\ref{2.1}) yields to the loss
\begin{equation}
R_{\alpha}\left(  f_{i}(y,\boldsymbol{\theta}),\widehat{g}_{i}\right)
=\frac{1}{\alpha+1}\log\left(  \int f_{i}(y,\boldsymbol{\theta})^{\alpha
+1}dy\right)  -\frac{1}{\alpha}\log f_{i}(Y_{i},\boldsymbol{\theta})^{\alpha
}+k. \label{2.2}%
\end{equation}
At $\alpha=0$ the RP loss is given by
\begin{equation}
R_{0}\left(  f_{i}(y,\boldsymbol{\theta}),\widehat{g}_{i}\right)
=\lim_{\alpha\downarrow0}R_{\alpha}\left(  f_{i}(y,\boldsymbol{\theta
}),\widehat{g}_{i}\right)  =-\log f_{i}(Y_{i},\boldsymbol{\theta})+k.
\label{2.3}%
\end{equation}
Now, Expression (\ref{2.2}) can be written as
\[
R_{\alpha}\left(  f_{i}(y,\boldsymbol{\theta}),\widehat{g}_{i}\right)
=-\frac{1}{\alpha}\log\frac{f_{i}(Y_{i},\boldsymbol{\theta})^{\alpha}}{\left(
\int f_{i}(y,\boldsymbol{\theta})^{\alpha+1}dy\right)  ^{\frac{\alpha}%
{\alpha+1}}}+k,
\]
and thus minimizing $R_{\alpha}\left(  f_{i}(y,\boldsymbol{\theta}),\widehat{g}%
_{i}\right)  $ in $\boldsymbol{\theta}$, for $\alpha>0,$ is equivalent to
maximizing%
\[
V_{i}(Y_{i},\boldsymbol{\theta)=}\frac{f_{i}(Y_{i},\boldsymbol{\theta
})^{\alpha}}{\left(  \int f_{i}(y,\boldsymbol{\theta})^{\alpha+1}dy\right)
^{\frac{\alpha}{\alpha+1}}}.
\]
In the following, we shall denote,%
\[
L_{\alpha}^{i}\left(  \boldsymbol{\theta}\right)  =\left(  \int f_{i}%
(y,\boldsymbol{\theta})^{\alpha+1}dy\right)  ^{\frac{\alpha}{\alpha+1}}.
\]
Based on this idea we are going to consider the objective function%
\begin{equation}
H_{n}^{\alpha}(\boldsymbol{\theta})=\frac{1}{n}%
%TCIMACRO{\tsum \limits_{i=1}^{n}}%
%BeginExpansion
{\textstyle\sum\limits_{i=1}^{n}}
%EndExpansion
\frac{f_{i}(Y_{i},\boldsymbol{\theta})^{\alpha}}{\left(  \int f_{i}%
(y,\boldsymbol{\theta})^{\alpha+1}dy\right)  ^{\frac{\alpha}{\alpha+1}}}%
=\frac{1}{n}%
%TCIMACRO{\tsum \limits_{i=1}^{n}}%
%BeginExpansion
{\textstyle\sum\limits_{i=1}^{n}}
%EndExpansion
\frac{f_{i}(Y_{i},\boldsymbol{\theta})^{\alpha}}{L_{\alpha}^{i}\left(
\boldsymbol{\theta}\right)  }=\frac{1}{n}%
%TCIMACRO{\tsum \limits_{i=1}^{n}}%
%BeginExpansion
{\textstyle\sum\limits_{i=1}^{n}}
%EndExpansion
V_{i}(Y_{i},\boldsymbol{\theta}) \label{2.4}%
\end{equation}
and then the minimum RP estimator, $\widehat{\boldsymbol{\theta}}_{\alpha},$ for the common $\boldsymbol{\theta,}$ is given by
%\begin{equation}
%\widehat{\boldsymbol{\theta}}_{\alpha}=\left\{
%\begin{array}
%[c]{cc}%
%\alpha>0 & \arg\max_{\theta\in\Theta}H_{n}^{\alpha}(\boldsymbol{\theta})\\
%\alpha=0 & \arg\max_{\theta\in\Theta}H_{n}^{0}(\boldsymbol{\theta})
%\end{array}
%\right.  \label{2.4.1}%
%\end{equation}
\begin{equation}
	\widehat{\boldsymbol{\theta}}_{\alpha} =  \arg\max_{\theta\in\Theta}H_{n}^{\alpha}(\boldsymbol{\theta}),
	\label{2.4.1}
\end{equation}
with $H_{n}^{\alpha}(\boldsymbol{\theta})$ defined in (\ref{2.4}) for $\alpha >0$ and
$H_{n}^{0}(\boldsymbol{\theta})=\frac{1}{n}{\textstyle\sum\limits_{i=1}^{n}}\log f_{i}(Y_{i},\boldsymbol{\theta}).$

It is interesting to observe that when $Y_{1},...,Y_{n}$ are independent and identically distributed random variables, the estimator $\widehat{\boldsymbol{\theta}}_{\alpha}$ coincides with the estimator $\widehat{\boldsymbol{\theta}}_{\alpha}^{\ast}$ given in (\ref{1.2}). In the following section we shall establish the consistency of $\widehat{\boldsymbol{\theta}}_{\alpha}$, as well as its asymptotic distribution.

\section{Consistency and asymptotic distribution \label{sec:CONS}}

We shall assume in the following that the true densities $g_{i}$  $i=1,...,n,$
belong to the assumed model, i.e., $g_{i}\equiv f_{i}(y,\boldsymbol{\theta}),$
$i=1,...,n,$ for some common parameter $\boldsymbol{\theta}$. We denote by
$\boldsymbol{\theta}^*$ the true value of the unknown parameter. In the
following we shall denote \  $E_{\boldsymbol{\theta}^*}\left[  Y\right]  ={\textstyle\int} yf_{i}(y,\boldsymbol{\theta}^*)dy,$ and we introduce the matrices
\begin{equation}
\boldsymbol{\Psi}_{n}=\frac{1}{n}%
%TCIMACRO{\tsum \limits_{i=1}^{n}}%
%BeginExpansion
{\textstyle\sum\limits_{i=1}^{n}}
%EndExpansion
\boldsymbol{J}^{\left(  i\right)  },\label{2.4.2}%
\end{equation}
with
\[
\boldsymbol{J}^{\left(  i\right)  }=\left(  -E_{\boldsymbol{\theta}^*%
}\left[  \frac{\partial^{2}V_{i}(Y;\boldsymbol{\theta})}{\partial\theta
_{j}\partial\theta_{k}}\right]  \right)  _{i,j=1,...,p}%
\]
and
\begin{equation}
\boldsymbol{\Omega}_{n}=\frac{1}{n}%
%TCIMACRO{\tsum \limits_{i=1}^{n}}%
%BeginExpansion
{\textstyle\sum\limits_{i=1}^{n}}
%EndExpansion
Var_{\boldsymbol{\theta}^*}\left[  \left(  \frac{\partial V_{i}%
(Y;\boldsymbol{\theta})}{\partial\theta_{j}}\right)  _{j=1,..,p}\right]
.\label{2.4.3}%
\end{equation}

In order to get the asymptotic results we shall assume the following conditions:

\begin{description}
\item[C1.\label{itm:C1}] The support, $\mathcal{X},$ of the density functions $f_{i}%
(y,\boldsymbol{\theta})$ is the same for all $i$ and does not depend on
$\boldsymbol{\theta}.$

\item[C2.\label{itm:C2}] There exists an open subset $\Theta^{\ast}$ of $\Theta$ containing
the true value of the parameter $\boldsymbol{\theta}^*$ such that for almost
all $y\in\mathcal{X}$ the density $f_{i}(y,\boldsymbol{\theta})$ admits all
third derivatives with respect to $\boldsymbol{\theta}$ and $i=1,...,n.$

\item[C3.\label{itm:C3}] For $i=1,2,...$ the integrals
\[%
%TCIMACRO{\tint }%
%BeginExpansion
{\textstyle\int}
%EndExpansion
f_{i}(y,\boldsymbol{\theta})^{1+\alpha}dy\text{ }%
\]
can be differentiated thrice with respect to $\boldsymbol{\theta}$ and we can
interchange integration and differentiation.

\item[C4.\label{itm:C4}] For $i=1,2,...$ the matrices $\boldsymbol{J}^{\left(  i\right)  }$ are positive
definite. We denote by $\lambda_{n}$ the minimum eigenvalue of
$\boldsymbol{\Omega}_{n}$ and 
\[
\lambda_{0}=\inf_{n}\lambda_{n}>0.
\]

\item[C5.\label{itm:C5}] There exists functions $M_{jkl}^{\left(  i\right)  }$ such that
\[
\left\vert \frac{\partial^{3}V_{i}(y;\boldsymbol{\theta})}{\partial\theta
_{j}\partial\theta_{k}\partial\theta_{l}}\right\vert \leq M_{jkl}^{\left(
i\right)  }\left(  y\right)  ,\text{ \qquad}\forall\boldsymbol{\theta\in
\Theta^{\ast},}\text{ }\forall j,k,l
\]
and
\[
E_{\boldsymbol{\theta}^*}\left[  M_{jkl}^{\left(  i\right)  }\left(
Y\right)  \right]  =m_{jkl}<\infty,\text{ }\forall j,k,l.
\]

\item[C6.\label{itm:C6}] For all $j,k,l$ the sequences $\left\{  \frac{\partial
V_{i}(Y;\boldsymbol{\theta})}{\partial\theta_{j}}\right\}  _{j=1,...,p}$
,$\left\{  \frac{\partial^{2}V_{i}(Y;\boldsymbol{\theta})}{\partial\theta
_{j}\partial\theta_{k}}\right\}  _{j,k=1,..,p}$ and $\left\{  \frac
{\partial^{3}V_{i}(Y;\boldsymbol{\theta})}{\partial\theta_{j}\partial
\theta_{k}\partial l}\right\}  _{j,k,l=1,..,p}$ are uniformly integrable in
the Ces\'{a}ro sense, i.e.

\begin{align*}
\lim_{N\rightarrow\infty}\left(  \sup_{n>1}\frac{1}{n}%
%TCIMACRO{\tsum \limits_{i=1}^{n}}%
%BeginExpansion
{\textstyle\sum\limits_{i=1}^{n}}
%EndExpansion
E_{\boldsymbol{\theta}^*}\left[  \left\vert \frac{\partial V_{i}%
(Y;\boldsymbol{\theta})}{\partial\theta_{j}}\right\vert {\LARGE I}_{\left\{
\frac{\partial V_{i}(Y;\boldsymbol{\theta})}{\partial\theta_{j}}>N\right\}
}(Y)\right]  \right)   &  =0,\\
\lim_{N\rightarrow\infty}\left(  \sup_{n>1}\frac{1}{n}%
%TCIMACRO{\tsum \limits_{i=1}^{n}}%
%BeginExpansion
{\textstyle\sum\limits_{i=1}^{n}}
%EndExpansion
E_{\boldsymbol{\theta}^*}\left[  \left\vert \frac{\partial^{2}%
V_{i}(Y;\boldsymbol{\theta})}{\partial\theta_{j}\partial\theta_{k}}\right\vert
{\LARGE I}_{\left\{  \frac{\partial^{2}V_{i}(Y;\boldsymbol{\theta})}%
{\partial\theta_{j}\partial\theta_{k}}>N\right\}  }(Y)\right]  \right)   &
=0,\\
\lim_{N\rightarrow\infty}\left(  \sup_{n>1}\frac{1}{n}%
%TCIMACRO{\tsum \limits_{i=1}^{n}}%
%BeginExpansion
{\textstyle\sum\limits_{i=1}^{n}}
%EndExpansion
E_{\boldsymbol{\theta}^*}\left[  \left\vert \frac{\partial^{3}%
V_{i}(Y;\boldsymbol{\theta})}{\partial\theta_{j}\partial\theta_{k}%
\partial\theta_{l}}\right\vert {\LARGE I}_{\left\{  \frac{\partial^{2}%
V_{i}(Y;\boldsymbol{\theta})}{\partial\theta_{j}\partial\theta_{k}%
\partial\theta_{l}}>N\right\}  }(Y)\right]  \right)   &  =0.
\end{align*}

\item[C7.\label{itm:C7}] For all $\varepsilon>0$%
\[
\lim_{n\rightarrow\infty}\left\{  \frac{1}{n}%
%TCIMACRO{\tsum \limits_{i=1}^{n}}%
%BeginExpansion
{\textstyle\sum\limits_{i=1}^{n}}
%EndExpansion
E_{\boldsymbol{\theta}^*}\left[  \left\Vert \boldsymbol{\Omega}_{n}%
^{-\frac{1}{2}}\frac{\partial V_{i}(Y,\theta)}{\partial\boldsymbol{\theta}%
}\right\Vert _{2}^{2}{\LARGE I}_{\left\{  \left\Vert \boldsymbol{\Omega}%
_{n}^{-\frac{1}{2}}\frac{\partial V_{i}(Y,\theta)}{\partial\boldsymbol{\theta
}}\right\Vert _{2}^{2}\right\}  }(Y)\right]  >\varepsilon\sqrt{n}\right\}
=0.
\]

\end{description}

Note that \ref{itm:C6} gives sufficient conditions for the weak law of large numbers with i.n.i.d.o. (\cite{Chandra}), while  \ref{itm:C7} is the assumption required for the multivariate central limit theorem for i.n.i.d.o. (\cite{Feller}). The following theorem states the consistency of $\widehat{\boldsymbol{\theta}}_{\alpha}$ and the second one establishes its asymptotic distribution. Proof of both theorems are given in Appendix \ref{ap:con} and Appendix \ref{ap:asy}, respectively. 

\begin{theorem}
\label{Th1}Let $Y_{1},...,Y_{n}$ be i.n.i.d.o. each with a density function
$f_{i}(y,\boldsymbol{\theta}),$ $\boldsymbol{\theta}\in\Theta\subset
\mathbb{R}^{p}.$ If conditions \ref{itm:C1}--\ref{itm:C6} holds, there exists a consistent sequence $\widehat{\boldsymbol{\theta}}_{n}$ of the system of equations
\begin{equation}
\frac{\partial H_{n}^{\alpha}(\boldsymbol{\theta})}{\partial\boldsymbol{\theta
}}=\boldsymbol{0}_p. \label{2.5}%
\end{equation}

\end{theorem}

%\begin{proof}
%See Appendix A1.
%\end{proof}

\begin{theorem}
\label{Th2}Let $Y_{1},...,Y_{n}$ be i.n.i.d.o. each with a density function
$f_{i}(y,\boldsymbol{\theta}),$ $\boldsymbol{\theta}\in\Theta\subset
\mathbb{R}^{p}.$ If conditions \ref{itm:C1}--\ref{itm:C7} are satisfied the asymptotic distribution of the minimum RP estimator is given by
\begin{align}\label{eq:asym_est}
\sqrt{n}\boldsymbol{\Omega}_{n}^{-\frac{1}{2}}\boldsymbol{\Psi}_{n}\left(
\widehat{\boldsymbol{\theta}}_{\alpha}-\boldsymbol{\theta}^*\right)
\underset{n\rightarrow\infty}{\overset{L}{\rightarrow}}N(\boldsymbol{0}%
_{p},\boldsymbol{I}_{p}),
\end{align}
being $\boldsymbol{I}_{p}$ the p-dimensional identity matrix and the matrices
$\boldsymbol{\Psi}_{n}$ and $\boldsymbol{\Omega}_{n}$ were defined in
(\ref{2.4.2}) and (\ref{2.4.3}), respectively.
\end{theorem}

\section{Wald-type test for i.n.i.d.o. \label{sec:Wald}}

\subsection{Wald-type tests for simple null hypotheses}
We define a family of Wald-type test statistics based on minimum RP estimators for testing the hypothesis
\begin{align}\label{eq:simple_null}
H_0:\boldsymbol{\theta}=\boldsymbol{\theta}^0 \quad \text{against} \quad H_1:\boldsymbol{\theta}\neq \boldsymbol{\theta}^0,
\end{align}
for a given $\boldsymbol{\theta}^0 \in \Theta$.

\begin{definition}
 Let $\widehat{\boldsymbol{\theta}}_{\alpha}$ be the minimum RP estimator of $\boldsymbol{\theta}$. The family  of proposed Wald-type test statistics for testing the null hypothesis (\ref{eq:simple_null}) is given by
 \begin{align}\label{eq:Wald}
 \boldsymbol{W}_n^0(\boldsymbol{\theta}^0)=n(\widehat{\boldsymbol{\theta}}_{\alpha}-\boldsymbol{\theta}^0)^T\boldsymbol{\Sigma}_{\alpha}^{-1}(\boldsymbol{\theta}^0)(\widehat{\boldsymbol{\theta}}_{\alpha}-\boldsymbol{\theta}^0),
 \end{align}
 where  $\boldsymbol{\Sigma}_{\alpha}(\boldsymbol{\theta}^0)=\lim_{n\rightarrow\infty} \boldsymbol{\Psi}_n(\boldsymbol{\theta}^0) \boldsymbol{\Omega}_n^{-1}(\boldsymbol{\theta}^0)\boldsymbol{\Psi}_n(\boldsymbol{\theta}^0)$.
\end{definition}

\begin{theorem}
\label{th:asym_test} The asymptotic distribution of the Wald-type test
statistics $\boldsymbol{W}_n^0(\boldsymbol{\theta}^0)$, defined in (\ref{eq:Wald}), under the null hypothesis (\ref{eq:simple_null}), is a chi-square distribution with $p$ degrees of freedom.
\end{theorem}

Based on Theorem \ref{th:asym_test}, we shall reject the null hypothesis in (\ref{eq:simple_null})  if
\begin{equation}\label{eq:test}
 \boldsymbol{W}_n^0(\boldsymbol{\theta}^0)>\chi^2_{p,\nu},
\end{equation}
being $\chi^2_{p,\nu}$ the upper $\nu$-th quantile of $\chi^2_{p}$.

\begin{theorem}\label{th:power1}
Let $\boldsymbol{\theta}^*$ be the true value of $\boldsymbol{\theta}$, with $\boldsymbol{\theta}^*\neq \boldsymbol{\theta}^0$, and let us denote
\begin{align*}
\ell(\boldsymbol{\theta})=(\boldsymbol{\theta}-\boldsymbol{\theta}^0)^T\boldsymbol{\Sigma}_{\alpha}^{-1}(\boldsymbol{\theta}-\boldsymbol{\theta}^0)
\end{align*}
and \ $\sigma^2_{ \boldsymbol{W}_n^0(\boldsymbol{\theta}^0)}(\boldsymbol{\theta}^*)=4(\boldsymbol{\theta}^*-\boldsymbol{\theta}^0)^T\left[\boldsymbol{\Sigma}_{\alpha}^{-1}(\boldsymbol{\theta}^0)\boldsymbol{\Sigma}_{\alpha}(\boldsymbol{\theta}^*)\boldsymbol{\Sigma}_{\alpha}^{-1}(\boldsymbol{\theta}^0) \right](\boldsymbol{\theta}^*-\boldsymbol{\theta}^0).$
Then, 
\begin{align*}
\sqrt{n}\left(\ell(\widehat{\boldsymbol{\theta}}_{\alpha})-\ell(\boldsymbol{\theta}^*)\right)
\underset{n\rightarrow\infty}{\overset{L}{\rightarrow}}N(\boldsymbol{0},\sigma^2_{ \boldsymbol{W}_n^0(\boldsymbol{\theta}^0)}(\boldsymbol{\theta}^*)).
\end{align*}
\end{theorem}

\begin{corollary}
Theorem \ref{th:power1},  makes it possible to have an approximation of the power function for the test given in (\ref{eq:test}). This is given by
\begin{align*}
\pi^{\alpha}_{\boldsymbol{W}_n^0(\boldsymbol{\theta}^0)}(\boldsymbol{\theta}^*)&=P(\text{Rejecting } H_0|\boldsymbol{\theta}=\boldsymbol{\theta}^*)=P(\boldsymbol{W}_n^0(\boldsymbol{\theta}^0)>\chi^2_{p,\nu}|\boldsymbol{\theta}=\boldsymbol{\theta}^*)\\
&=P\left( \ell(\widehat{\boldsymbol{\theta}}_{\alpha})-\ell(\boldsymbol{\theta}^*)>\frac{\chi^2_{p,\nu}}{n}-\ell(\boldsymbol{\theta}^*)\right)\\
&=1-\Phi_n\left(\frac{\sqrt{n}}{\sigma_{ \boldsymbol{W}_n^0(\boldsymbol{\theta}^0)}(\boldsymbol{\theta}^*)}\left(\frac{\chi^2_{p,\nu}}{n}-\ell(\boldsymbol{\theta}^*) \right) \right),
\end{align*}
where $\Phi_n(\cdot)$  is a sequence of distribution functions which tends uniformly to the standard normal distribution function $\Phi(\cdot)$. We can observe that
\begin{align*}
\underset{n\rightarrow\infty}{\text{lim}}\pi^{\alpha}_{\boldsymbol{W}_n^0(\boldsymbol{\theta}^0)}(\boldsymbol{\theta}^*)=1 \quad \forall \alpha \geq 0,
\end{align*}
so the Wald-type tests are consistent in the sense of Fraser. 
\end{corollary}

From here, it can be deduced that the necessary sample size for the Wald-type tests to have a
predetermined power, $\pi^{\alpha}_{\boldsymbol{W}_n^0(\boldsymbol{\theta}^0)}(\boldsymbol{\theta}^*)\approx \pi^{*}$, is given by
\begin{align*}
n=\left[\frac{A+B+\sqrt{A(A+2B)}}{2\ell^2(\boldsymbol{\theta}^*)} \right],
\end{align*} 
where $A=\sigma^2_{ \boldsymbol{W}_n^0(\boldsymbol{\theta}^0)}(\boldsymbol{\theta}^*)(\Phi^{-1}(1-\pi^{\ast}))^{2}$, and $B=2\ell(\boldsymbol{\theta}^*)\chi^2_{p,\nu}$.
%\color{black}

\subsection{Wald-type tests  for composite null hypotheses}
We may be also interested on testing a set of $r<p$ redundant restrictions on the parameter vector $\boldsymbol{\theta}$. In this context, we are interested in testing
\begin{align}\label{eq:composite_null}
H_0:\boldsymbol{M}^T\boldsymbol{\theta}=\boldsymbol{m} \quad \text{against} \quad H_1:\boldsymbol{M}^T\boldsymbol{\theta}\neq\boldsymbol{m},
\end{align}
where $\boldsymbol{M}$ is $p\times r$  full rank matrix with $r <p$ and $\boldsymbol{m}$ is a $r$-vector.

\begin{definition}
Let $\widehat{\boldsymbol{\theta}}_{\alpha}$ be the minimum RP estimator of $\boldsymbol{\theta}$. The family  of proposed Wald-type test statistics for testing the null hypothesis (\ref{eq:composite_null}) is given by
 \begin{align}\label{eq:Wald_composite}
 \boldsymbol{W}_n(\widehat{\boldsymbol{\theta}}_{\alpha})=n(\boldsymbol{M}^T\widehat{\boldsymbol{\theta}}_{\alpha}-\boldsymbol{m})^T\left[\boldsymbol{M}^T\boldsymbol{\Sigma}_{\alpha}(\widehat{\boldsymbol{\theta}}_{\alpha})\boldsymbol{M}\right]^{-1}(\boldsymbol{M}^T\widehat{\boldsymbol{\theta}}_{\alpha}-\boldsymbol{m}).
 \end{align}
\end{definition}

\begin{theorem}
\label{th:asym_test_composite} The asymptotic distribution of the Wald-type test
statistics $\boldsymbol{W}_n(\widehat{\boldsymbol{\theta}}_{\alpha})$, defined in (\ref{eq:Wald_composite}), under the null hypothesis (\ref{eq:composite_null}), is a chi-square distribution with $r$ degrees of freedom.
\end{theorem}

Based on Theorem \ref{th:asym_test_composite}, we shall reject the null hypothesis in (\ref{eq:composite_null})  if
\begin{equation}\label{eq:test_composite}
\boldsymbol{W}_n(\widehat{\boldsymbol{\theta}}_{\alpha})>\chi^2_{r,\nu}.
\end{equation}

We could generalize our results to a more general restricted space $\Theta_0 \subset \Theta$ defined by a set of $r<p$ non-redundant restrictions of the form $h(\boldsymbol{\theta}) = \boldsymbol{0},$ by substituting the matrix $\boldsymbol{M}^T\boldsymbol{\Sigma}_{\alpha}(\widehat{\boldsymbol{\theta}}_{\alpha})\boldsymbol{M}$ by   $\boldsymbol{H}^T\boldsymbol{\Sigma}_{\alpha}(\widehat{\boldsymbol{\theta}}_{\alpha})\boldsymbol{H}$ with  $\boldsymbol{H} = \frac{\partial \boldsymbol{h}(\boldsymbol{\theta})}{\partial \boldsymbol{\theta}^T}$  in (\ref{eq:Wald_composite}). The asymptotic distribution stated in Theorem \ref{th:asym_test_composite} still holds.

\subsection{Contiguous alternatives hypothesis}
The previous results provides an asymptotic approximation to the power function of the proposed Wald-type tests. We now consider the particular set of contiguous alternatives hypothesis of the form 
\begin{equation}\label{eq:contiguous_hypothesis} H_1: \boldsymbol{\theta}_n = \boldsymbol{\theta}^0+n^{-1/2}\boldsymbol{d},
\end{equation}
where $\boldsymbol{d}$ is a fixed vector in $\mathbb{R}^p$ such that $\boldsymbol{d} \in \boldsymbol{\Theta}$ and $\boldsymbol{\theta}^0$ is an element of $\Theta_0.$ Note that the alternative hypothesis move towards $\boldsymbol{\theta}^0$ and it get closer when the sample size $n$ increases. 

\begin{theorem} \label{thm:contiguous}
	Under the contiguous alternative hypotheses given in (\ref{eq:contiguous_hypothesis}), the asymptotic distribution of the Wald-type test statistics, $\boldsymbol{W}_n(\widehat{\boldsymbol{\theta}}_{\alpha}),$ defined in (\ref{eq:composite_null}) is a non-central chi-square distribution with $r$ degrees of freedom and non-centrality parameter 
	$$\delta = \boldsymbol{d}^T \boldsymbol{M}[\boldsymbol{M}^T\boldsymbol{\Sigma}_{\alpha}(\widehat{\boldsymbol{\theta}}_{\alpha})\boldsymbol{M}]^{-1} \boldsymbol{M}^T \boldsymbol{d}.$$
\end{theorem}

\section{Influence function analysis \label{sec:IF}}

In this section we shall obtain the influence function (IF) of the minimum RP functional for the non-homogeneous case. We shall denote by $G_i$ the true distribution function associated to the observation $Y_i$ whose density function is denoted by $g_i$  and by $\boldsymbol{T}_{\alpha}(G_1,\dots,G_n)$ the minimum RP functional defined as the minimizer of 

\begin{align}\label{eq:IF1}
\sum_{i=1}^n H_n^{\alpha,(i)}=\frac{1}{n}\sum_{i=1}^n\left\{\frac{1}{1+\alpha}\log\left(  \int f_{i}(y,\boldsymbol{\theta})^{\alpha
+1}dy\right)  -\frac{1}{\alpha}\log\left(  \int f_{i}(y,\boldsymbol{\theta
})^{\alpha}dG_i(y)\right)  \right\}
\end{align}
or, fixed a value of $\alpha$, under appropriate differentiability conditions, as the solution of the system of equations obtained after differentiating (\ref{eq:IF1}) and equalling to zero

\begin{align}\label{eq:IF2}
\frac{1}{n} \sum_{i=1}^n \left\{\frac{\int f_i(y,\boldsymbol{\theta})^{\alpha+1}\boldsymbol{u}_i(y,\boldsymbol{\theta})g_i(y)dy}{\int f_i(y,\boldsymbol{\theta})^{\alpha+1}g_i(y)dy} - \frac{\int f_i(y,\boldsymbol{\theta})^{\alpha}\boldsymbol{u}_i(y,\boldsymbol{\theta})g_i(y)dy}{\int f_i(y,\boldsymbol{\theta})^{\alpha}g_i(y)dy} \right\}=\boldsymbol{0}_p.
\end{align}
By $\boldsymbol{u}_i(y,\boldsymbol{\theta})$ we are denoting
\begin{align*}
\boldsymbol{u}_i(y,\boldsymbol{\theta})=\frac{\partial \log (f_i(y,\boldsymbol{\theta}))}{\partial \boldsymbol{\theta}}
\end{align*}
and by 
\begin{align*}
g_{i,\varepsilon}=(1-\varepsilon)g_i+\varepsilon \Delta _{t_i}, \quad i=1,\dots,n
\end{align*}
the contaminated density where $\Delta _{t_i}$ is the degenerated distribution at point $t_i$. 

Let $\boldsymbol{\theta}=\boldsymbol{T}_{\alpha}(G_1,\dots,G_n)$ and we denote by
\begin{align*}
\boldsymbol{\theta}_{\varepsilon}^{i_0}=\boldsymbol{T}_{\alpha}(G_1,\dots, G_{i_0-1}, G_{i_0, \varepsilon}, G_{i_0+1}, \dots,G_n)
\end{align*}
the minimum Renyi pseudodistance functional with contamination only in the $i_0$-th direction, where $G_{i_0, \varepsilon}$ is the distribution function associated to the denisty function 
\begin{align*}
g_{i_0,\varepsilon}(y)=(1-\varepsilon)f_i(y,\boldsymbol{\theta})+\varepsilon \Delta _{i_0}(y)
\end{align*}
and
\begin{align*}
g_i(y)= \left\{ \begin{array}{ll}
         f_i(y,\boldsymbol{\theta}) & \mbox{if $i \neq i_0$};\\
        (1-\varepsilon) f_i(y,\boldsymbol{\theta})+\varepsilon \Delta _{i_0}  & \mbox{if $i=i_0$}.\end{array} \right.
\end{align*}

It is also possible to contaminate in all the directions and in this case we shall denote by
\begin{align*}
\boldsymbol{\theta}_{\varepsilon}=\boldsymbol{T}_{\alpha}(G_{1,\varepsilon},\dots G_{n,\varepsilon})
\end{align*}
the minimum RP functional with contamination in all directions. 

In the following theorem we present the expressions of the IF. See Appendix \ref{app:IF} for details.

\begin{theorem}\label{th:IF}
The influence function in the $i_0$-th direction is given by
\begin{align*}
IF(t_{i_0}, \boldsymbol{T}_{\alpha},G_1,\dots,G_n)=(\boldsymbol{M}_{n,\alpha}(\boldsymbol{\theta}))^{-1}\dfrac{-\boldsymbol{\ell}_{i_0,\alpha}(\boldsymbol{\theta})}{\left( \int f_{i_0}(y,\boldsymbol{\theta})^{\alpha} g_{i_0}(y)dy\right)^2},
\end{align*}
where
\begin{align*}
\boldsymbol{\ell}_{i_0,\alpha}(\boldsymbol{\theta})=f_{i_0}(y,\boldsymbol{\theta}) \int f_{i_0}(y,\boldsymbol{\theta})^{\alpha+1} \boldsymbol{u}_{i_0}(y,\boldsymbol{\theta})dy -f_{i_0}(y,\boldsymbol{\theta})\boldsymbol{u}_{i_0}(y,\boldsymbol{\theta}) \int f_{i_0}(y,\boldsymbol{\theta})^{\alpha+1} dy,
\end{align*}
\begin{align}\label{eq:MIF}
\boldsymbol{M}_{n,\alpha}(\boldsymbol{\theta})=\frac{1}{n}{\sum\limits_{i=1}^{n}}\left[\dfrac{\boldsymbol{A}_{i,\alpha}(\boldsymbol{\theta})}{\left( \int f_{i}(y,\boldsymbol{\theta})^{1+\alpha} dy\right)^2}-\dfrac{\boldsymbol{A}^*_{i,\alpha}(\boldsymbol{\theta})}{\left( \int f_{i}(y,\boldsymbol{\theta})^{\alpha} g_{i}(y)dy\right)^2}\right],
\end{align}
and 
\begin{small}
\begin{align*}
\boldsymbol{A}_{i,\alpha}(\boldsymbol{\theta})=&\left[(1+\alpha)  \int f_{i}(y,\boldsymbol{\theta})^{\alpha+1} \boldsymbol{u}^T_{i}(y,\boldsymbol{\theta})\boldsymbol{u}_{i}(y,\boldsymbol{\theta})dy + \int f_{i}(y,\boldsymbol{\theta})^{\alpha+1} \frac{\partial \boldsymbol{u}_{i}(y,\boldsymbol{\theta})}{\partial \boldsymbol{\theta}}dy\right]\int f_{i}(y,\boldsymbol{\theta})^{\alpha+1} dy \\
&- (1+\alpha) \left(\int f_{i}(y,\boldsymbol{\theta})^{\alpha+1} \boldsymbol{u}_{i}(y,\boldsymbol{\theta}) dy \right)\left(\int f_{i}(y,\boldsymbol{\theta})^{\alpha+1} \boldsymbol{u}_{i}(y,\boldsymbol{\theta}) dy \right)^T,\\
\boldsymbol{A}^*_{i,\alpha}(\boldsymbol{\theta})=&\left[\alpha  \int f_{i}(y,\boldsymbol{\theta})^{\alpha}g_i(y) \boldsymbol{u}^T_{i}(y,\boldsymbol{\theta})\boldsymbol{u}_{i}(y,\boldsymbol{\theta})dy + \int f_{i}(y,\boldsymbol{\theta})^{\alpha}g_i(y) \frac{\partial \boldsymbol{u}_{i}(y,\boldsymbol{\theta})}{\partial \boldsymbol{\theta}}dy\right]\int f_{i}(y,\boldsymbol{\theta})^{\alpha}g_i(y) dy \\
&- \alpha \left(\int f_{i}(y,\boldsymbol{\theta})^{\alpha} g_i(y)\boldsymbol{u}_{i}(y,\boldsymbol{\theta}) dy \right)\left(\int f_{i}(y,\boldsymbol{\theta})^{\alpha}g_i(y) \boldsymbol{u}_{i}(y,\boldsymbol{\theta}) dy \right)^T.
\end{align*}
\end{small}
Similarly, the influence function in all directions is given by
\begin{align*}
IF(t_1,\dots,t_n, \boldsymbol{T}_{\alpha},G_1,\dots,G_n)=(\boldsymbol{M}_{n,\alpha}(\boldsymbol{\theta}))^{-1} \sum\limits_{i=1}^n\dfrac{-\boldsymbol{\ell}_{i,\alpha}(\boldsymbol{\theta})}{\left( \int f_{i}(y,\boldsymbol{\theta})^{\alpha} g_{i}(y)dy\right)^2}.
\end{align*}
\end{theorem}

\begin{remark}\label{remark:IFs}
When  the true distribution $g_i$ belongs to the model so that $g_i(y)=f_i(y,\boldsymbol{\theta})$ for $i=1,\dots,n$, then $\boldsymbol{M}_{n,\alpha}(\boldsymbol{\theta})$  given in (\ref{eq:MIF}) coincides with $\boldsymbol{\Psi}_n(\boldsymbol{\theta})$ given in (\ref{2.4.2}) and  \begin{align*}
IF(t_{i_0}, \boldsymbol{T}_{\alpha},F_{1,\boldsymbol{\theta}},\dots,F_{n,\boldsymbol{\theta}})&=\boldsymbol{\Psi}_n^{-1}(\boldsymbol{\theta})\boldsymbol{D}_{i_0,\alpha}(\boldsymbol{\theta}),\\
IF(t_1,\dots,t_n, \boldsymbol{T}_{\alpha},F_{1,\boldsymbol{\theta}},\dots,F_{n,\boldsymbol{\theta}})&=\boldsymbol{\Psi}_n^{-1}(\boldsymbol{\theta})\sum_{i=1}^n\boldsymbol{D}_{i,\alpha}(\boldsymbol{\theta}),
\end{align*}
with $\boldsymbol{D}_{i,\alpha}(\boldsymbol{\theta})=\dfrac{-\boldsymbol{\ell}_{i,\alpha}(\boldsymbol{\theta})}{\left( \int f_{i}(y,\boldsymbol{\theta})^{\alpha+1} dy\right)^2}$.
\end{remark}

\begin{remark}
In particular, letting $t_i=t$, $G_i=G$, $f_{i}(y,\boldsymbol{\theta})=f(y,\boldsymbol{\theta})$ for $i=1,\dots,n$  (this situation corresponds to the case of independent and identically distributed, i.i.d., random variables) and $g(y)=f(y,\boldsymbol{\theta})$, we have
\begin{align*}
IF(t, \boldsymbol{T}_{\alpha},G)=(\boldsymbol{M}_{\alpha}(\boldsymbol{\theta}))^{-1}\left[f(y,\boldsymbol{\theta})^{\alpha}\boldsymbol{u}(y,\boldsymbol{\theta})-\boldsymbol{c}_{\alpha}(\boldsymbol{\theta})f(y,\boldsymbol{\theta})^{\alpha} \right],
\end{align*}
where
\begin{align*}
\boldsymbol{c}_{\alpha}(\boldsymbol{\theta})=\frac{\int f(y,\boldsymbol{\theta})^{\alpha+1} \boldsymbol{u}(y,\boldsymbol{\theta})dy}{\int f(y,\boldsymbol{\theta})^{\alpha+1} dy}
\end{align*}
and 
\begin{align*}
\boldsymbol{M}_{\alpha}(\boldsymbol{\theta})=\frac{ 1}{\int f(y,\boldsymbol{\theta})^{\alpha+1}dy} \left[ \int f(y,\boldsymbol{\theta})^{\alpha+1}dy \int f(y,\boldsymbol{\theta})^{\alpha+1} \boldsymbol{u}(y,\boldsymbol{\theta})\boldsymbol{u}^T(y,\boldsymbol{\theta})  dy \right. \\
\left. - \left( \int f(y,\boldsymbol{\theta})^{\alpha+1} \boldsymbol{u}(y,\boldsymbol{\theta}) dy \right)\left( \int f(y,\boldsymbol{\theta})^{\alpha+1} \boldsymbol{u}(y,\boldsymbol{\theta}) dy \right)^T \right] ,
\end{align*}
as in \cite{Broniatowski}.
\end{remark}

\subsection{Influence function of the Wald-type test statistics}
Once we have computed the IF for the minimum RP estimators, we can define and study the IF for the Wald-type test statistics. First, we  define the associated statistical functional, evaluated at $G_1,\dots,G_n$ as
\begin{align}
W^0_{\alpha}(G_1,\dots,G_n)=(\boldsymbol{T}_{\alpha}(G_1,\dots,G_n)-\boldsymbol{\theta}^0)^T\boldsymbol{\Sigma}_{\alpha}^{-1}(\boldsymbol{\theta}^0)(\boldsymbol{T}_{\alpha}(G_1,\dots,G_n)-\boldsymbol{\theta}^0),
\end{align}
corresponding to (\ref{eq:simple_null}) for the simple null hypothesis. Let us consider first the  contamination only in one direction, say $i_0$-th direction. The corresponding IF is then defined as
\begin{align} \label{eq:IF_wald}
IF(t_{i_0},W^0_{\alpha},G_1,\dots,G_n)&=\frac{\partial}{\partial \varepsilon} \left. W^0_{\alpha}(G_1,\dots,G_{i_0-1},G_{i_0,\varepsilon},G_{i_0+1},\dots,G_n)\right|_{\varepsilon=0}\\
&=2\boldsymbol{T}_{\alpha}(G_1,\dots,G_n)-\boldsymbol{\theta}^0)^T\boldsymbol{\Sigma}_{\alpha}^{-1}(\boldsymbol{\theta}^0) IF(t_{i_0},\boldsymbol{T}_{\alpha},G_1,\dots,G_n). \notag
\end{align}
However, if we evaluate (\ref{eq:IF_wald}) at the null distribution $G_i=F_{i,\boldsymbol{\theta}^0}$, it becomes identically zero. Therefore, it becomes necessary to consider the second order IF of the proposed Wald-type test functional
\begin{align*}
IF^{(2)}(t_{i_0},W^0_{\alpha},G_1,\dots,G_n)&=\frac{\partial^2}{\partial^2 \varepsilon} \left. W^0_{\alpha}(G_1,\dots,G_{i_0-1},G_{i_0,\varepsilon},G_{i_0+1},\dots,G_n)\right|_{\varepsilon=0}\\
&=2IF(t_{i_0},\boldsymbol{T}_{\alpha},G_1,\dots,G_n)\boldsymbol{\Sigma}_{\alpha}^{-1}(\boldsymbol{\theta}^0) IF(t_{i_0},\boldsymbol{T}_{\alpha},G_1,\dots,G_n). \notag
\end{align*}
Similarly, we can consider contamination in all directions, obtaining that the second order influence function of the proposed Wald-type tests functional for testing simple null hypothesis is given by
\begin{align*}
&IF^{(2)}(t_1,\dots,t_n,W^0_{\alpha},G_1,\dots,G_n)\\
& \quad =\frac{\partial^2}{\partial^2 \varepsilon} \left. W^0_{\alpha}(G_1,\dots,G_{i_0-1},G_{i_0,\varepsilon},G_{i_0+1},\dots,G_n)\right|_{\varepsilon=0} \notag \\
&\quad  =2IF(t_1,\dots,t_n,\boldsymbol{T}_{\alpha},G_1,\dots,G_n)\boldsymbol{\Sigma}_{\alpha}^{-1}(\boldsymbol{\theta}^0) IF(t_1,\dots,t_n,\boldsymbol{T}_{\alpha},G_1,\dots,G_n). \notag
\end{align*}

\begin{remark}
When  the true distribution  belongs to the model, then the second order influence functions of the proposed Wald-type tests functional for testing simple null hypothesis in (\ref{eq:simple_null}) are given by
 \begin{align*}
&IF^{(2)}(t_{i_0}, W_{\alpha},F_{1,\boldsymbol{\theta}^0},\dots,F_{n,\boldsymbol{\theta}^0})\\
& \quad \quad \quad =2 \left[\boldsymbol{D}_{i_0,\alpha}(\boldsymbol{\theta}^0)\right]^T\left[\boldsymbol{\Psi}_n^{-1}(\boldsymbol{\theta}^0)\boldsymbol{\Sigma}_{\alpha}^{-1}(\boldsymbol{\theta}^0)\boldsymbol{\Psi}_n^{-1}(\boldsymbol{\theta}^0)\right]\left[\boldsymbol{D}_{i_0,\alpha}(\boldsymbol{\theta}^0)\right],\\
&IF^{(2)}(t_{1},\dots.t_n, W_{\alpha},F_{1,\boldsymbol{\theta}^0},\dots,F_{n,\boldsymbol{\theta}^0})\\
&\quad \quad \quad =2 \left[\sum_{i=1}^n\boldsymbol{D}_{i,\alpha}(\boldsymbol{\theta}^0)\right]^T\left[\boldsymbol{\Psi}_n^{-1}(\boldsymbol{\theta}^0)\boldsymbol{\Sigma}_{\alpha}^{-1}(\boldsymbol{\theta}^0)\boldsymbol{\Psi}_n^{-1}(\boldsymbol{\theta}^0)\right]\left[\sum_{i=1}^n\boldsymbol{D}_{i,\alpha}(\boldsymbol{\theta}^0)\right].
\end{align*}
\end{remark}
\begin{remark}
In a similar manner, when the true distribution  belongs to the model, the second order influence functions of the proposed Wald-type tests functionals for testing composite null hypothesis in (\ref{eq:composite_null}) are given by
 \begin{align*}
&IF^{(2)}(t_{i_0}, W^0_{\alpha},F_{1,\boldsymbol{\theta}^0},\dots,F_{n,\boldsymbol{\theta}^0})\\
& \quad \quad =2 \left[\boldsymbol{\Psi}_n^{-1}(\boldsymbol{\theta}^0)\boldsymbol{D}_{i_0,\alpha}(\boldsymbol{\theta}^0)\right]^T\boldsymbol{M}\left[\boldsymbol{M}^T\boldsymbol{\Sigma}_{\alpha}(\boldsymbol{\theta}^0)\boldsymbol{M}\right]^{-1}\boldsymbol{M}^T\left[\boldsymbol{\Psi}_n^{-1}(\boldsymbol{\theta}^0)\boldsymbol{D}_{i_0,\alpha}(\boldsymbol{\theta}^0)\right],\\
&IF^{(2)}(t_{1},\dots.t_n, W^0_{\alpha},F_{1,\boldsymbol{\theta}^0},\dots,F_{n,\boldsymbol{\theta}^0})\\
& \quad \quad =2 \left[\boldsymbol{\Psi}_n^{-1}(\boldsymbol{\theta}^0)\sum_{i=1}^n\boldsymbol{D}_{i_0,\alpha}(\boldsymbol{\theta}^0)\right]^T\boldsymbol{M}\left[\boldsymbol{M}^T\boldsymbol{\Sigma}_{\alpha}(\boldsymbol{\theta}^0)\boldsymbol{M}\right]^{-1}\boldsymbol{M}^T\left[\boldsymbol{\Psi}_n^{-1}(\boldsymbol{\theta}^0)\sum_{i=1}^n\boldsymbol{D}_{i_0,\alpha}(\boldsymbol{\theta}^0)\right].
\end{align*}
\end{remark}

\section{Multiple linear regression model} \label{sec:MRM}

Consider the MLRM
\begin{equation}\label{eq:linear}
Y_{i}=\boldsymbol{X}_{i}^{T}\boldsymbol{\beta}+\varepsilon_{i}, \quad i=1,\dots,n,
\end{equation}
where the errors $\varepsilon_{i}^{\prime }s$ are i.i.d. normal random variables with mean zero and variance $\sigma ^{2}$, $\boldsymbol{X}_{i}^{T}=(X_{i1},...,X_{ip})$  is the vector of independent variables corresponding to the $i$-th condition and $\boldsymbol{\beta} = \left( \beta_{1},...,\beta_{p}\right) ^{T}$ is the vector of regression coefficients to be estimated.  We will consider that, for each $i$, $\boldsymbol{X}_{i}$ is fixed, yielding to independent but not identically distributed  $Y_{i}'s$ (i.n.i.d.o.), with $Y_i\sim \mathcal{N}(\boldsymbol{X}_{i}^{T}\boldsymbol{\beta},\sigma^2)$.  Under the previous notation, with $f_i\equiv \mathcal{N}(\boldsymbol{X}_{i}^{T}\boldsymbol{\beta},\sigma^2)$, we have, for $\alpha>0$,
\begin{align*}
V_i(Y_{i};\boldsymbol{\theta},\boldsymbol{X}_{i})=& \frac{\frac{1}{(2\pi)^{\alpha/2}\sigma^{\alpha}}\exp\left(\frac{-\alpha (Y_{i}-\boldsymbol{X}_{i}^{T}\boldsymbol{\beta})^2}{2\sigma^2}\right)}{\left((2\pi)^{\alpha/2}\sigma^{\alpha}\sqrt{1+\alpha}\right)^{-\frac{\alpha}{\alpha+1}}}\\
=& \left(\frac{1+\alpha}{2\pi}\right)^{\frac{\alpha}{2(\alpha+1)}}\sigma ^{-\frac{\alpha }{\alpha +1}}\exp \left( -\frac{\alpha }{2}\left( \frac{Y_{i}-\boldsymbol{X}_{i}^{T}\boldsymbol{\beta }}{\sigma }\right) ^{2}\right).
\end{align*}
Thus, our objective function to be minimized becomes
\begin{align*}
\frac{1}{n}\sum_{i=1}^nV_i(Y_{i};\boldsymbol{\theta},\boldsymbol{X}_{i})=& \left(\frac{1+\alpha}{2\pi}\right)^{\frac{\alpha}{2(\alpha+1)}}\frac{1}{n}\sum_{i=1}^n\sigma ^{-\frac{\alpha }{\alpha +1}}\exp \left( -\frac{\alpha }{2}\left( \frac{Y_{i}-\boldsymbol{X}_{i}^{T}\boldsymbol{\beta }}{\sigma }\right) ^{2}\right).
\end{align*}
Taking into account that the term $\left(\frac{1+\alpha}{2\pi}\right)^{\frac{\alpha}{2(\alpha+1)}}$ does not depend on the model parameters, we have that, for $\alpha>0$
\begin{align}
\left( \widehat{\boldsymbol{\beta }}_{\alpha },\widehat{\sigma }_{\alpha
}\right) = \arg \max_{\boldsymbol{\beta },\sigma }{\textstyle \sum\limits_{i=1}^{n}}\sigma ^{-\frac{\alpha }{\alpha +1}}\exp \left( -\frac{\alpha }{2}\left( \frac{Y_{i}-\boldsymbol{X}_{i}^{T}\boldsymbol{\beta }}{ \sigma }\right) ^{2}\right).  
\end{align}
Derivating with respect to $\boldsymbol{\beta }$ and $\sigma$  we see that the estimators $\widehat{\boldsymbol{\beta }}_{\alpha}$ and $\widehat{\sigma}_{\alpha}$ are solutions of the system 
\begin{equation}\label{eq:estimating}
%\left\{
\begin{array}{l}
{\textstyle\sum\limits_{i=1}^{n}}\exp \left( -\frac{\alpha }{2}\left( \frac{Y_{i}-\boldsymbol{X}_{i}^{T}\boldsymbol{\beta }}{\sigma }\right) ^{2}\right)
\left( \frac{Y_{i}-\boldsymbol{X}_{i}^{T}\boldsymbol{\beta }}{\sigma }\right) \boldsymbol{X}_{i}=\boldsymbol{0}_{p} \\
{\textstyle\sum\limits_{i=1}^{n}}\exp \left( -\frac{\alpha }{2}\left( \frac{Y_{i}-\boldsymbol{X}_{i}^{T}\boldsymbol{\beta }}{\sigma }\right) ^{2}\right)\left\{ \left( \frac{Y_{i}-\boldsymbol{X}_{i}^{T}\boldsymbol{\beta }}{\sigma}\right) ^{2}-\frac{1}{1+\alpha }\right\} =0
\end{array},
%\right. . 
\end{equation}
which is exactly the same as the one suggested by Castilla et al. (2020) for the case of homogeneous data. If $\alpha=0$, we have

\begin{equation}
\left( \widehat{\boldsymbol{\beta }}_{\alpha=0 },\widehat{\sigma }_{\alpha=0
}\right) =\arg \max_{\boldsymbol{\beta },\sigma }\frac{1}{(2\pi \sigma
^{2})^{n/2}}\exp \left( -\frac{\left\Vert \boldsymbol{Y}-\mathbb{X}%
\boldsymbol{\beta }\right\Vert _{2}^{2}}{2\sigma ^{2}}\right)
\end{equation}
and we get the system necessary  to get the MLE of $\boldsymbol{\beta}$ and $\sigma$, whose
well-known solution is given by
\begin{equation*}
\widehat{\boldsymbol{\beta }}_{0}=(\mathbb{X}^{T}\mathbb{X}\mathbf{)}%
^{-1}\mathbb{X}^{T}\mathbf{Y}\text{ \ and \ }\widehat{\sigma }_{0}^{2}=\frac{1}{n}{\textstyle\sum\limits_{i=1}^{n}}\left( Y_{i}-\boldsymbol{%
X}_{i}^{T}\widehat{\boldsymbol{\beta }}_{0}\right) ^{2}, 
\end{equation*}
where $\mathbb{X}^{T}=(\boldsymbol{X}_1,\dots,\boldsymbol{X}_n)_{p\times n}$ is the matrix of explanatory variables.

\begin{lemma}\label{lemma:MLRM}
Consider the set-up of the MLRM with i.n.i.d.o. defined in (\ref{eq:linear}) and assume that the true data generating density belongs to the model family. If the following mild conditions about the explanatory variables hold
\begin{description}
\item[M1.\label{itm:M1}] The values of the explanatory variables are such that, for all $j$, $k$ and $l$
\begin{align*}
&\underset{n>1\  1\leq i\leq n}{\text{sup max } }  |X_{ij}|=O(1), \quad \text{and} \quad \underset{n>1\  1\leq i\leq n}{\text{sup max } } |X_{ij}X_{ij}|=O(1)
\end{align*}
and
$$
\frac{1}{n}\sum_{i=1}^n|x_{ij}x_{ik}x_{il}|=O(1).
$$
\item[M2.\label{itm:M2}] The matrix $\mathbb{X}^{T}$ satisfies
\begin{align*}
&\underset{n}{\text{inf }} \left[\text{min eigenvalue of } \frac{1}{n} \mathbb{X}^{T}\mathbb{X}\right]>0,\\
&n \times \underset{1\leq i\leq n}{\text{max }} \left[\boldsymbol{X}_i^T (\mathbb{X}^{T}\mathbb{X})^{-1}\boldsymbol{X}_i \right]=O(1),
\end{align*}
\end{description}
then \ref{itm:C1}--\ref{itm:C7} are satisfied.
\end{lemma}

On the other hand, after some heavy computations we follow that  expressions (\ref{2.4.2}) and (\ref{2.4.3}) are given by 
\begin{equation}
   \boldsymbol{\Psi}_{n}= \begin{bmatrix}
    \frac{-1}{\sigma^2(\alpha+1)^{3/2}}\left(\frac{1}{n}\sum_{i=1}^n\boldsymbol{X}_i\boldsymbol{X}^T_i\right) & 0 \\
    0 &  \frac{-2}{\sigma^2(\alpha+1)^{5/2}}
    \end{bmatrix},
\end{equation}
and
 \begin{equation}
    \boldsymbol{\Omega}_{n}= \begin{bmatrix}
    \frac{\left(\frac{1}{n}\sum_{i=1}^n\boldsymbol{X}_i\boldsymbol{X}^T_i\right)}{(2\alpha+1)^{3/2}} & 0 \\
    0 &  \frac{(3\alpha^2+4\alpha+2)}{\sigma^2(\alpha+1)^2(2\alpha+1)^{5/2}}
    \end{bmatrix}.
\end{equation}

\begin{theorem} \label{thm:asymptotic}
Consider the set-up of the MLRM with i.n.i.d.o. defined in (\ref{eq:linear}) and assume that the true data generating density belongs to the model family and  the observed explanatory variables satisfy conditions \ref{itm:M1} and \ref{itm:M2}. Then, 
\begin{enumerate}
\item There exists a consistent sequence as $\widehat{\boldsymbol{\theta}}_{\alpha}=(\widehat{\boldsymbol{\beta}}_{\alpha},\widehat{\sigma}_{\alpha})$ of solutions to the minimum R\'enyi estimating equations (\ref{eq:estimating}).
\item  $\widehat{\boldsymbol{\beta}}_{\alpha}$ and $\widehat{\sigma}_{\alpha}$ are asymptotically independent   and  their asymptotic  distribution is given by
\begin{align*}
\sqrt{n}\left(\widehat{\boldsymbol{\theta}}_{\alpha}-\boldsymbol{\theta}^*\right)
\underset{n\rightarrow\infty}{\overset{L}{\rightarrow}}N(\boldsymbol{0}%
,\boldsymbol{\Sigma}_{\alpha}),
\end{align*}
with 
\begin{align*}
&\boldsymbol{\Sigma}_{\alpha}=\lim_{n\rightarrow \infty}\boldsymbol{\Sigma}_{n}, \\
    &\boldsymbol{\Sigma}_{n} = \begin{bmatrix}
    \frac{\sigma^2(\alpha+1)^3}{(2\alpha+1)^{3/2}}\left(\frac{1}{n}\sum_{i=1}^n\boldsymbol{X}_i\boldsymbol{X}^T_i\right)^{-1} & \boldsymbol{0} \\
    \boldsymbol{0} &  \frac{\sigma^2(\alpha+1)^3(3\alpha^2+4\alpha+2)}{4(2\alpha+1)^{5/2}}
    \end{bmatrix}.
\end{align*}

\end{enumerate}
\end{theorem}

We could now apply the theory stated to test any simple or composite hypothesis on the linear regression parameters. The asymptotic distribution under the null hypothesis of the Wald type tests defined in (\ref{eq:Wald_composite}) is given in Theorem \ref{thm:asymptotic} and the asymptotic distribution of the Wald-type test statistics under contiguous alternative hypothesis is given in Theorem \ref{thm:contiguous}. The non-centrality parameter in Theorem \ref{thm:contiguous} can be expressed as 
	$$\delta = \boldsymbol{d^\ast}^T[\boldsymbol{M}^T\boldsymbol{\Sigma}_n \boldsymbol{M}]^{-1} \boldsymbol{d^\ast},$$
with $\boldsymbol{d}^\ast = \boldsymbol{M}^T \boldsymbol{d}.$ 
If we consider the composite null hypothesis (\ref{eq:composite_null}), then
$$\delta = \frac{(2\alpha+1)^{3/2}}{\sigma(\alpha+1)^3} \boldsymbol{d^\ast}^T[\boldsymbol{M}^T \left(\frac{1}{n}\sum_{i=1}^n\boldsymbol{X}_i\boldsymbol{X}^T_i\right) \boldsymbol{M}]^{-1} \boldsymbol{d^\ast}.$$

\begin{figure}
	\center
	\begin{tabular}{cc}
	\includegraphics[scale =0.4]{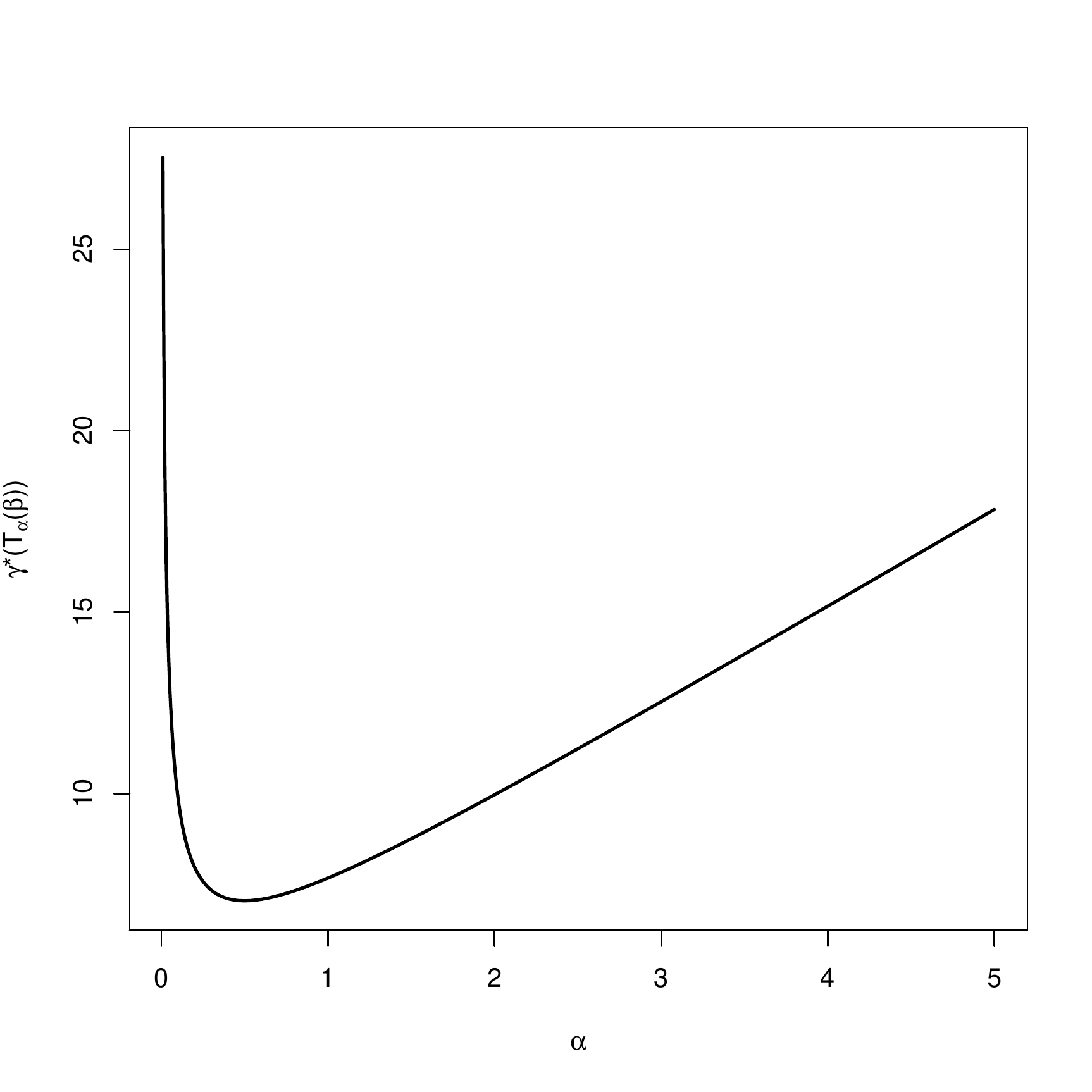} & \includegraphics[scale =0.4]{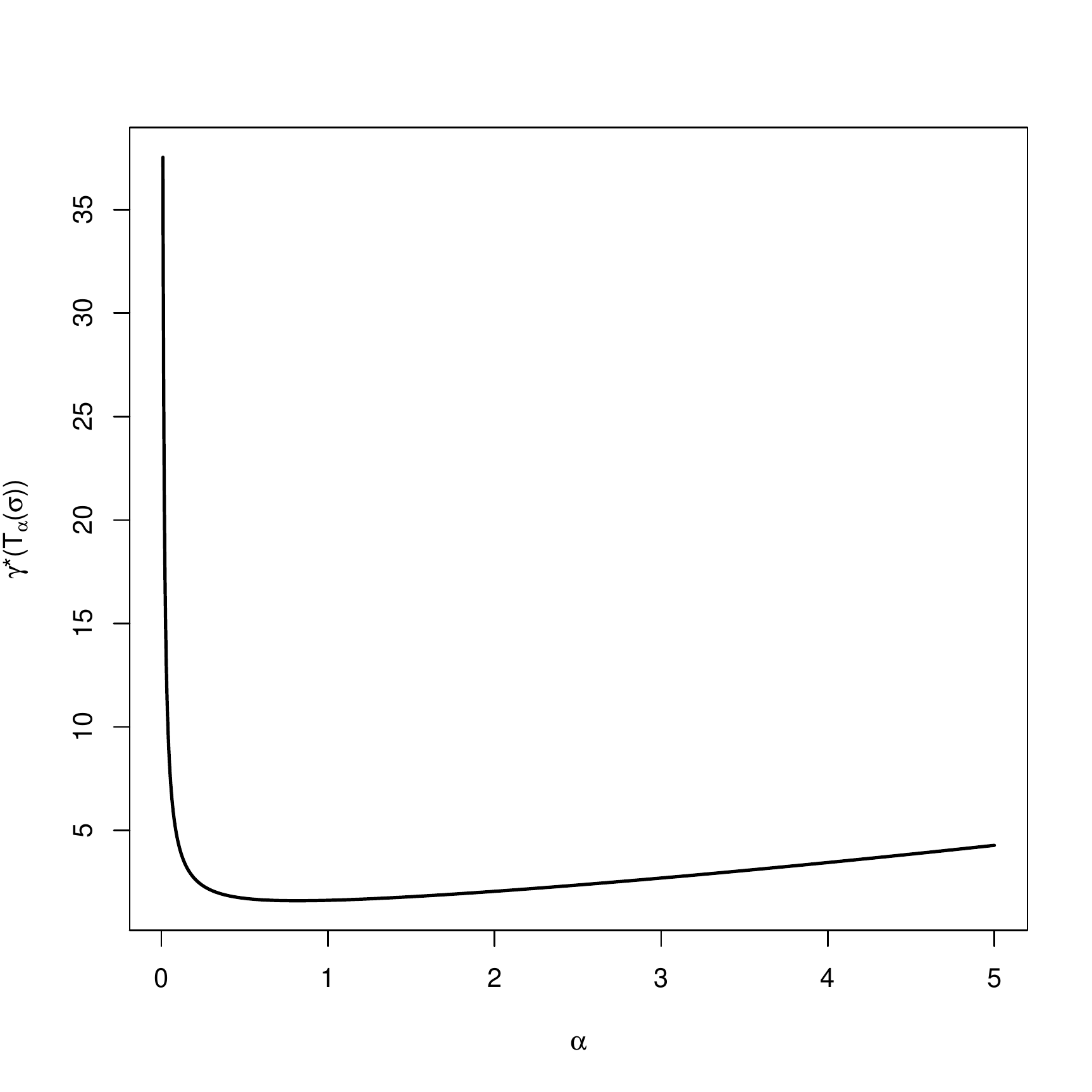}\\
\end{tabular}
	\caption{Gross error function for $\boldsymbol{T}_{\alpha}(\boldsymbol{\beta})$ (left) and $\boldsymbol{T}_{\alpha}(\sigma)$ (right).\label{figure:gross_error}}
\end{figure}

\begin{figure}
	\center
	\begin{tabular}{cc}
		\includegraphics[scale=0.42]{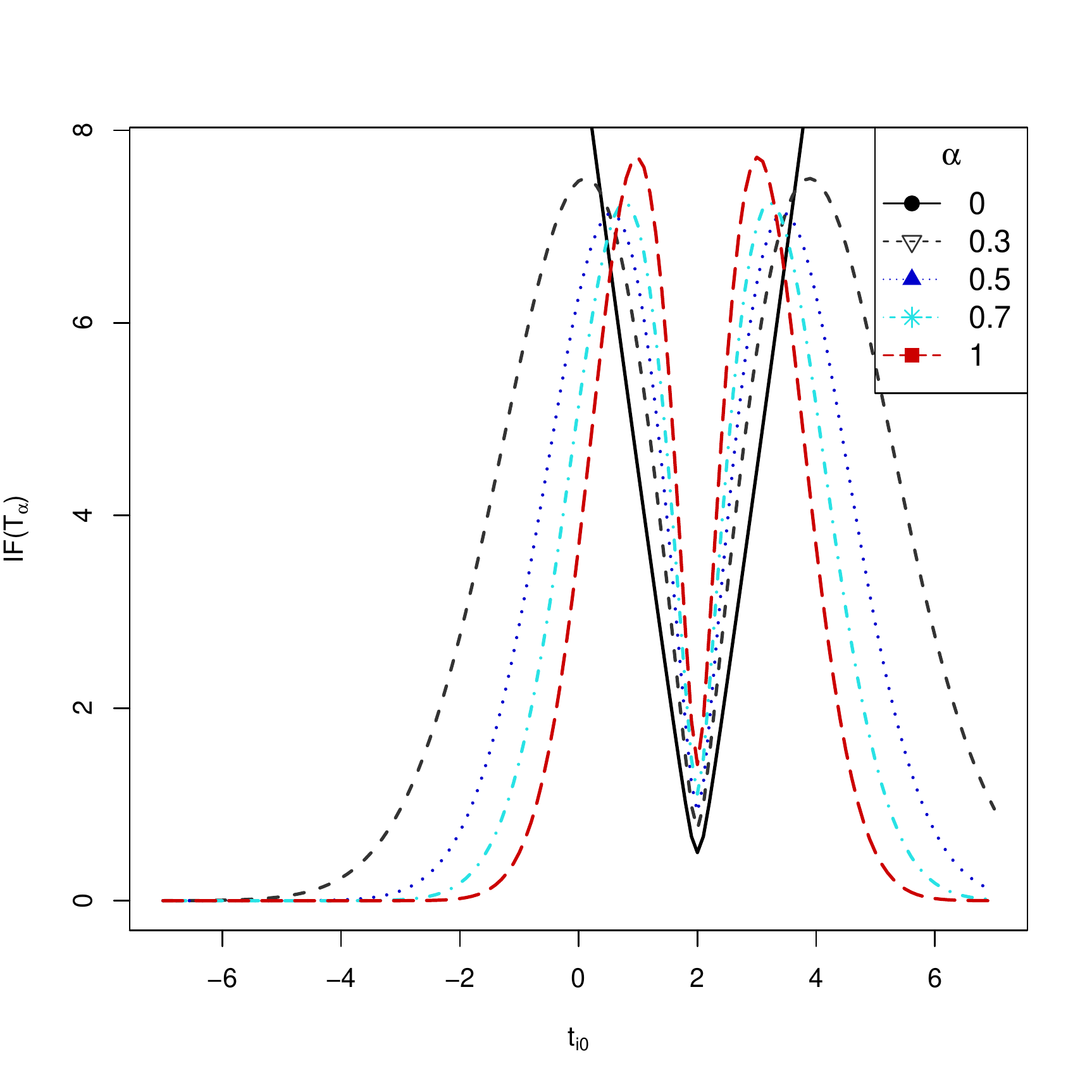}& \includegraphics[scale=0.42]{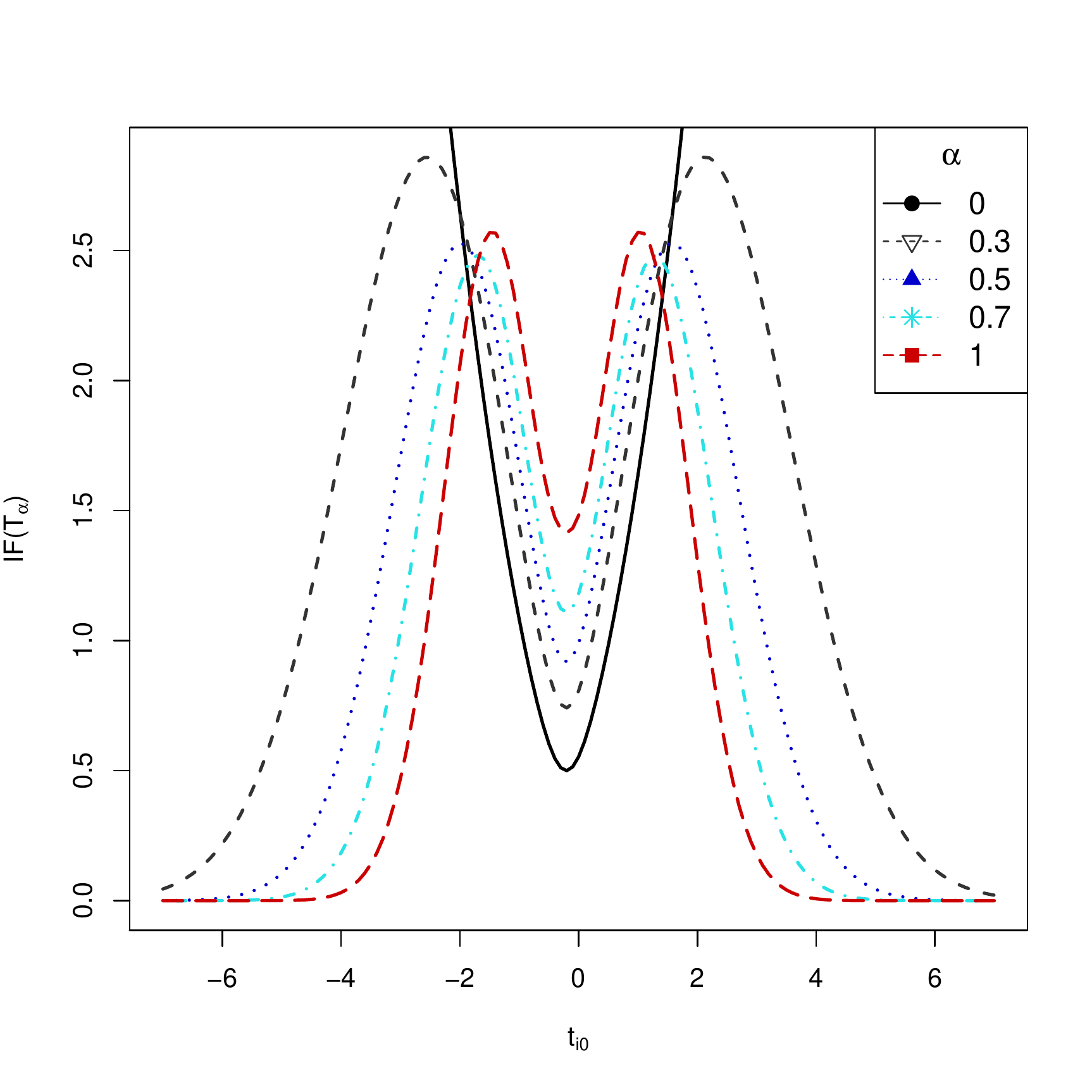}\\
		\includegraphics[scale=0.42]{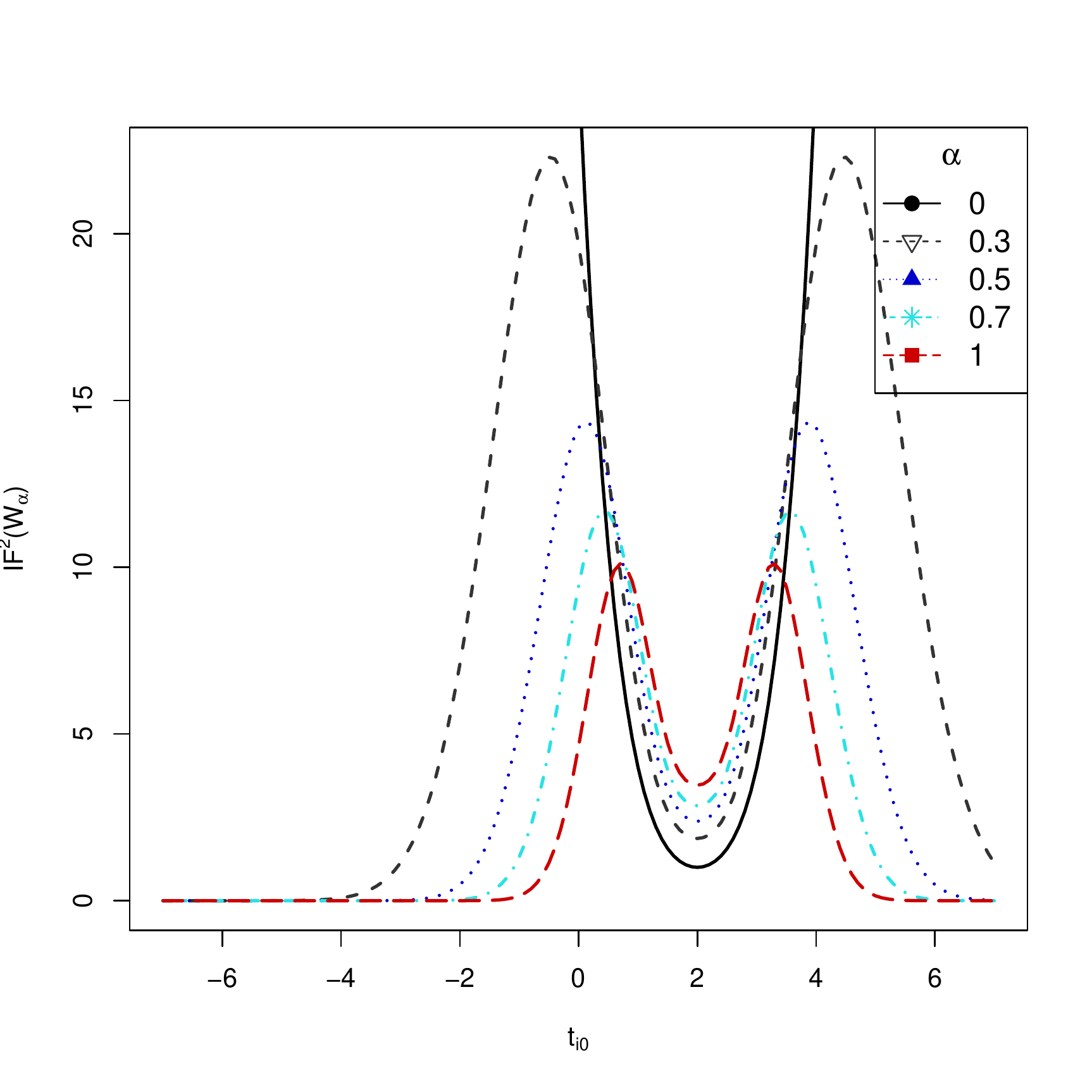}& \includegraphics[scale=0.42]{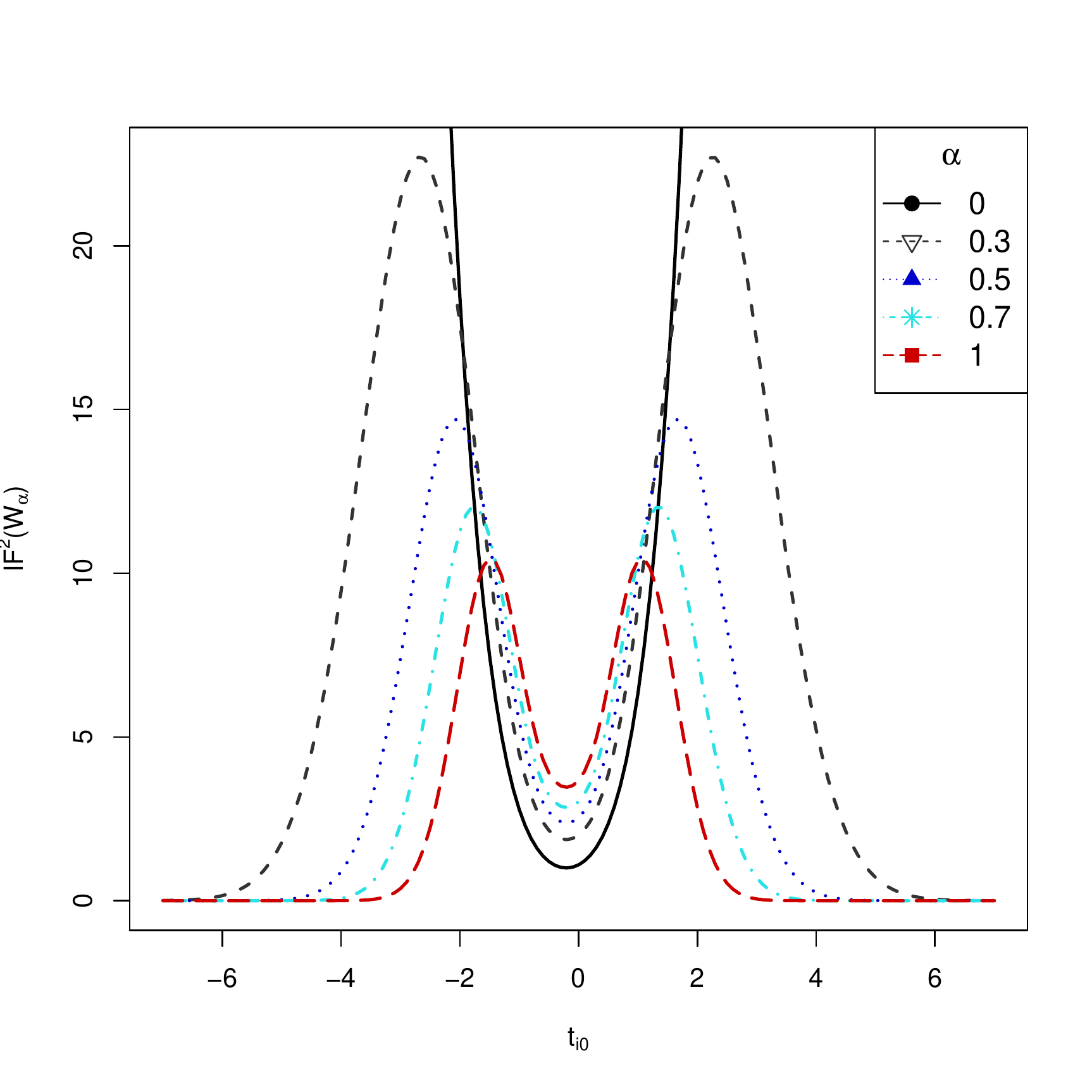}
	\end{tabular}
	\caption{$\ell_2$-norm of the first order IF of the minimum RP estimator (top) and second order IF of Wald type test estimators for testing (\ref{eq:simple_null}) with $\boldsymbol{\theta}_0 = (1,1,1)^T$ (bottom) with fixed \ref{itm:Design1} (left) and \ref{itm:Design2} (right), and contamination in the direction $i_0 = 1.$ \label{figure:IF}}
\end{figure}

Now, based on Remark \ref{remark:IFs} we can obtain the IF of the functional associated to the minimum RP estimator of $\boldsymbol{\theta}$. These are given by
\begin{align*}
IF(t_{i_0}, \boldsymbol{T}_{\alpha},F_{1,\boldsymbol{\theta}},\dots,F_{n,\boldsymbol{\theta}})&=\boldsymbol{\Psi}_n^{-1}(\boldsymbol{\theta})\left(\boldsymbol{D}^T_{i_0,\alpha}(\boldsymbol{\beta}),\boldsymbol{D}_{i_0,\alpha}(\sigma)\right)^T,\\
IF^{(2)}(t_{i_0}, W_{\alpha},F_{1,\boldsymbol{\theta}},\dots,F_{n,\boldsymbol{\theta}})
&  =2 \left(\boldsymbol{D}^T_{i_0,\alpha}(\boldsymbol{\beta}),\boldsymbol{D}_{i_0,\alpha}(\sigma)\right)\left[\boldsymbol{\Psi}_n^{-1}(\boldsymbol{\theta})\boldsymbol{\Sigma}_{n}^{-1}(\boldsymbol{\theta})\boldsymbol{\Psi}_n^{-1}(\boldsymbol{\theta})\right]\\
& \quad \times \left(\boldsymbol{D}^T_{i_0,\alpha}(\boldsymbol{\beta}),\boldsymbol{D}_{i_0,\alpha}(\sigma)\right)^T,\\
IF^{(2)}(t_{i_0}, W^0_{\alpha},F_{1,\boldsymbol{\theta}},\dots,F_{n,\boldsymbol{\theta}})
& =2 \left[\boldsymbol{\Psi}_n^{-1}(\boldsymbol{\theta})\left(\boldsymbol{D}^T_{i_0,\alpha}(\boldsymbol{\beta}),\boldsymbol{D}_{i_0,\alpha}(\sigma)\right)^T\right]^T\boldsymbol{M}\left[\boldsymbol{M}^T\boldsymbol{\Sigma}_{n}(\boldsymbol{\theta})\boldsymbol{M}\right]^{-1}\\
& \quad \times \boldsymbol{M}^T\left[\boldsymbol{\Psi}_n^{-1}(\boldsymbol{\theta})\left(\boldsymbol{D}^T_{i_0,\alpha}(\boldsymbol{\beta}),\boldsymbol{D}_{i_0,\alpha}(\sigma)\right)\right],
\end{align*}
with 
\begin{align*}
\boldsymbol{D}_{i_0,\alpha}(\boldsymbol{\beta})&=\frac{-1}{\sigma}\exp\left( -\frac{\alpha}{2}\left( \frac{t_{i_0}-\boldsymbol{x}_{i_0}^{T}\boldsymbol{\beta}}{\sigma}\right) ^{2}\right)\left( \frac{t_{i_0}-\boldsymbol{x}_{i_0}^{T}\boldsymbol{\beta}}{\sigma}\right) \boldsymbol{x}_{i_0},\\
\boldsymbol{D}_{i_0,\alpha}(\sigma)&=\frac{-1}{\sigma^2}\exp\left( -\frac{\alpha}{2}\left( \frac {t_{i_0}-\boldsymbol{x}_{i_0}^{T}\boldsymbol{\beta}}{\sigma}\right) ^{2}\right) \left[
\left( \frac{t_{i_0}-\boldsymbol{x}_{i_0}^{T}\boldsymbol{\beta}}{\sigma }\right)
^{2}-\frac{1}{\alpha+1}\right]. 
\end{align*}
%%% IMPORTANTE CREO QUE SOBRA LA ESPERANZA E[XX^T] PORQUE ADEMÁS LOS NUESTROS SON FIJOS
Note that, since $\boldsymbol{\Psi}_n(\boldsymbol{\theta})$ and $\boldsymbol{\Sigma}_{n}(\boldsymbol{\theta})$ are diagonal matrices, we could express separately the IF of the functional $\boldsymbol{T}_{\alpha}(\boldsymbol{\beta})$ and $\boldsymbol{T}_{\alpha}(\sigma)$  associated to the minimum RP estimator, $\widehat{\boldsymbol{\beta}}_\alpha$ and $\widehat{\sigma}_\alpha$ respectively. Following \cite{Basu2018}, we consider two different fixed design matrices for the univariate lineal regression model:
\begin{description}[leftmargin=50pt]%ojo! paquete enumitem (subido al dropbox)
	\item[Design 1\label{itm:Design1}] Two-points design.We fix $\boldsymbol{x}_i = (1,x_{i1})^T$, with $x_{i1} = a, i=1,..,n/2$ and $x_{i1} = b, i=n/2+1,..,n.$
	\item[Design 2\label{itm:Design2}] Fixed-Normal design. We fix $\boldsymbol{x}_i = (1,x_{i1})^T$, where $x_{i1}, i=1,...,n$ are prefixed independent and identically distributed observations from a $\mathcal{N}(\mu_x=0,\sigma_x=1).$ 
\end{description}
Figure \ref{figure:IF} presents the $\ell_2$-norm of the first order IF of the minimum RP estimator and second order IF of Wald type test estimators for testing (\ref{eq:simple_null}) with $\boldsymbol{\theta}_0 = (1,1,1)^T$ with both fixed designs and contamination in one direction for different values of $\alpha$. Clearly, the IF is bounded for positives values of the parameter $\alpha$ and is unbounded at the MLE, highlighting it lack of robustness.
Moreover, the supremum of the $\ell_2$-norm of the IF indicates the robustness of the estimator. Hence, we could study the optimal parameter of $\alpha$ trough the gross error sensitivity function. We define the gross error sensitivity of the functional $\boldsymbol{T}_{\alpha}$ considering contamination in the $i_0$ direction as 
\begin{equation}
\gamma^\ast(\boldsymbol{T}_{\alpha}, F_{1,\boldsymbol{\theta}},\dots,F_{n,\boldsymbol{\theta}}) = \operatorname{sup}_{t_{i0}}\left\{|| IF_{i_0} (t_{i_0}, \boldsymbol{T}_{\alpha},F_{1,\boldsymbol{\theta}},\dots,F_{n,\boldsymbol{\theta}})\right || \}.
\end{equation}
Considering separately the influence function of the functionals $\boldsymbol{T}_{\alpha}(\boldsymbol{\beta})$ and $\boldsymbol{T}_{\alpha}(\sigma),$ it is easy to show that 
\begin{equation}
\begin{aligned}
\gamma^\ast(\boldsymbol{T}_{\alpha}(\boldsymbol{\beta}), F_{1,\boldsymbol{\theta}},\dots,F_{n,\boldsymbol{\theta}}) &= \sigma  \frac{(\alpha+1)^{3/2}}{\alpha^{1/2}} \exp\left(-\frac{1}{2}\right) \bigg|\bigg|\left(\frac{1}{n}\sum_{i=1}^n\boldsymbol{X}_i\boldsymbol{X}^T_i\right)^{-1}\boldsymbol{x}_{i_0} \bigg|\bigg|,\\
\gamma^\ast(\boldsymbol{T}_{\alpha}(\sigma), F_{1,\boldsymbol{\theta}},\dots,F_{n,\boldsymbol{\theta}}) &=  \frac{(\alpha+1)^{5/2}}{\alpha} \exp\left( -\frac{3\alpha+2}{2(\alpha+1)}\right).
\end{aligned}
\end{equation}
Figure \ref{figure:gross_error} represents the gross error sensitivity functions depending on the parameter $\alpha,$ using  \ref{itm:Design1} and fixing the true standard error $\sigma =1$. The optimal value of $\alpha$ depends on the functional, being $\alpha = 1/2$ and $\alpha = \sqrt{2/3}$ for $\boldsymbol{T}_{\alpha}(\boldsymbol{\beta})$ and $\boldsymbol{T}_{\alpha}(\sigma)$ respectively. Therefore, a global optimal value of $\alpha,$ in terms of robustness, should varies between values $\alpha = 0.5$ to $\alpha =0.82$ if the true standard error is $\sigma =1.$

Finally, we study the Asymptotic Relative Efficiency (ARE) of the proposed minimum RP estimators with respect to the MLE, which is B.A.N. (Best Asymptotically Normal). The ARE of is computed as the ratio of their asymptotic variances. Note that this ratio does not depend on the regression parameters, but is only determined by $\alpha.$
\begin{equation}
\begin{aligned}
\text{ARE}(\widehat{\beta}_\alpha) &= \frac{(2\alpha+1)^{3/2}}{(\alpha+1)^3},\\
\text{ARE}(\widehat{\sigma}_\alpha) &= \frac{2(2\alpha+1)^{5/2}}{(\alpha+1)^3(3\alpha^2+4\alpha+2)}.
\end{aligned}
\end{equation}
Table \ref{table:ARE} represents the ARE of the minimum RP estimator, $(\widehat{\beta}_\alpha,\widehat{\sigma}_\alpha)$. As shown, the increment of $\alpha$ leads to an efficiency loss, which is heightened for the standard error estimator. Therefore, to ensure  sufficing efficiency, the parameter $\alpha$ should be chosen from low values. However, the efficiency reduction might worth in contrast with the robustness advantage. In view of the error sensitivity function study, values above $\alpha=0.82$ are not advocated.

\begin{table}[ht]
	\center
	\caption{ARE of the minimum RP estimator with respect to the MLE for the multiple linear regression model for different values of $\alpha$.}
	\label{table:ARE}
	\begin{tabular}{rrrrrrrrrr}
		\hline
		$\alpha$ & 0 & 0.1 & 0.2 & 0.3 & 0.4 & 0.5 & 0.8 & 1 & 1.5 \\ 
		\hline
		$ \text{ARE}(\widehat{\beta}_\alpha)$ \footnotesize{$(\times 100)$}  & 100.00 & 98.76 & 95.86 & 92.12 & 88.01 & 83.81 & 71.89 & 64.95 & 51.20 1 \\
		$\text{ARE}(\widehat{\sigma}_\alpha)$  \footnotesize{$(\times 100)$}  & 100.00 & 97.54 & 91.92 & 84.95 & 77.65 & 70.57 & 52.50 & 43.30 & 27.77  \\ 
		\hline
	\end{tabular}
\end{table}

\section{Numerical results \label{sec:simulation}}
%\subsection{Simulation Study}
We empirically evaluate the performance of the proposed Wald-type test statistics based on minimum RP estimator for MLRM through an extensive simulation study. We consider the univariate regression model with fixed design matrix
$$ y_i= \beta_0 + \beta_1 x_{1i} + \varepsilon_i, \hspace{0.3cm} i=1,..,n$$
and the two different design matrices presented in Section \ref{sec:MRM}. We generate the response variable from the linear regression model (\ref{eq:linear}) with regression parameters $\boldsymbol{\beta}^0=(1,1)$ and $\sigma^0 = 1$. To introduce contamination on the data, we swap the true regression vector to $\boldsymbol{\beta}^0=(1.5,2)$ for a $10\%$ of the sample size.
We analyse the performance of Wald-type tests for simple null hypothesis on both regression parameters, $\boldsymbol{\beta}$ and $\sigma$, at different values of the tuning parameter $\alpha$. Note that, for the proposed design matrices, the matrix $\frac{1}{n}\sum_{i=1}^n\boldsymbol{X}_i\boldsymbol{X}^T_i$ is finitely defined and is positive definite.

\begin{figure}
\begin{tabular}{cc}
	\includegraphics[scale=0.4]{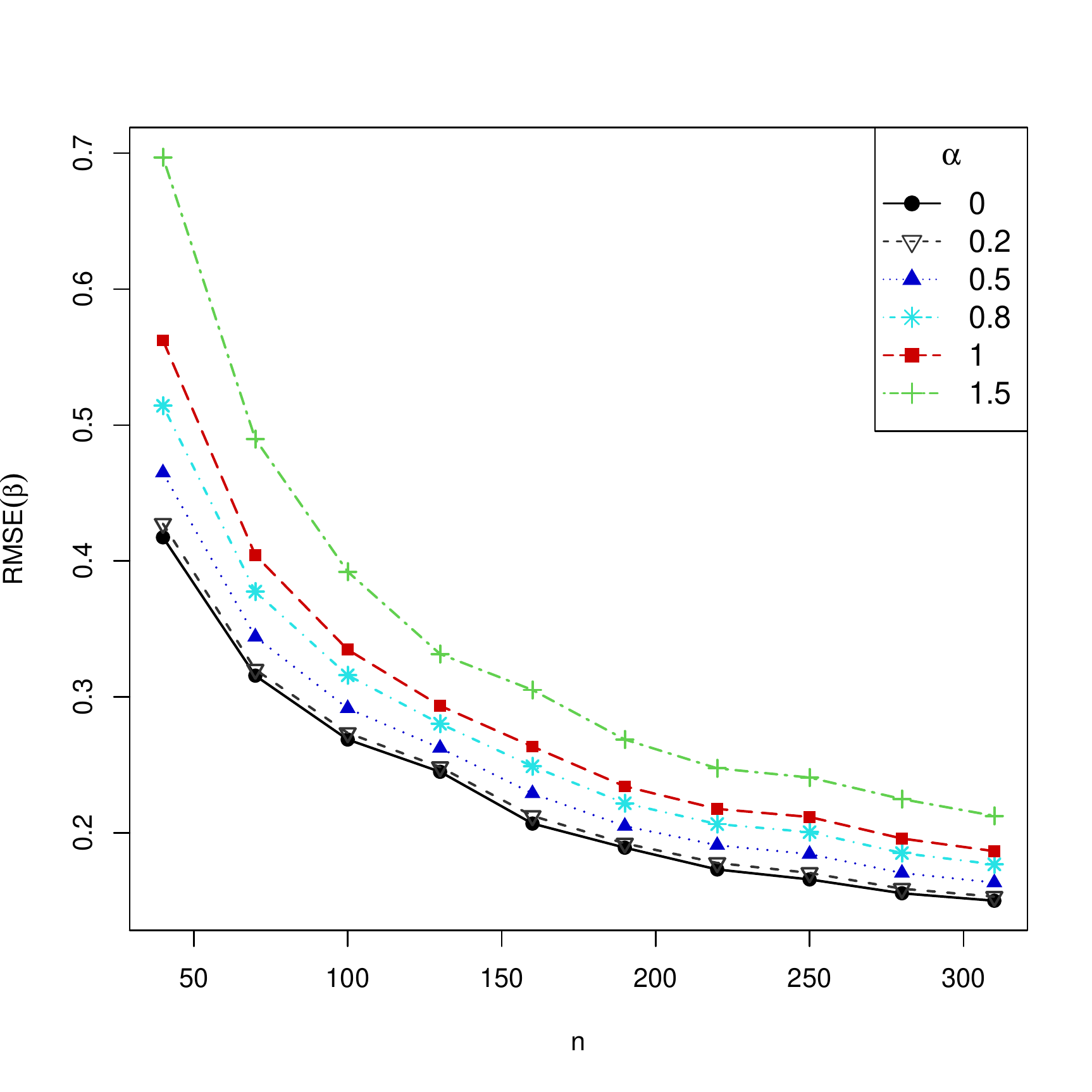}&
	\includegraphics[scale=0.4]{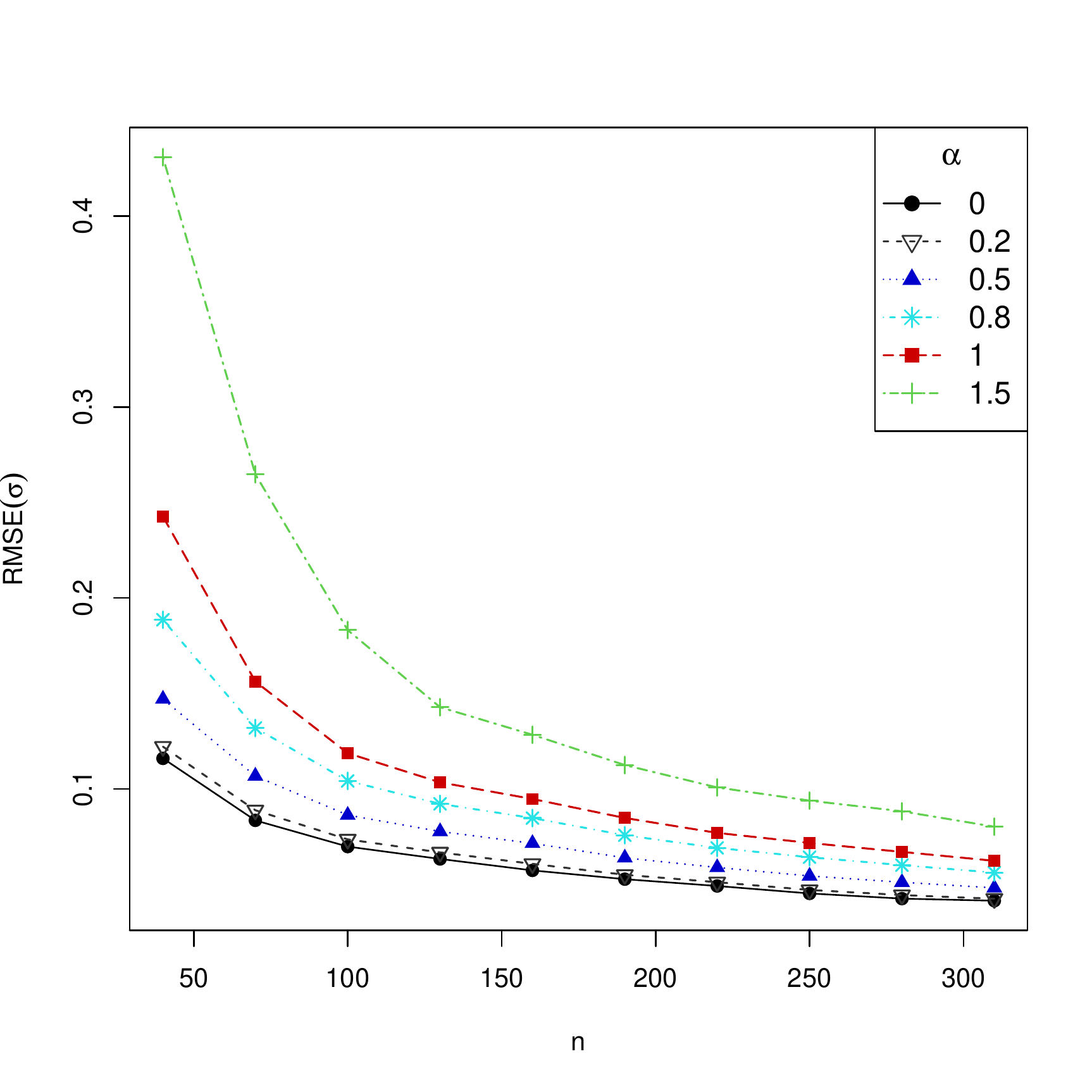}\\
	\includegraphics[scale=0.4]{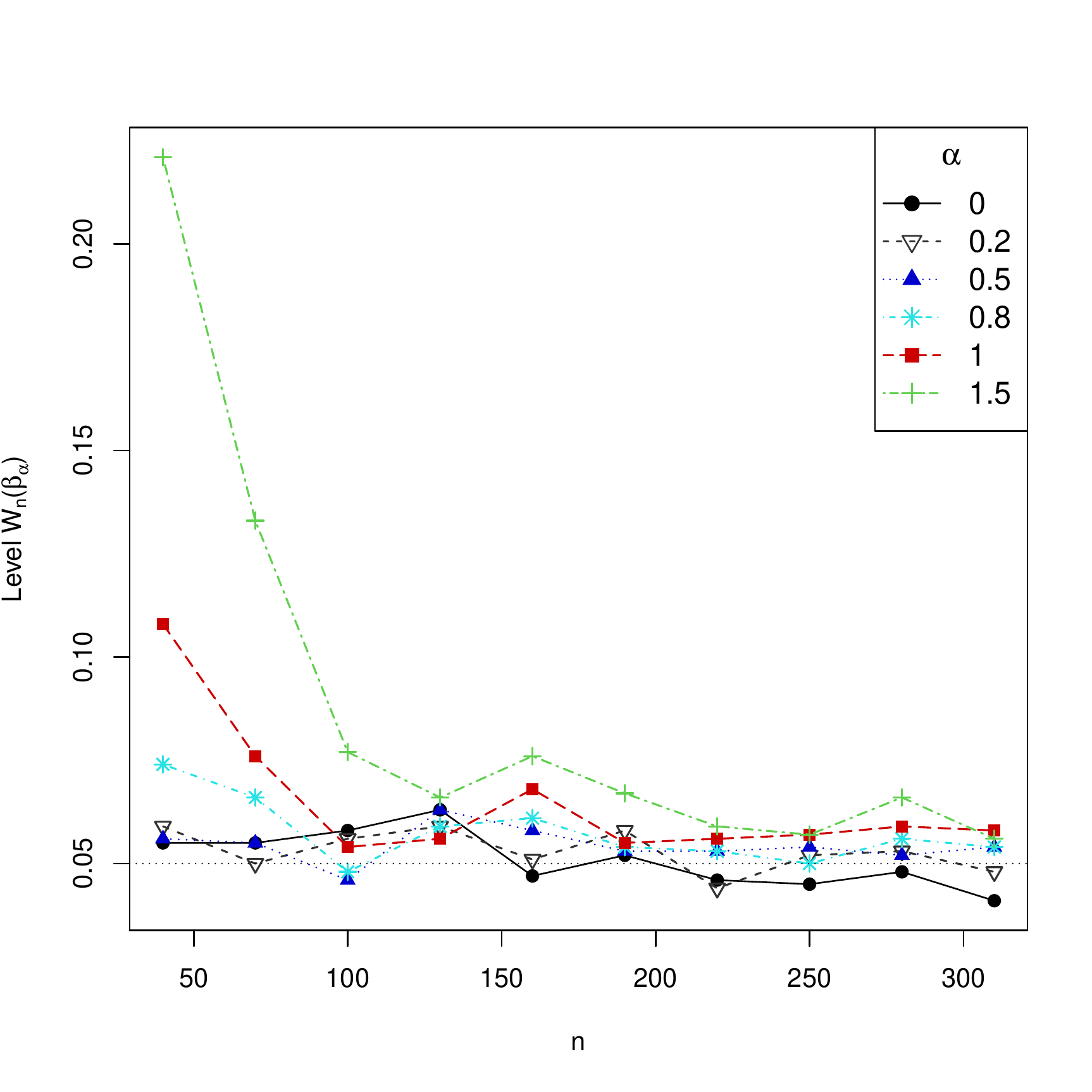}&
	\includegraphics[scale=0.4]{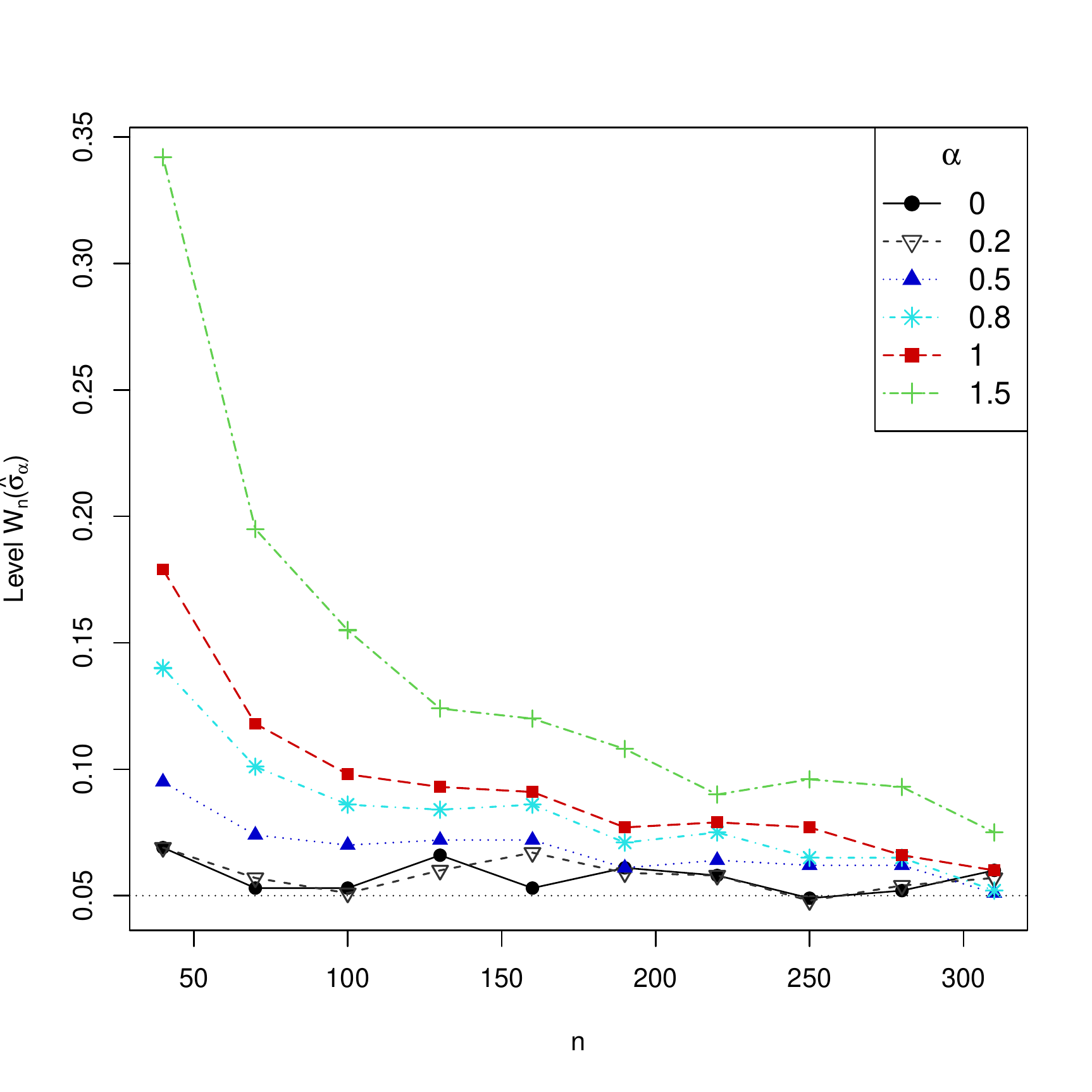}\\
	\includegraphics[scale=0.4]{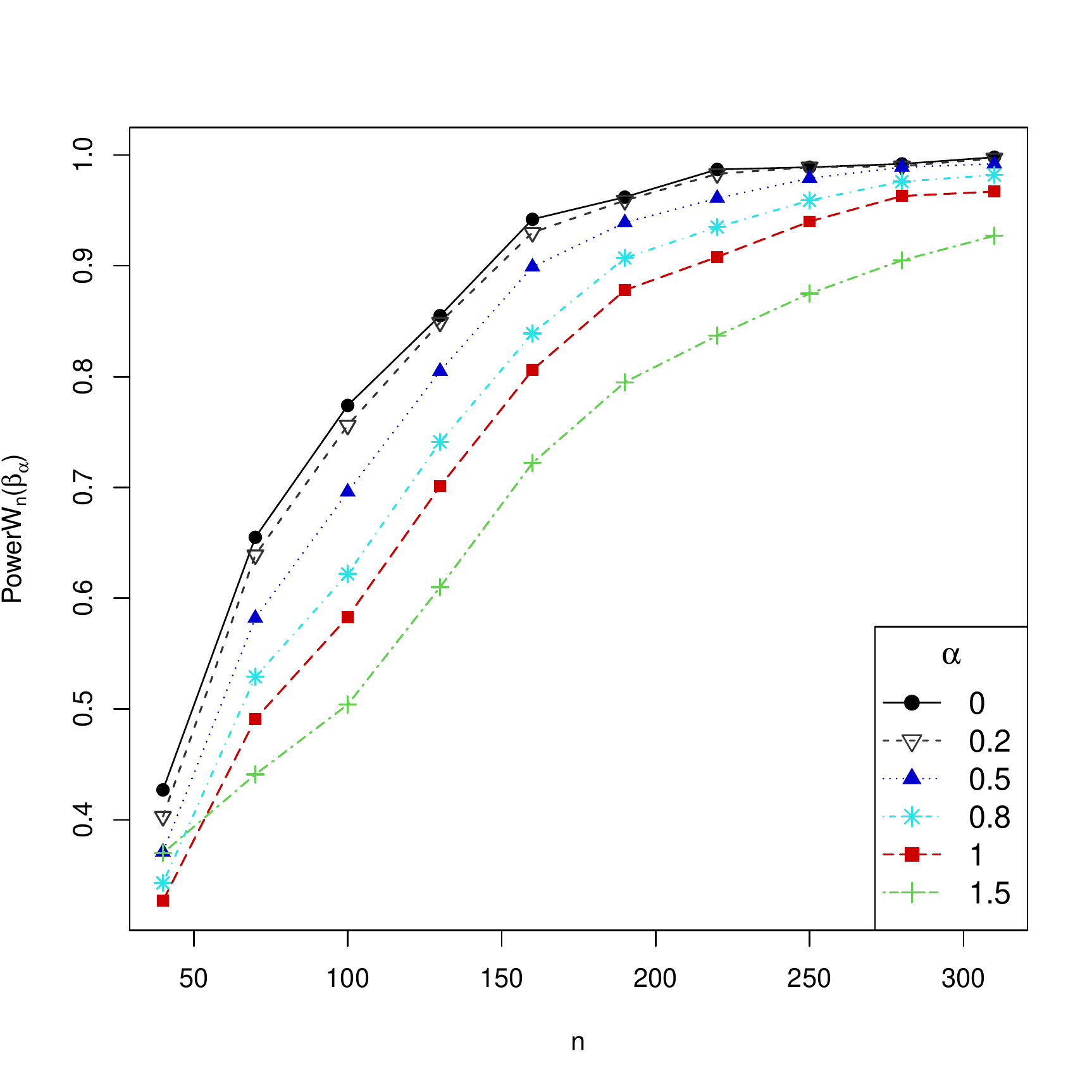}&
	\includegraphics[scale=0.4]{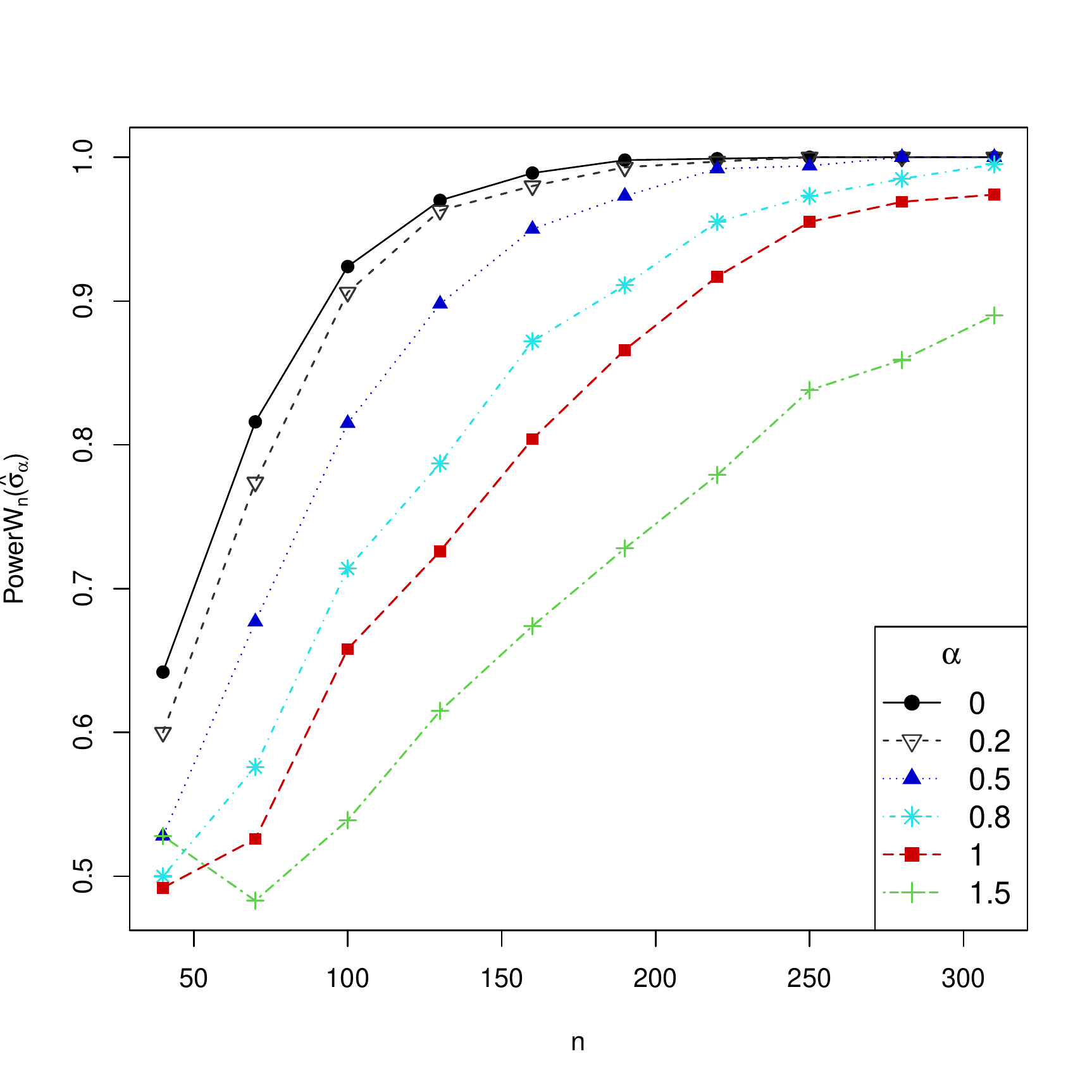}
	\end{tabular} 
	\caption{RMSE (top) empirical level (middle) and empirical power(bottom) against sample size for the null hypothesis (\ref{eq:hypothesis2}) (left) and (\ref{eq:hypothesis3}) (right) for the corresponding Wald type tests with pure data and \ref{itm:Design1}.}
	\label{d1e0}
\end{figure}

\begin{figure}
\begin{tabular}{cc}
	\includegraphics[scale=0.4]{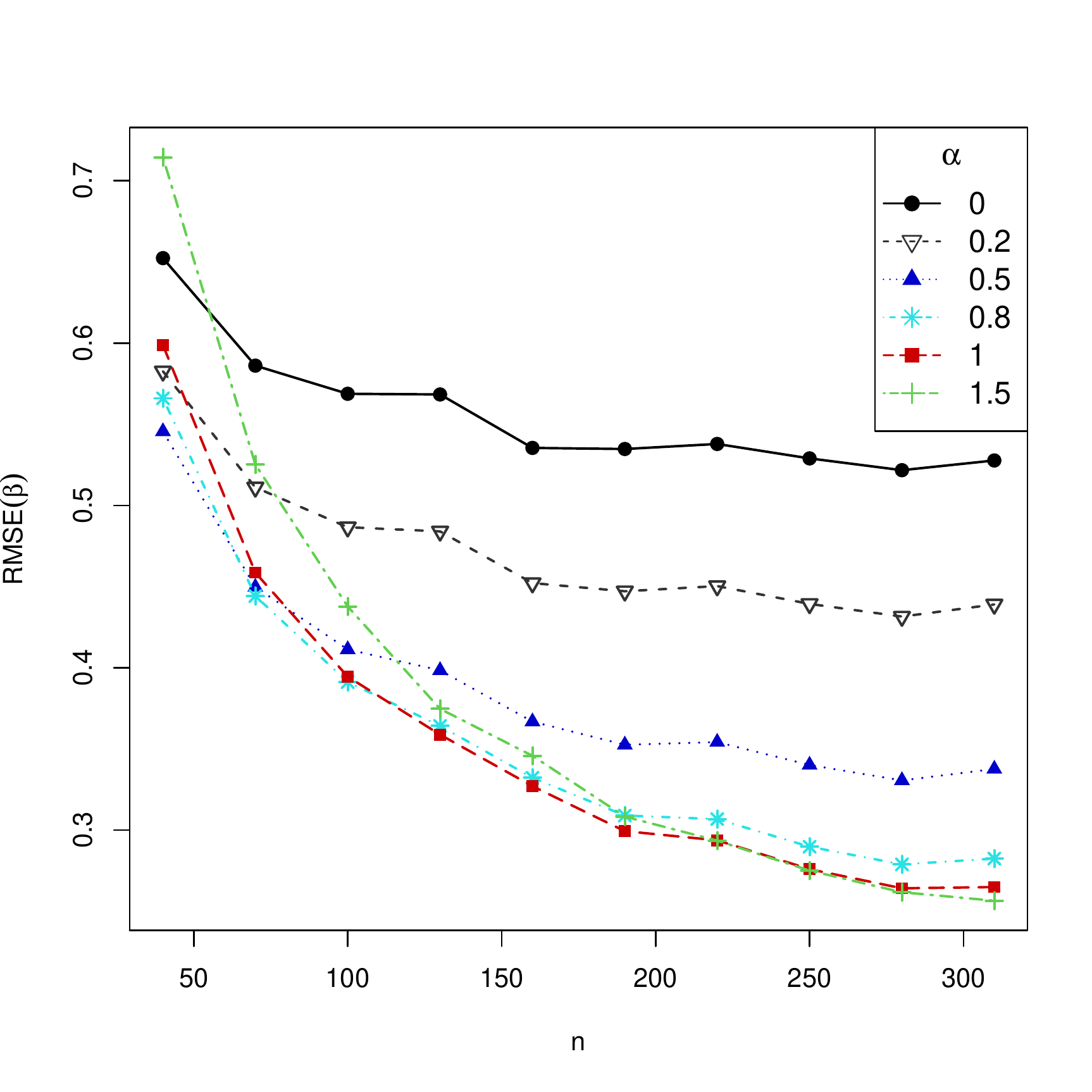}&
	\includegraphics[scale=0.4]{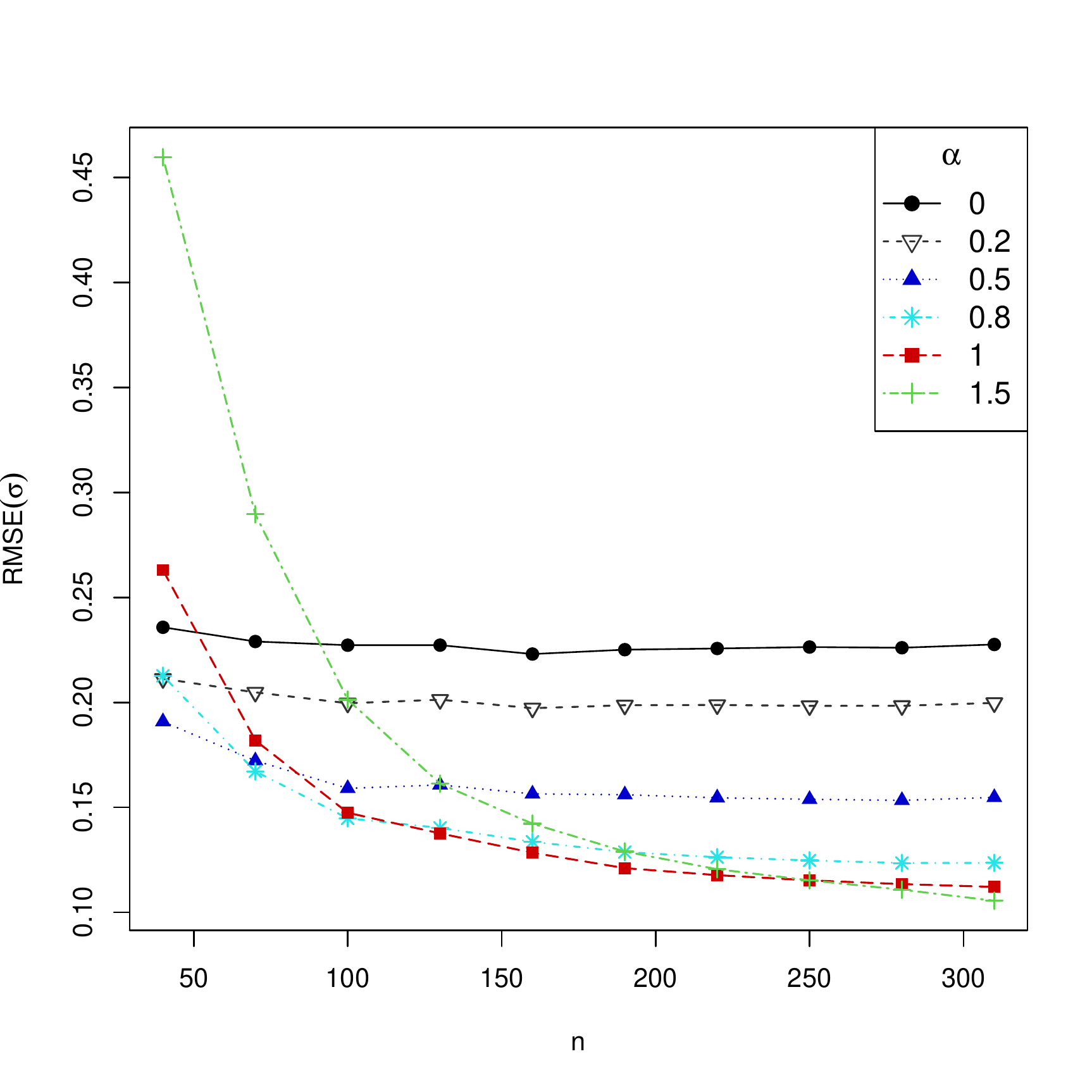}\\
	\includegraphics[scale=0.4]{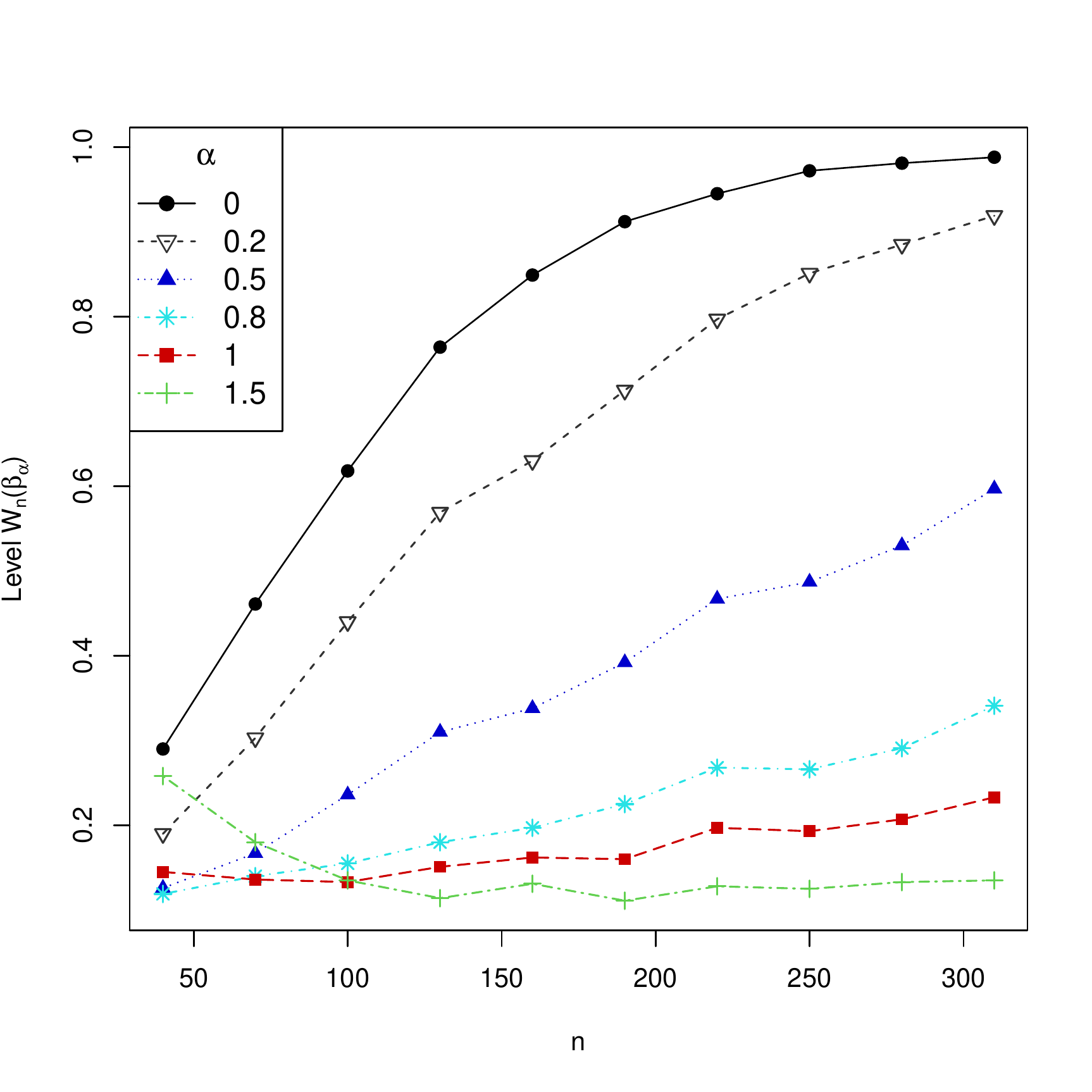}&
	\includegraphics[scale=0.4]{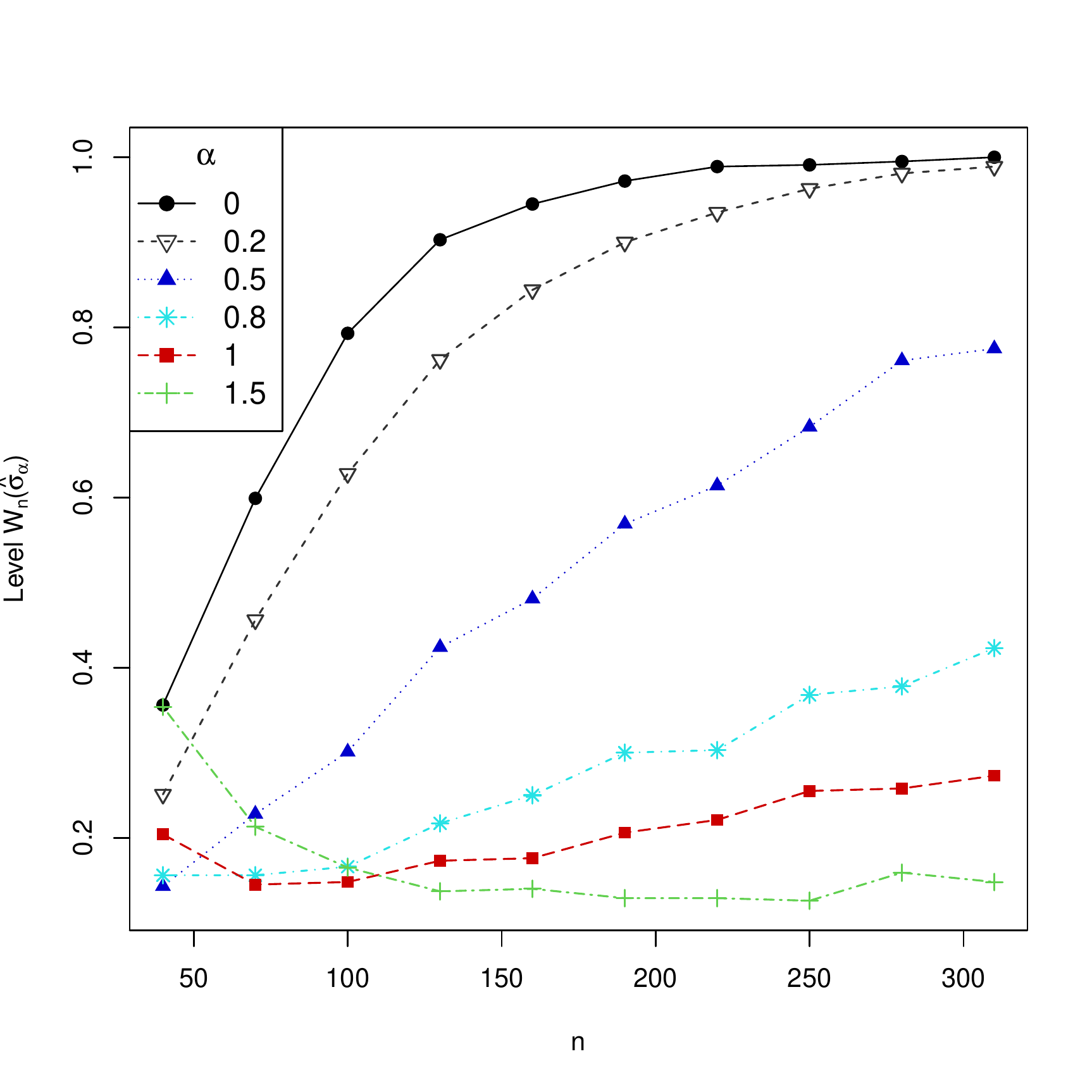}\\
	\includegraphics[scale=0.4]{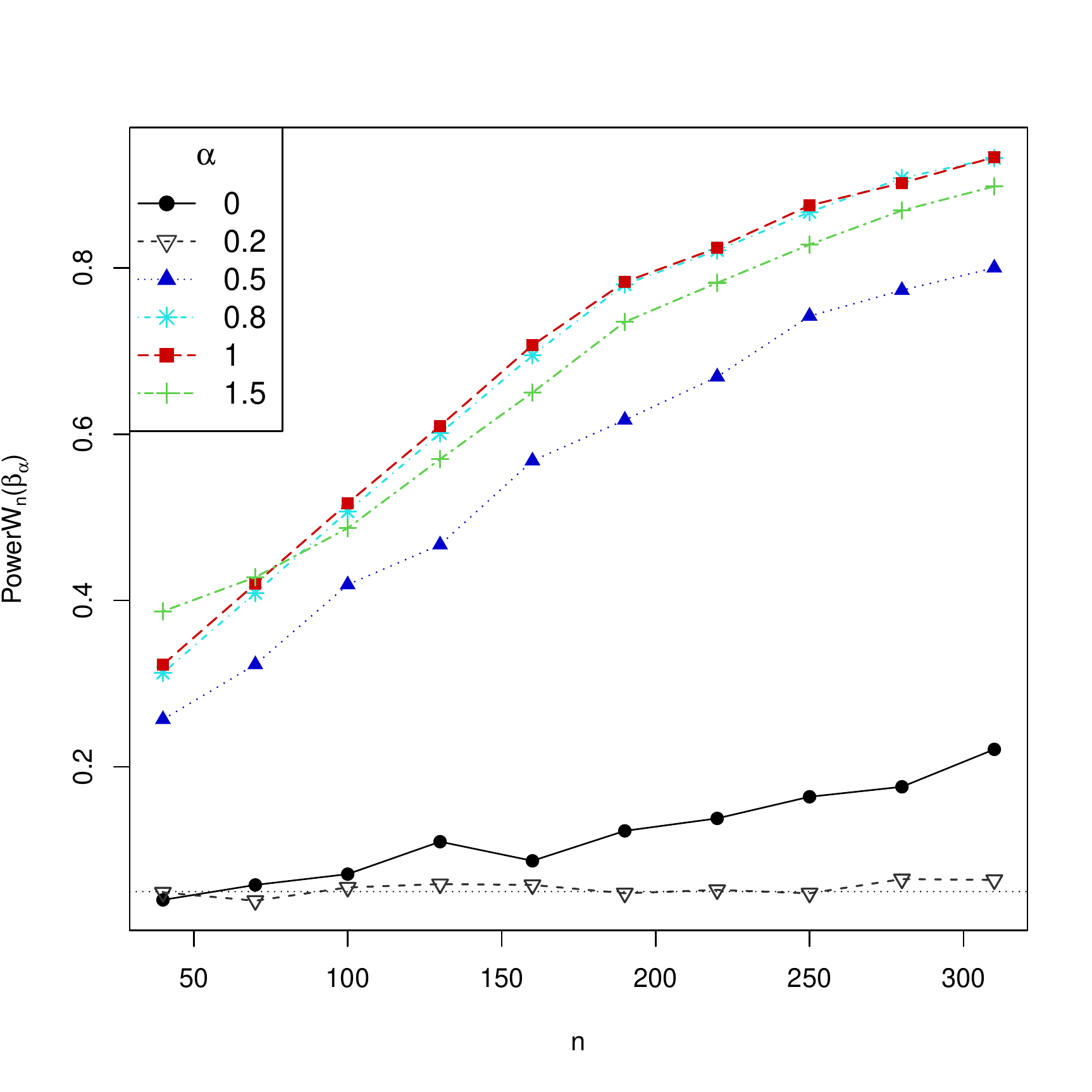}&
	\includegraphics[scale=0.4]{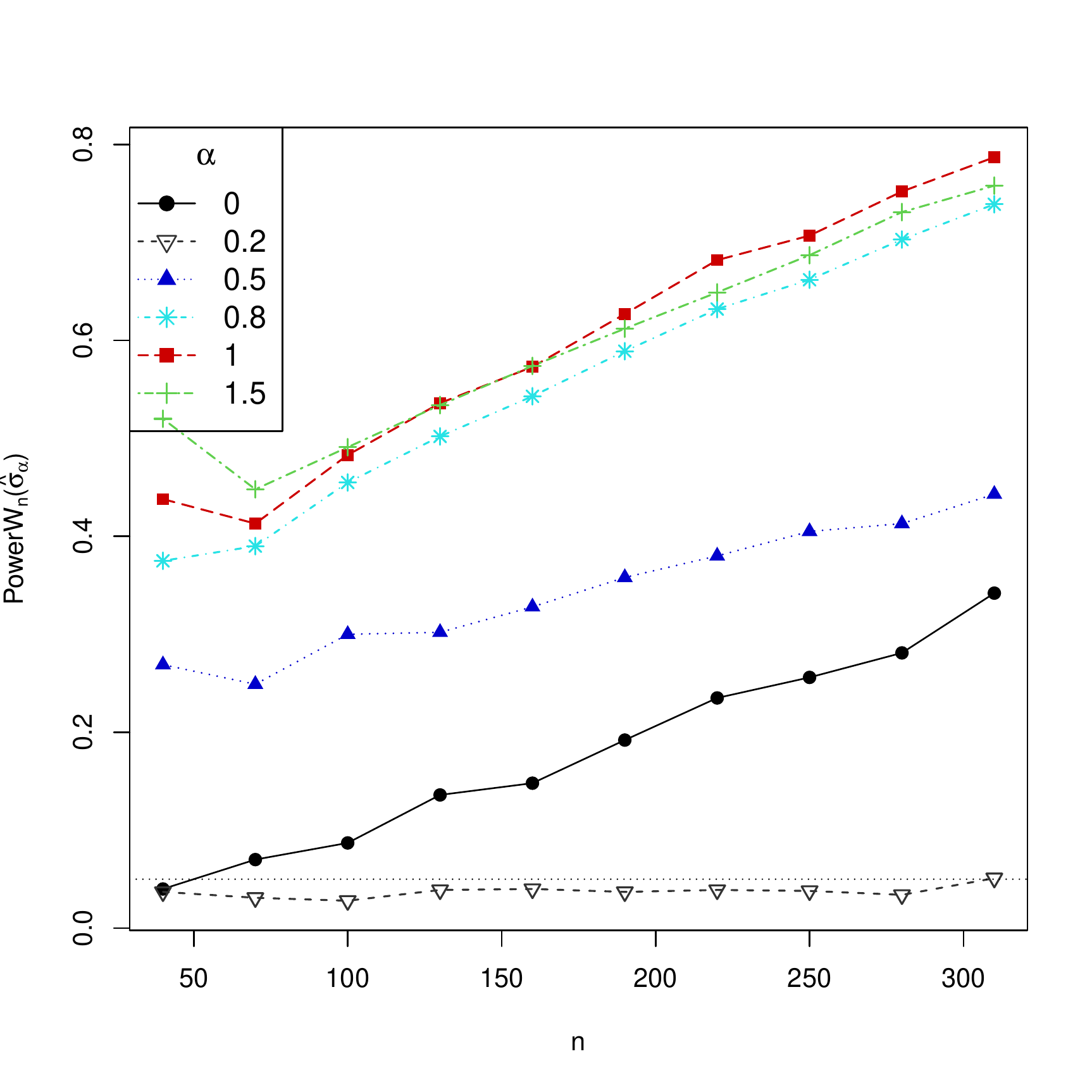}
		\end{tabular} 
	\caption{RMSE (top) empirical level (middle) and empirical power(bottom) against sample size for the null hypothesis  (\ref{eq:hypothesis2}) (left) and (\ref{eq:hypothesis3}) (right) for the corresponding Wald-type tests with $10\%$ of outliers and \ref{itm:Design2}.}
	\label{d1e1}
\end{figure}

\begin{figure}
\begin{tabular}{cc}
	\includegraphics[scale=0.4]{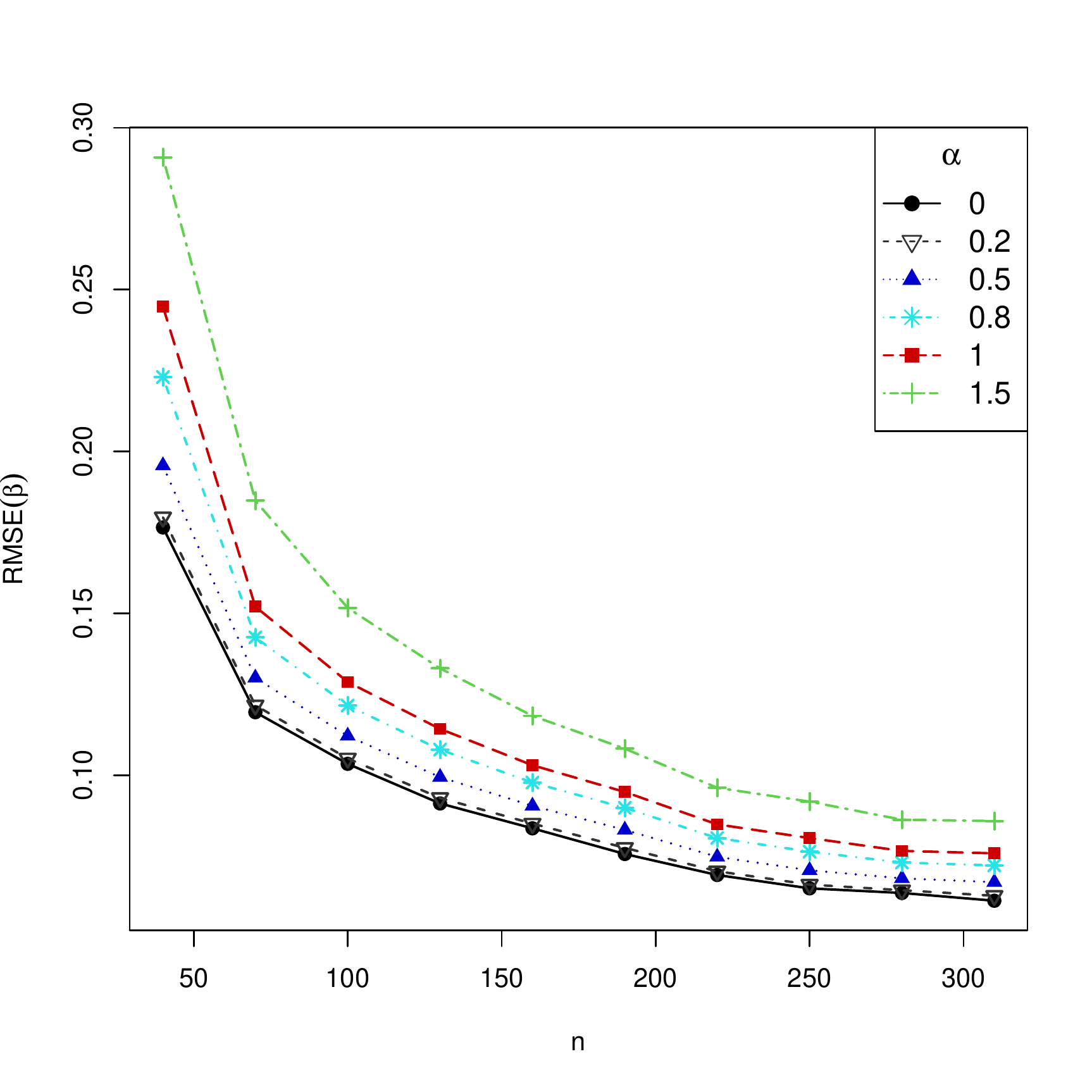}&
	\includegraphics[scale=0.4]{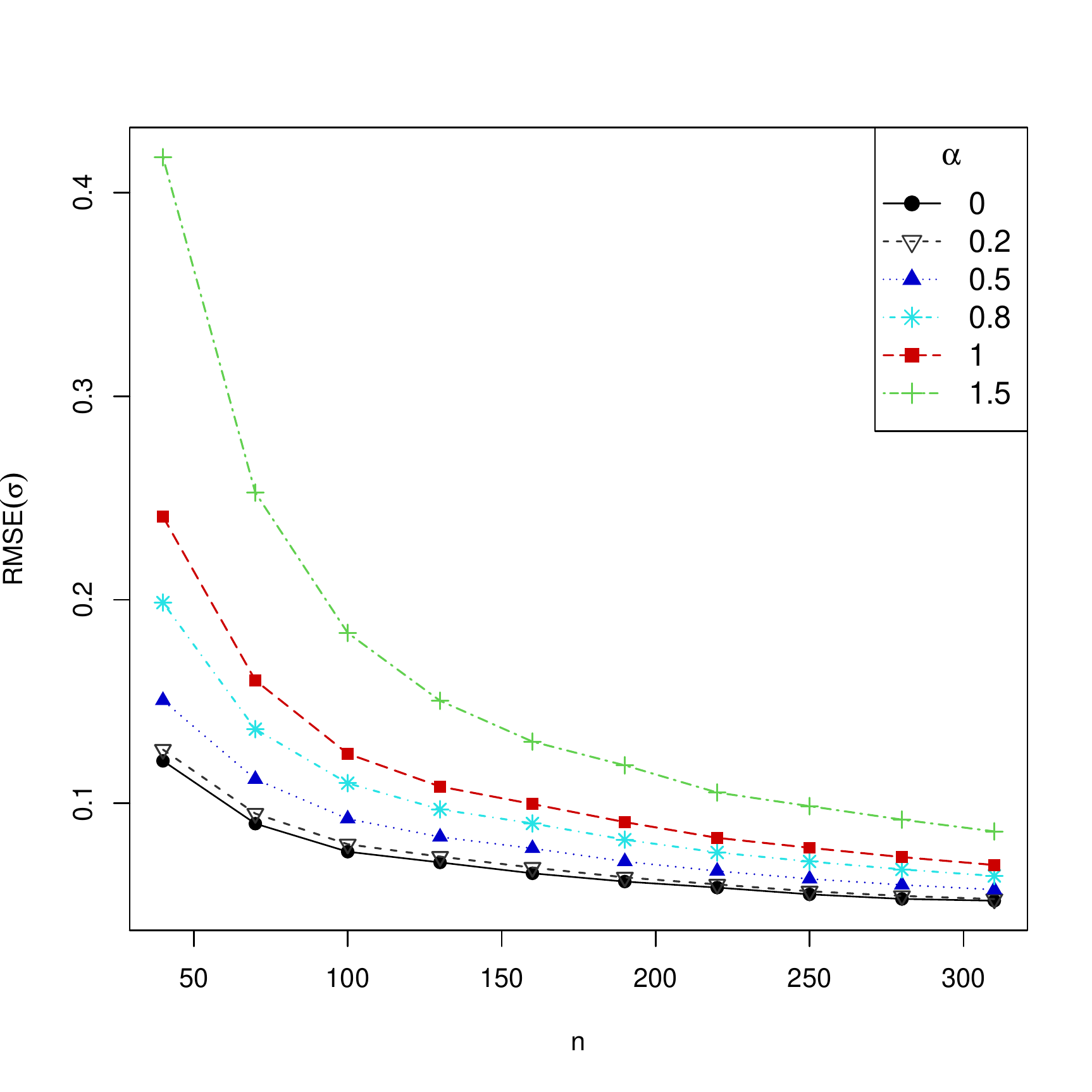}\\
	\includegraphics[scale=0.4]{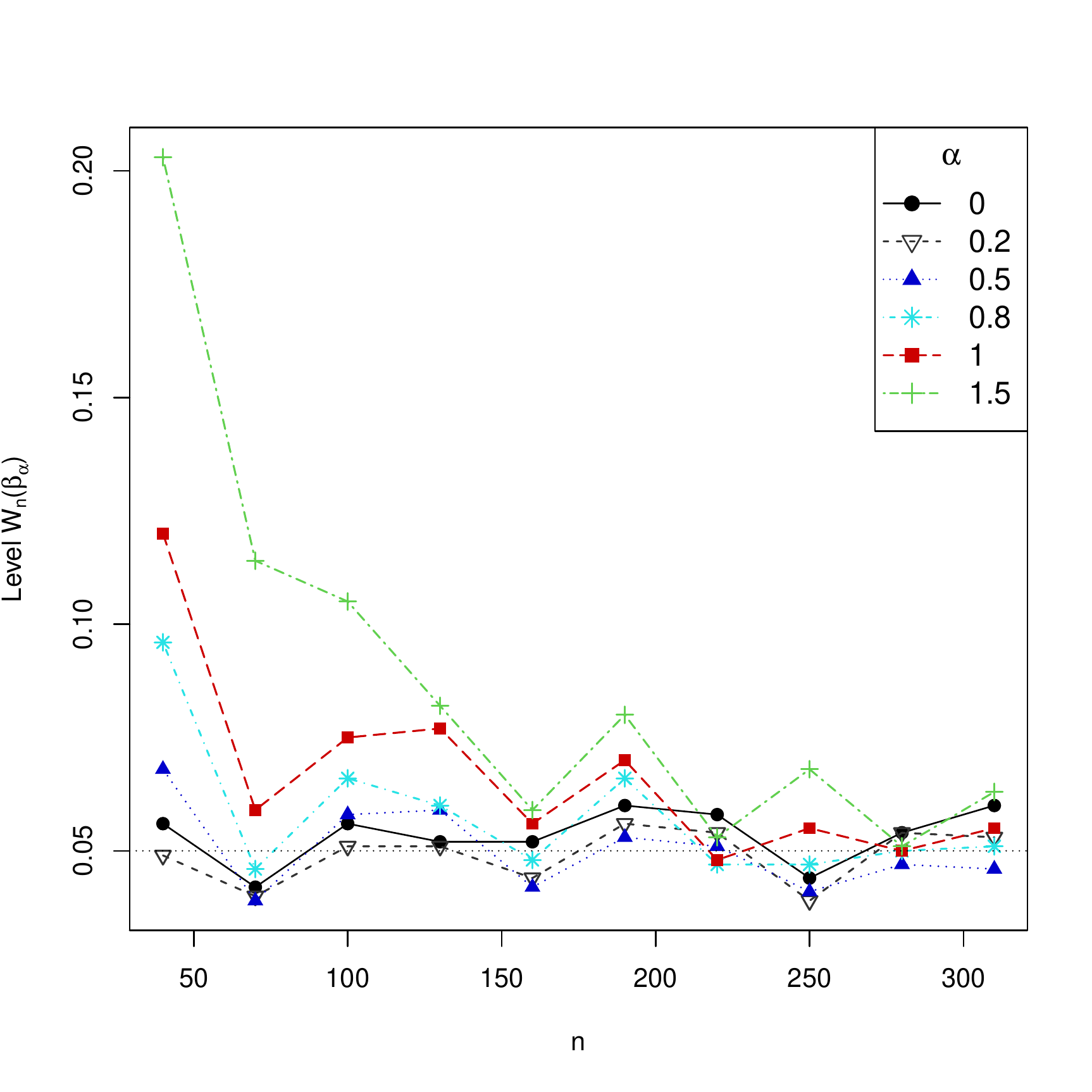}&
	\includegraphics[scale=0.4]{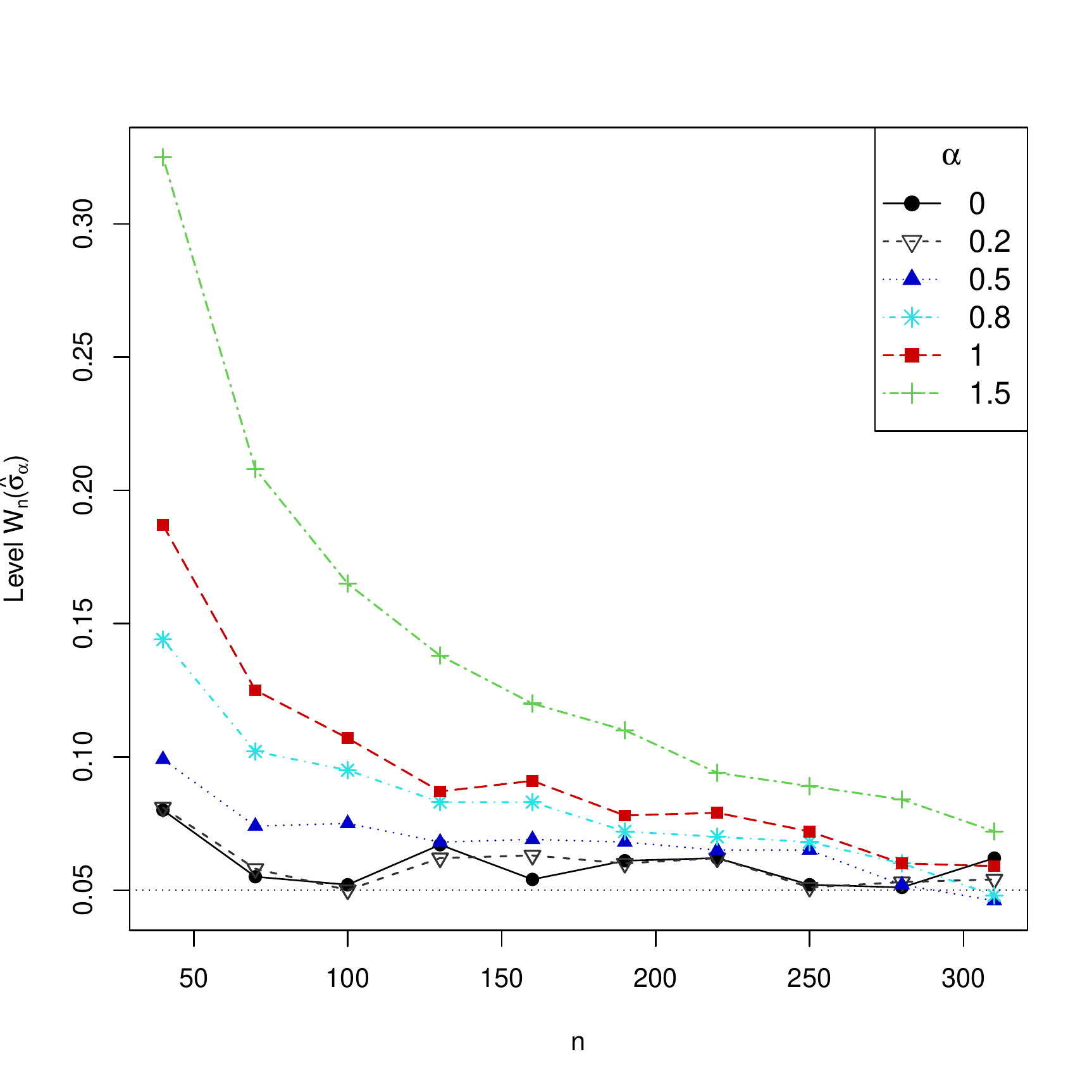}\\
	\includegraphics[scale=0.4]{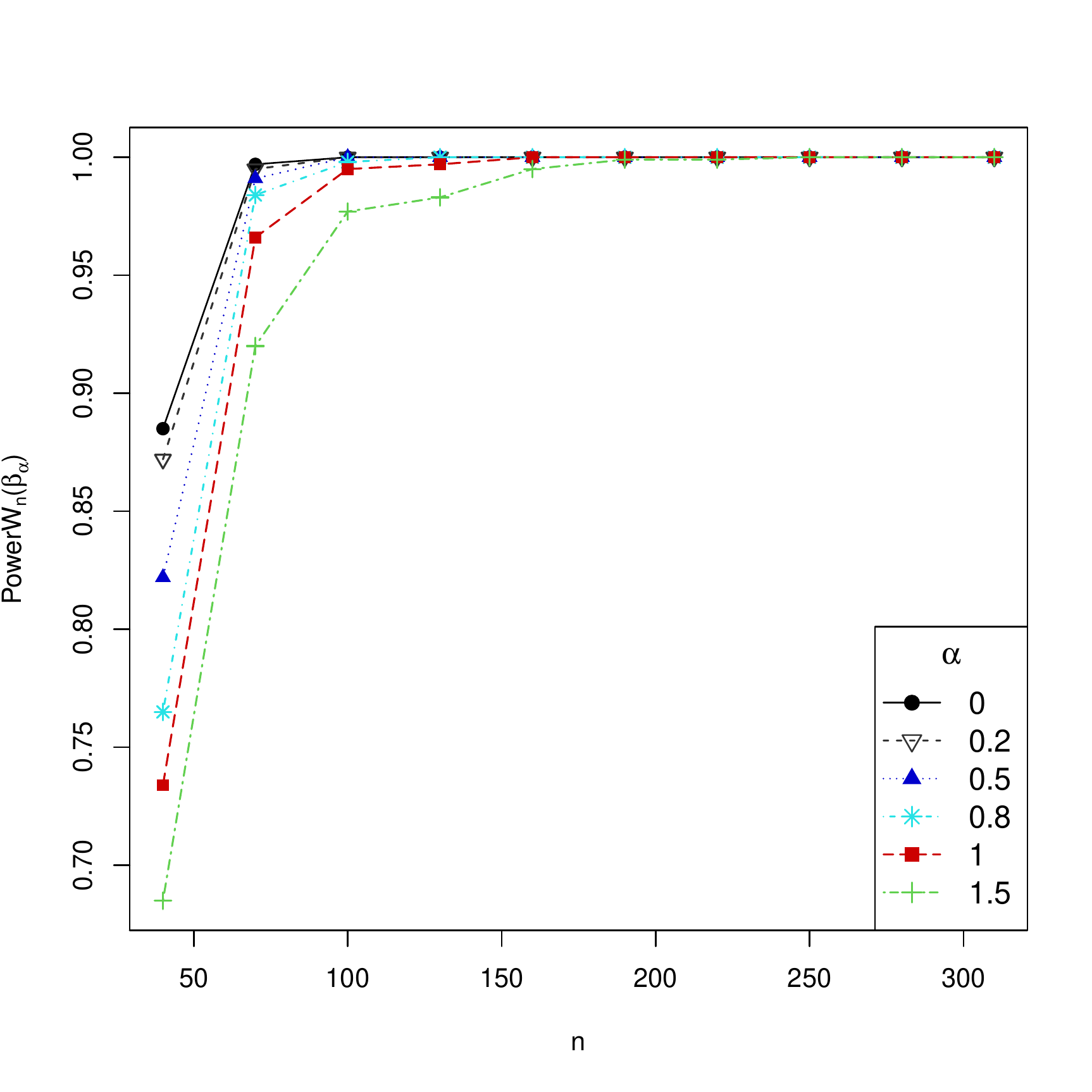}&
	\includegraphics[scale=0.4]{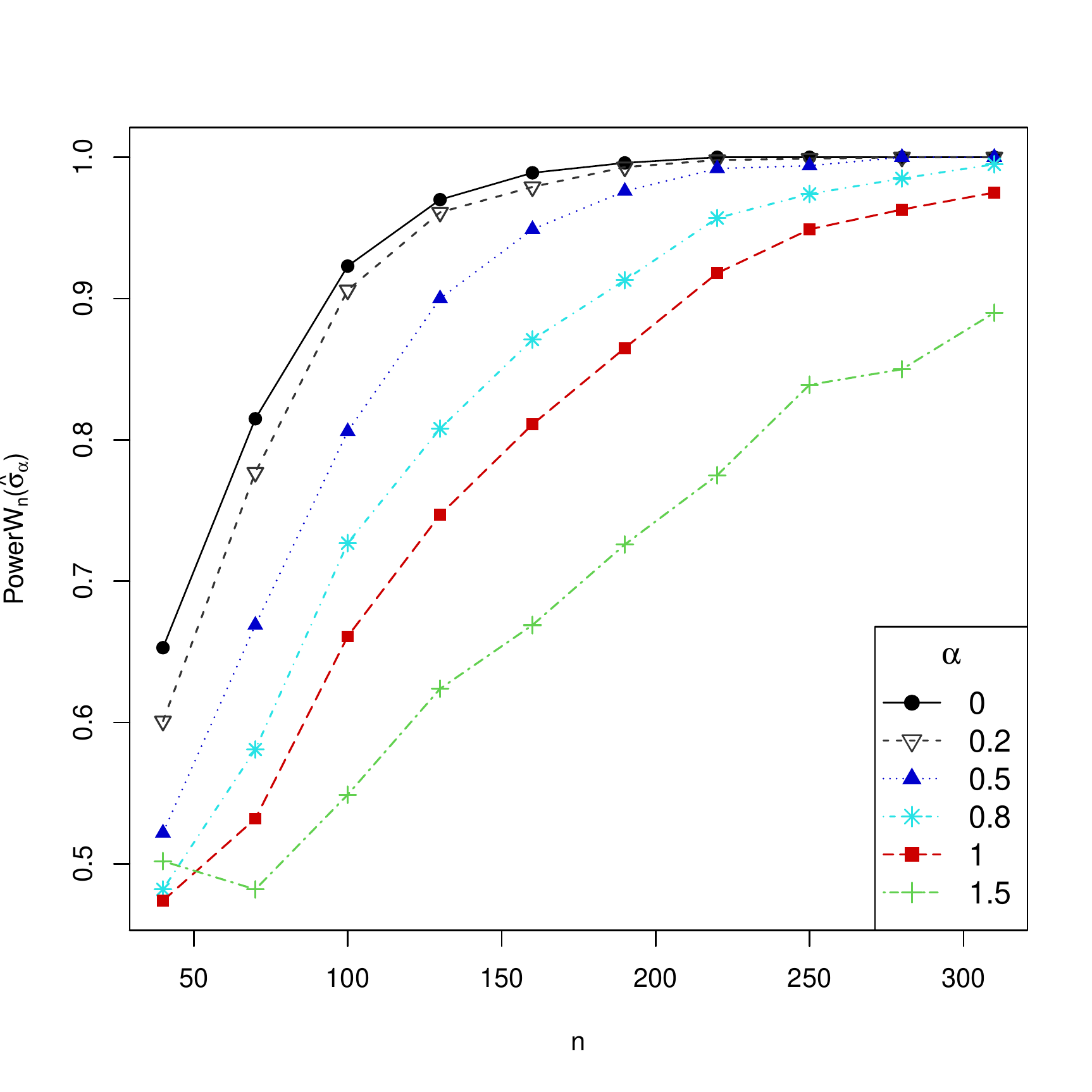}
			\end{tabular} 
	\caption{RMSE (top) empirical level (middle) and empirical power(bottom) against sample size for the null hypothesis (\ref{eq:hypothesis2}) (left) and (\ref{eq:hypothesis3}) (right) for the corresponding Wald-type tests with pure data and \ref{itm:Design2}.}
	\label{d2e0}
\end{figure}

\begin{figure}
\begin{tabular}{cc}
	\includegraphics[scale=0.4]{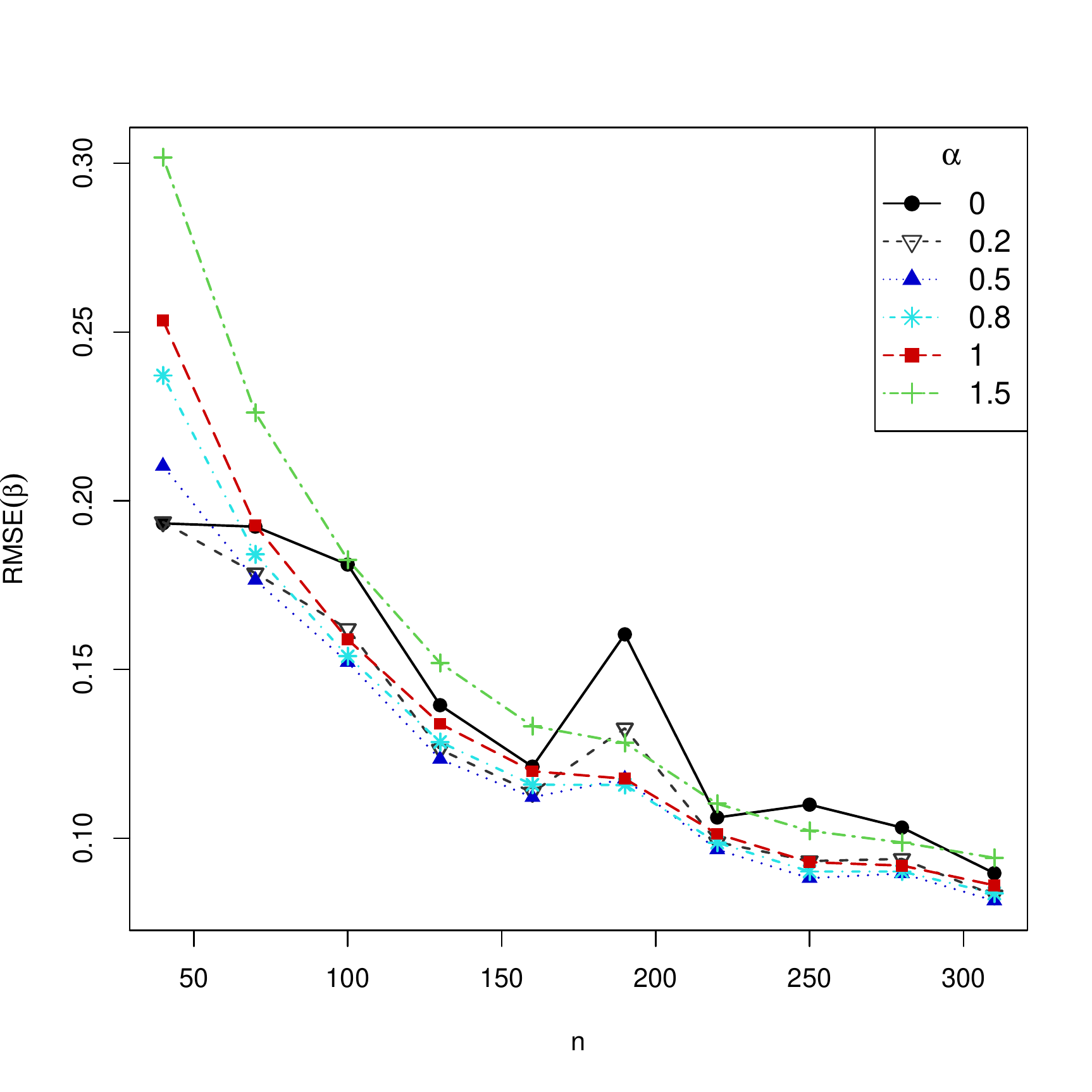}&
	\includegraphics[scale=0.4]{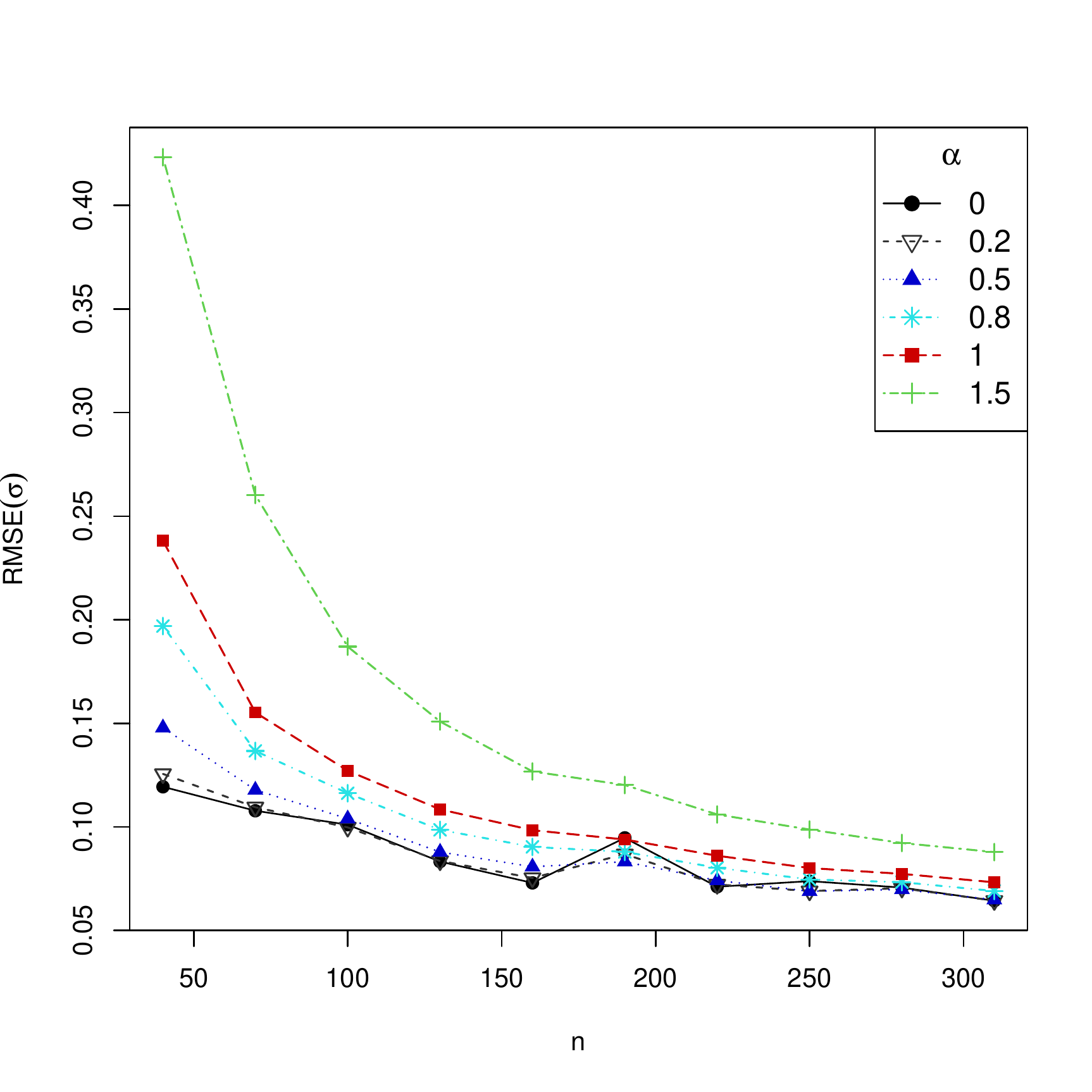}\\
	\includegraphics[scale=0.4]{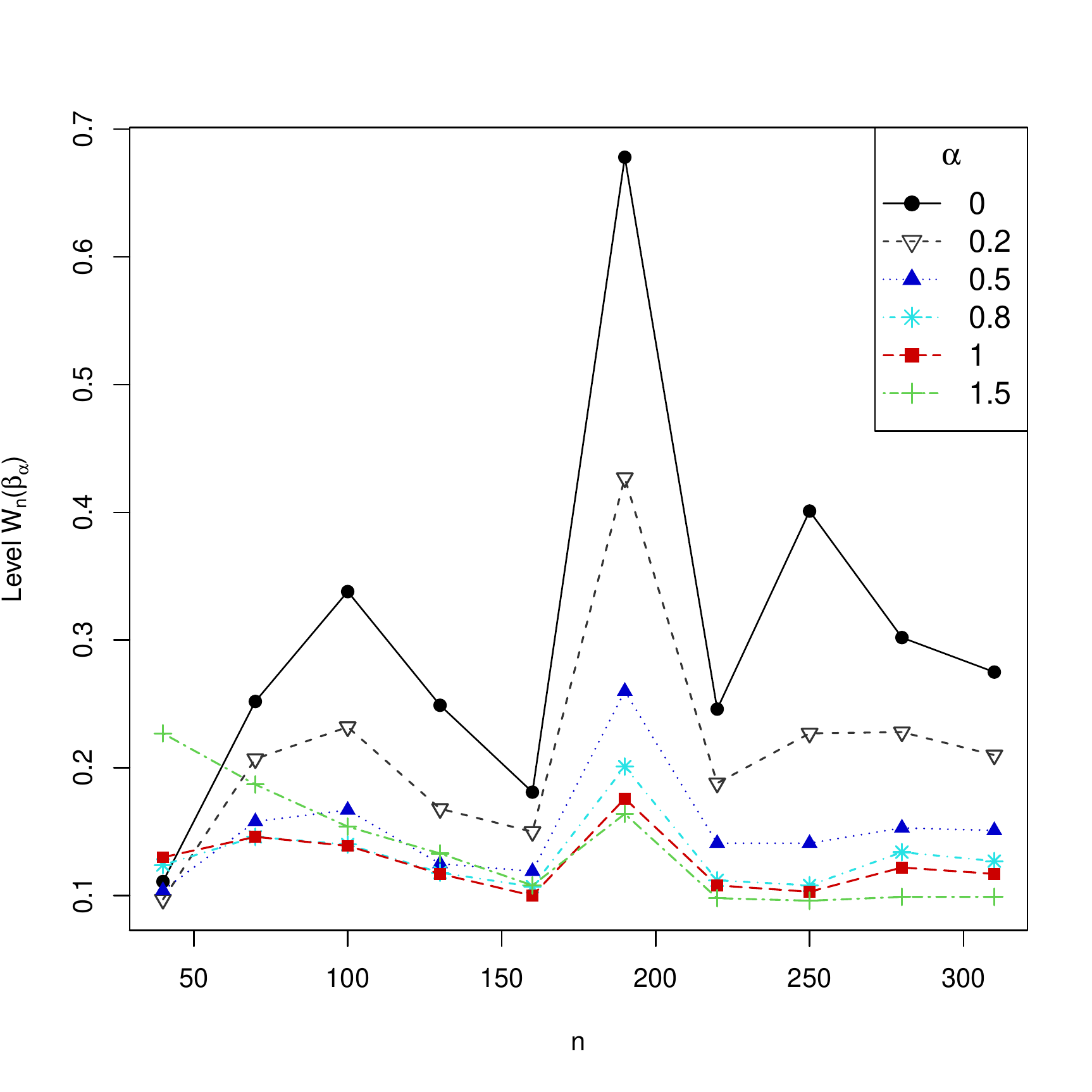}&
	\includegraphics[scale=0.4]{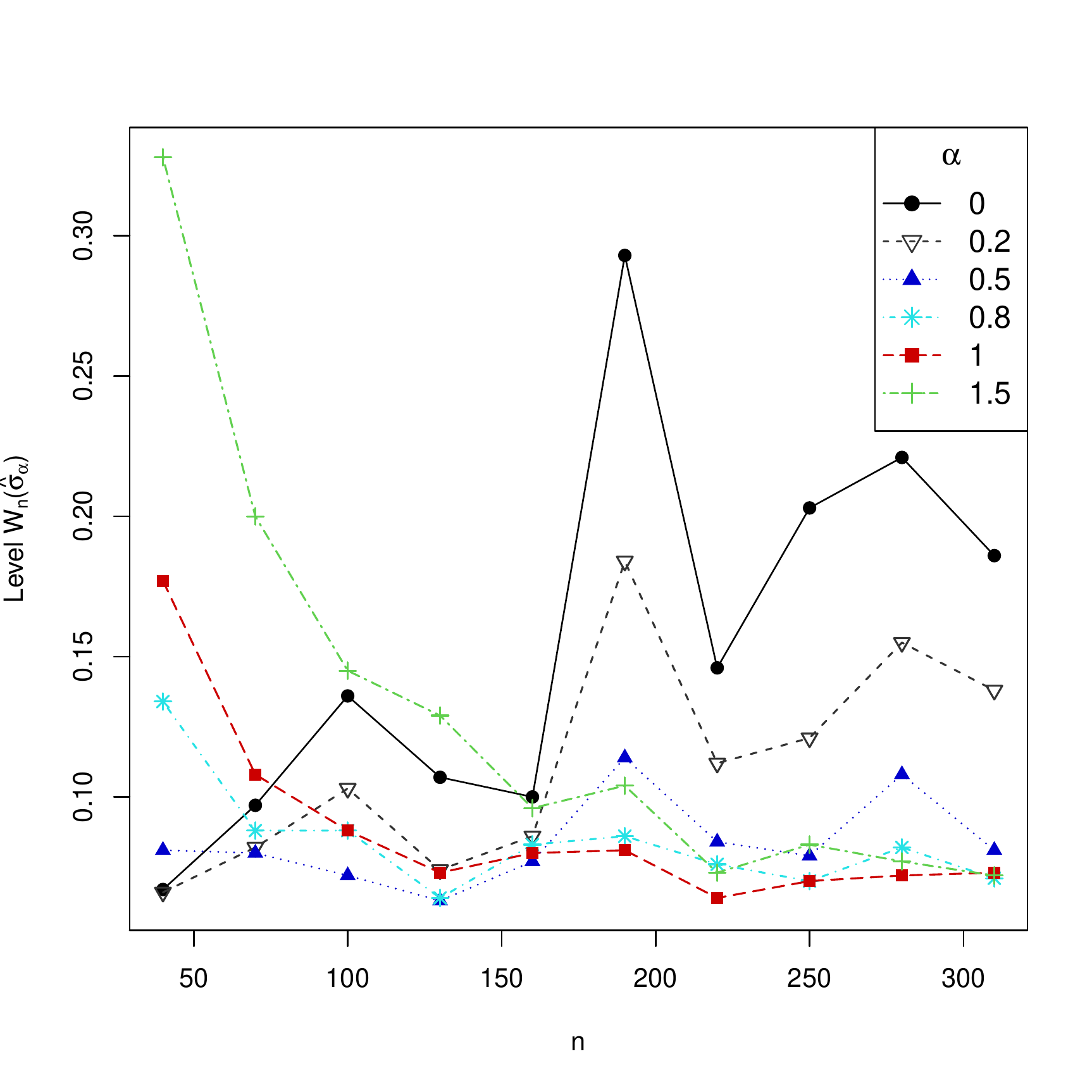}\\
	\includegraphics[scale=0.4]{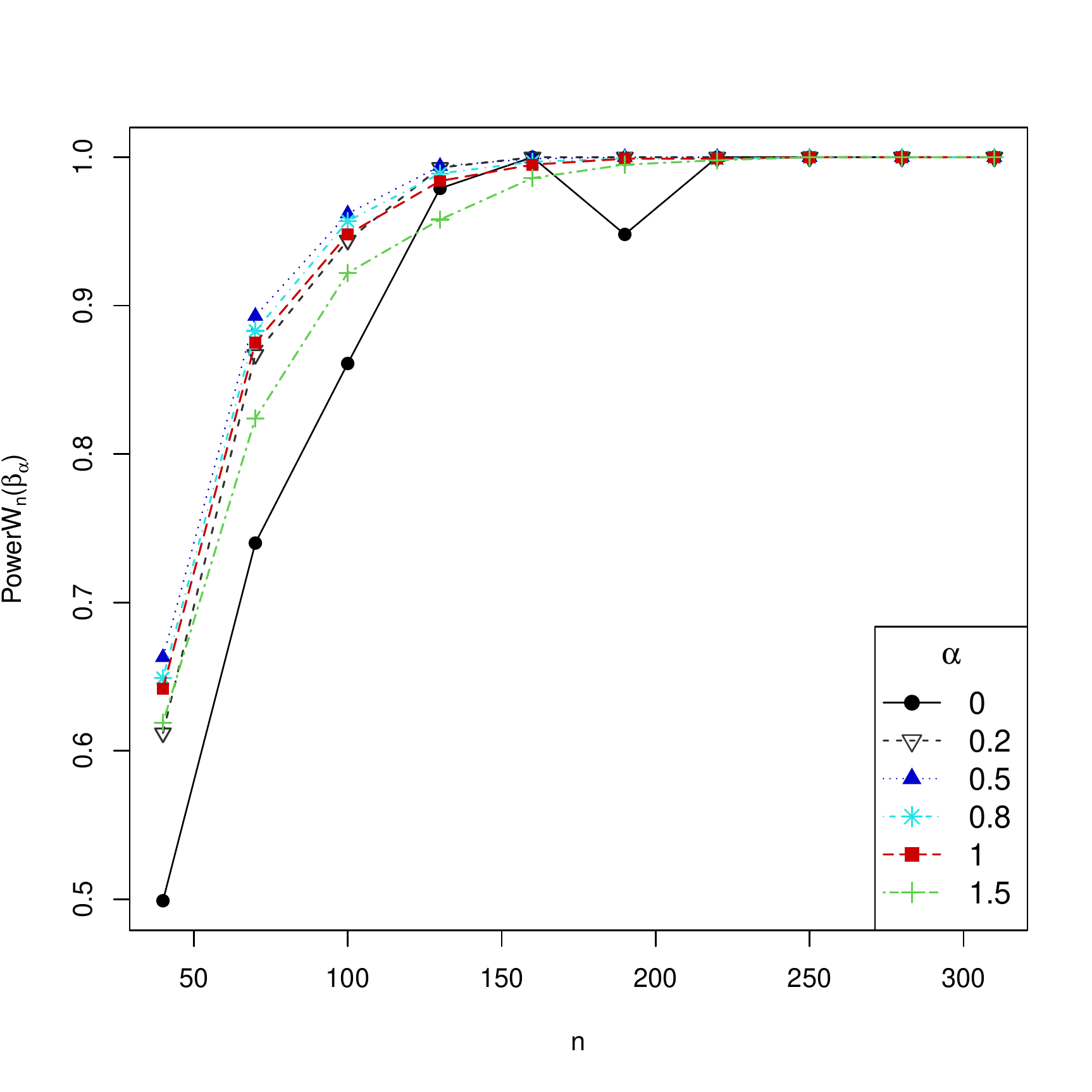}&
	\includegraphics[scale=0.4]{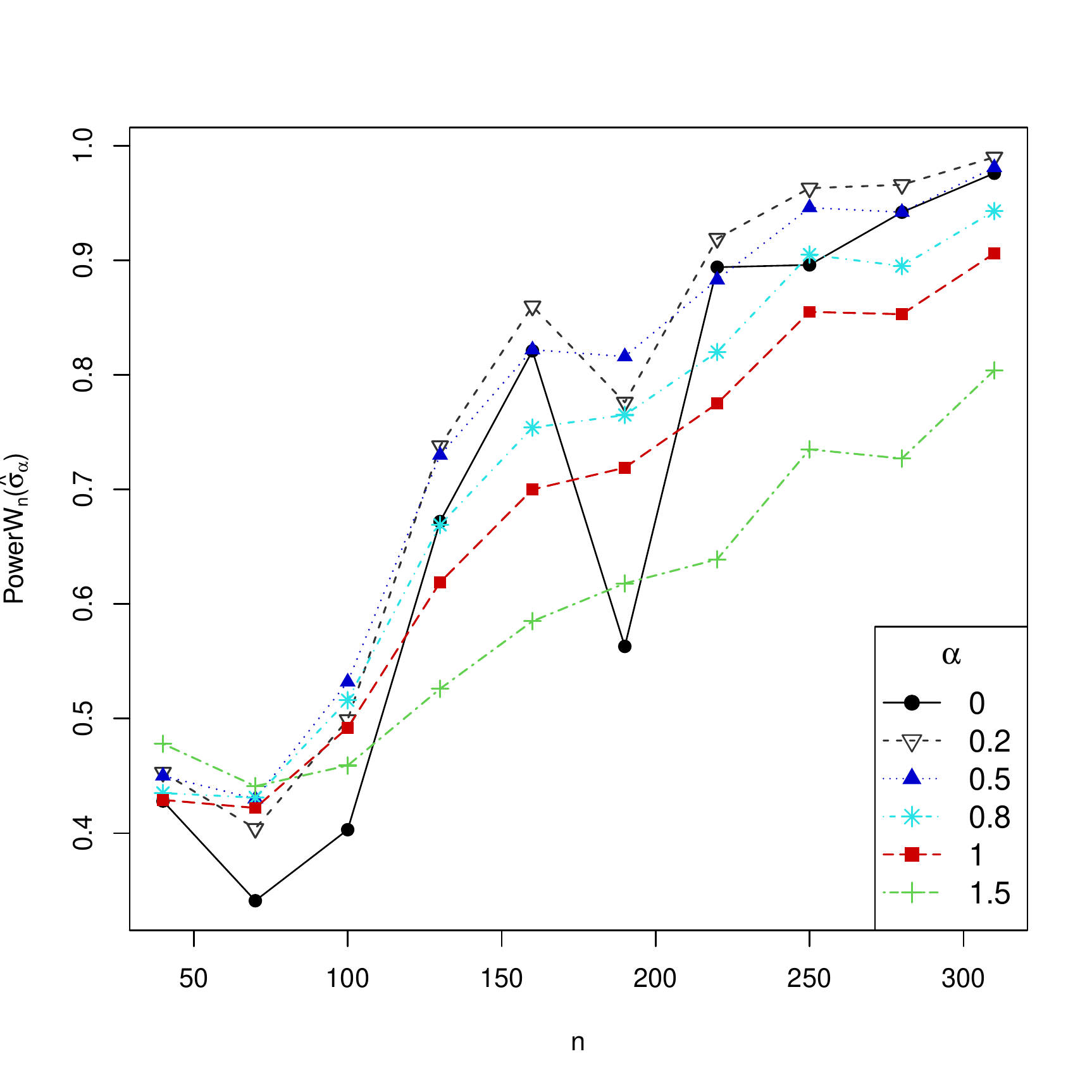}
				\end{tabular} 
	\caption{RMSE (top) empirical level (middle) and empirical power(bottom) against sample size for the null hypothesis  (\ref{eq:hypothesis2}) (left) and (\ref{eq:hypothesis3}) (right) for the corresponding Wald-type tests with $10\%$ of outliers and \ref{itm:Design2}.}
	\label{d2e1}
\end{figure}

We consider two different null hypotheses
\begin{align}
%\label{eq:hypothesis1} \operatorname{H}_0&: (\beta_0,\beta_1) = (1,1)\\ 
\label{eq:hypothesis2} \operatorname{H}_0&: \beta_1 = 1,\\ 
%\hspace{0.3cm} \text{and} \hspace{0.3cm} 
\label{eq:hypothesis3} \operatorname{H}_0&: \sigma = 1 ,
\end{align}
corresponding with  the composite null hypothesis
$$\operatorname{H}_0: \boldsymbol{M}^T\boldsymbol{\beta} = \boldsymbol{m},$$ with 
\begin{align*}
%\boldsymbol{M}_{\boldsymbol{\beta}} &= \begin{bmatrix} 1& 0&0\\ 0&1&0\\  \end{bmatrix}  \hspace{0.3cm}\text{and} \hspace{0.5cm} \boldsymbol{m}_{\boldsymbol{\beta}} = (1,1), \\
\boldsymbol{M}_{\beta_1} &= [0,1,0]  \hspace{0.5cm} \text{and} \hspace{0.5cm} \boldsymbol{m}_{\beta_1} = 1, \\
\boldsymbol{M}_\sigma &= [0,0,1] \hspace{0.5cm} \text{and}  \hspace{0.5cm}\boldsymbol{m}_\sigma = 1,
\end{align*}
%$$ \boldsymbol{M}_\beta = \begin{bmatrix} 1& 0&0\\ 0&1&0\\  \end{bmatrix} \hspace{0.3cm} \text{and} \hspace{0.3cm} \boldsymbol{m}_\beta = (1,1)$$ and $\boldsymbol{M}_\sigma = [0,0,1]$  and  $\boldsymbol{m}_\sigma = 1$ 
respectively. The Wald-type test statistics for testing (\ref{eq:hypothesis2})-(\ref{eq:hypothesis3}), which we will denote $W_n(\widehat{\beta_1}_{\alpha})$ and $W_n(\widehat{\sigma}_{\alpha})$, are given in (\ref{eq:Wald_composite}) by substituting the corresponding matrices. So as to investigate the trade-off between efficiency and  robustness depending on the tuning parameter $\alpha$, we compute the empirical levels for the proposed Wald-type tests and powers when the true parameter values are $\beta_1^0 = 0.45$  and $\sigma^0=0.8$, respectively. These levels and powers are computed as the number of times that the null hypothesis is rejected over the total simulated samples $R=1000.$ 
%For the seek of completeness we also estimate the levels and powers for the Wald type test based on the density power divergence (DPD), which also depends on a parameter $\alpha$ controlling the trade-off between efficiency and robustness. 
Figures \ref{d1e0}-\ref{d2e1} contain the root mean square error (RMSE), empirical level and power results for the null hypothesis tests (\ref{eq:hypothesis2}) and (\ref{eq:hypothesis3}), for a $5\%$ significance level. 
%The results for the Wald type test for $\boldsymbol{\beta}$ with both robust losses, RP and DPD, are very similar, so we do not include them for the seek of brevity.
The results show the clear improvement in robustness when $\alpha$ increases, in detriment to the efficiency. 
%Additionally, both robust methods performs similarly in the hypothesis test for $\boldsymbol{\beta},$ but the proposed Wald test statistic based on the minimum RP estimator produces higher powers for the hypothesis test respect to the standard error $\sigma.$
The MLE produces the best performance with pure data, showing its major efficiency, and the behaviour of the minimum RP estimator improves when $\alpha$ decreases, i.e., estimators based in low values of the parameter enjoy greater efficiency. However, in presence of data contamination, the RMSE and empirical level of the Wald-type test statistics rise for low values of $\alpha$, highlighting its lack of robustness. The most revealing setting is \ref{itm:Design1}, at which the empirical level and power of the Wald-type tests based on the MLE reaches their worst results, but the proposed Wald-type test based on RP loss statistics continues to perform adequately for sufficiently high values of $\alpha$.


For the first hypothesis test (\ref{eq:hypothesis2}), we could apply Theorem \ref{thm:contiguous} to obtain the power under the contiguous alternative hypothesis (\ref{eq:contiguous_hypothesis}). The distribution of the Wald-type tests $W_n(\widehat{\beta}_1)$ is given by a chi-squared with 1 degree of freedom and non-centrality parameter
$$\delta = \frac{\sigma(\alpha+1)^3}{(2\alpha+1)^{3/2}} \boldsymbol{d^\ast}^T\left(\frac{1}{n}\sum_{i=1}^n X_{1i}^2\right)^{-1} \boldsymbol{d^\ast},$$
depending on the standard deviation error $\sigma$, the tuning parameter $\alpha$ and the fixed value $d_x= \boldsymbol{d^\ast}^T\left(\frac{1}{n}\sum_{i=1}^n X_{1i}^2\right)^{-1} \boldsymbol{d^\ast}.$ The choice $d_x = 0$ corresponds with the level of the test. Table \ref{table:contigous} summarizes the empirical power results over different values of $\alpha$ and $d_x,$ with $\sigma =1.$

\begin{table}[htb]
	\center
	\caption{Empirical power values of the null hypothesis (\ref{eq:hypothesis2}) under contiguous hypothesis.}
	\label{table:contigous}
	\begin{tabular}{r|rrrrrrrr}
		\hline
		\multicolumn{9}{c}{$d_x$}\\
		\hline
		$\alpha$& 0 & 2 & 5 & 10 & 15 & 20 & 25 & 30 \\ 
		\hline
		0 & 0.05 & 0.28 & 0.59 & 0.88 & 0.97 & 0.99 & 1.00 & 1.00 \\ 
		0.2 & 0.05 & 0.27 & 0.58 & 0.86 & 0.97 & 0.99 & 1.00 & 1.00 \\ 
		0.5 & 0.05 & 0.25 & 0.52 & 0.81 & 0.94 & 0.98 & 1.00 & 1.00 \\ 
		0.8 & 0.05 & 0.22 & 0.44 & 0.75 & 0.90 & 0.97 & 0.99 & 1.00 \\ 
		1 & 0.05 & 0.21 & 0.41 & 0.71 & 0.87 & 0.95 & 0.98 & 0.99 \\ 
		1.5 & 0.05 & 0.17 & 0.35 & 0.60 & 0.78 & 0.89 & 0.95 & 0.97 \\ 
		\hline
	\end{tabular}
\end{table} 
Note that greater values of $d_x$ produces greater power values as expected, and empirical power decreases with $\alpha.$ However, the efficiency loss is not very significant in comparison to the robustness advantage.

\section{Real Data examples \label{sec:data}}

\subsection{Brain and Weight Data}
\begin{figure}[h]
\centering
	\includegraphics[scale=0.65]{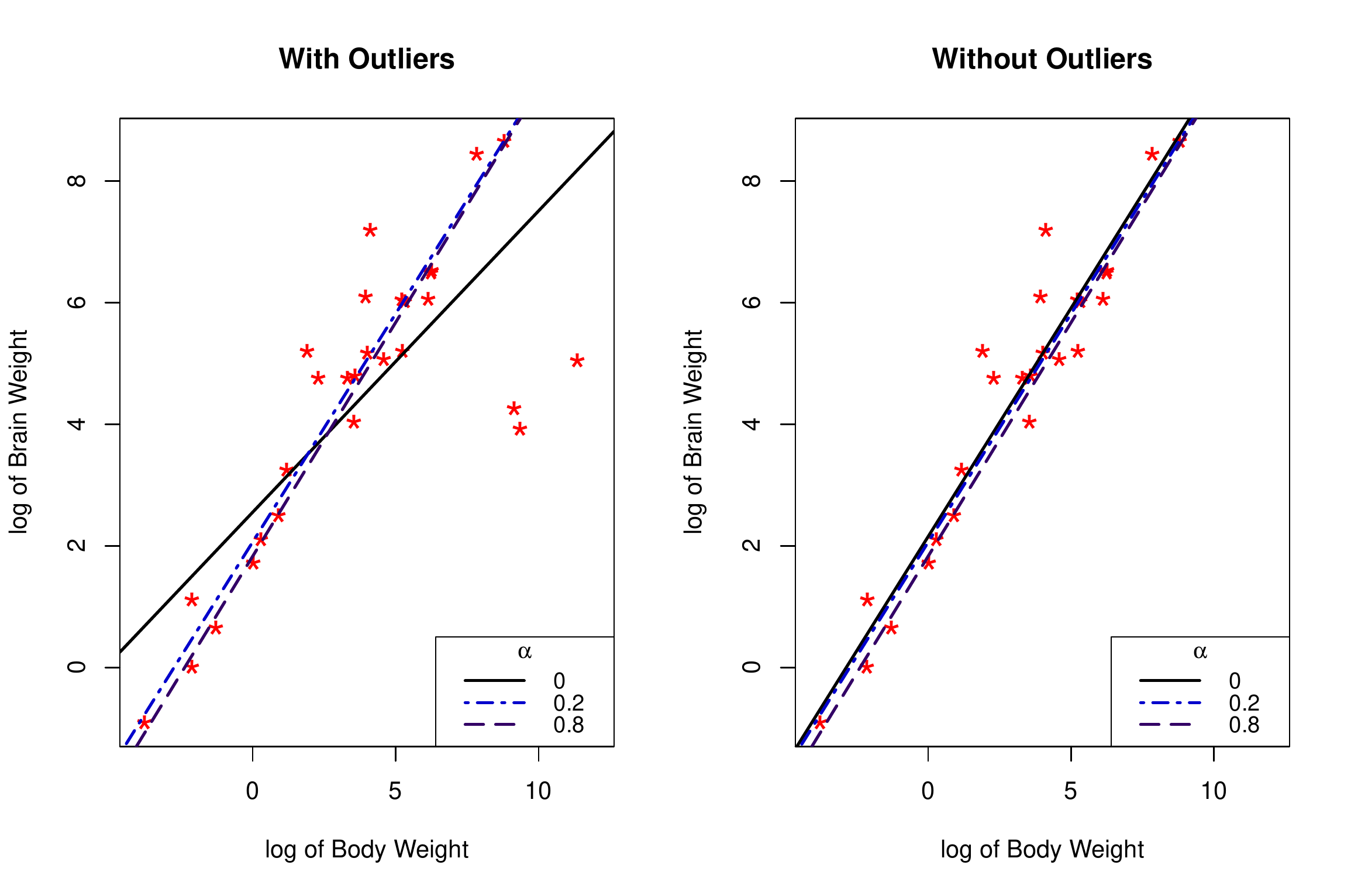}
	\caption{Plots of the data-points and fitted regression lines for the Brain and Weight Data using several minimum RP estimators before and after deleting the outliers.}
	\label{fig:animals1}
\end{figure}

These data, adapted from a larger data set in \cite{Weisberg}, were presented in Rousseeuw and Leroy (\cite{Rousseeuw}, pp. 57) as an example of the unrobustness of the classical MLE in simple linear regression. In this sample, the body weight (in kilograms) and the brain weight (in grams) of $n=28$ animals are compared, to investigate if a larger brain is required to govern a heavier body. As suggested in in \cite{Rousseeuw}, a transformation should be done to clearly represent either  the larger or smaller measurements. In this case, we take the Napierian logarithm of both brain and body weights.  Observations $6,16$ and $25$, those corresponding to dinosaurs,  posses an unusual small brain as compared with a heavy body, which clearly affects to the slope of the classical estimation method ($\alpha=0$) as can be seen in Figure \ref{fig:animals1} and Table \ref{table:animals1}. The estimates of the regression coefficients  and the error variance obtained from the minimum RP estimation for various $\alpha$  are also presented here, observing how the estimation based on $\alpha>0$ is more robust to the presence of these outliers.

%\begin{table}[h]
%\center
%	\caption{The parameter estimates of the linear regression model for the Brain and Weight Data using several minimum RP estimators.}
%	\label{table:animals1}	
%\begin{tabular}{c rrr rrr r}
%\hline
%$\alpha$     & 0      & 0.2    & 0.4    & 0.6    & 0.8    & 1      &  1.2\\ \hline
%$\sigma$ & 1.4714 & 0.6410 & 0.4929 & 0.4092 & 0.3640 & 0.3378 & 0.3184\\
%$\beta_1$ & 2.5523 & 2.0617 & 1.9378 & 1.8616 & 1.8265 & 1.8142 & 1.8092\\
%$\beta_2$ & 0.4958 & 0.7509 & 0.7560 & 0.7634 & 0.7694 & 0.7731 & 0.7756 \\ \hline
%\end{tabular}
%\end{table}

\begin{table}[h]
\center
	\caption{The parameter estimates of the linear regression model for the Brain and Weight Data using several minimum RP estimators.}
	\label{table:animals1}
\begin{tabular}{r rrr r rrr}
\hline
& \multicolumn{3}{c}{With outliers} & &\multicolumn{3}{c}{Without outliers} \\
\cline{2-4}  \cline{6-8} 
$\alpha$ & $\sigma$  & $\beta_0$  & $\beta_1$  && $\sigma$  & $\beta_0$  & $\beta_1$  \\ \hline
0     & 1.4714 & 2.5523 & 0.4958 && 0.6962 & 2.1504 & 0.7522 \\
0.2   & 0.6410 & 2.0617 & 0.7509 && 0.6309 & 2.0580 & 0.7519 \\
0.4   & 0.4929 & 1.9378 & 0.7560 && 0.4929 & 1.9378 & 0.7560 \\
0.6   & 0.4092 & 1.8616 & 0.7634 && 0.4092 & 1.8616 & 0.7634 \\
0.8   & 0.3640 & 1.8265 & 0.7694 && 0.3640 & 1.8265 & 0.7694 \\
1     & 0.3378 & 1.8142 & 0.7731 && 0.3378 & 1.8142 & 0.7731 \\
%1.2   & 0.3184 & 1.8092 & 0.7756 && 0.3184 & 1.8092 & 0.7756 \\
\hline
\end{tabular}
\end{table}

In order to compare the performance of the Wald-type test with different values of the tunning parameter $\alpha$, we consider the following tests
\begin{align}
\label{testanimal1} \operatorname{H}_0 &: \beta_0 = 1.98,\\ 
\label{testanimal2} \operatorname{H}_0 &: \beta_1 = 0.73,\\  
\label{testanimal3} \operatorname{H}_0 &: (\beta_0, \beta_1) = (1.98, 0.73) , 
\end{align}
where the values $1.98$ and $0.73$ are respectively the mean value of the estimated coefficients $\beta_0$ and $\beta_1$ with the different values of $\alpha$ and using the original data (with outliers) listed in Table \ref{table:animals1}. Table \ref{table:pvalueanimals} shows the p-values obtained by using the corresponding Wald-type test statistics, $W_n(\widehat{\beta}_0)$, $W_n(\widehat{\beta}_1)$ and $W_n(\widehat{\boldsymbol{\beta}}).$

\begin{table}[h]
	\centering
	\label{table:pvalueanimals}
	\caption{p-value obtained for the tests (\ref{testanimal1})-(\ref{testanimal3}) using the corresponding Wald-type test statistics.}
	\begin{tabular}{r rrr r rrr}
		\hline
		& \multicolumn{3}{c}{With outliers} & &\multicolumn{3}{c}{Without outliers} \\
		\cline{2-4}  \cline{6-8} 
		$\alpha$ & $W_n(\widehat{\beta}_0)$ & $W_n(\widehat{\beta}_1)$ & $W_n(\widehat{\boldsymbol{\beta}})$ 
		&& $W_n(\widehat{\beta}_0)$ & $W_n(\widehat{\beta}_1)$ & $W_n(\widehat{\boldsymbol{\beta}})$ \\ 
		\hline
	0 & 0.080 & 0.000 & 0.000 && 0.452 & 0.542 & 0.072 \\ 
	0.2 & 0.713 & 0.556 & 0.204 && 0.723 & 0.537 & 0.197 \\ 
	0.4 & 0.833 & 0.437 & 0.358 && 0.833 & 0.437 & 0.358 \\ 
	0.6 & 0.539 & 0.310 & 0.305 && 0.539 & 0.310 & 0.305 \\ 
	0.8 & 0.423 & 0.236 & 0.235 && 0.423 & 0.236 & 0.235 \\ 
	1 & 0.393 & 0.203 & 0.203 && 0.393 & 0.203 & 0.203 \\ 
		\hline
	\end{tabular}
\end{table}
As shown, the robustness of the test increases with $\alpha$, showing the robustness improvement of the proposed Wald-type test statistics. Note that the major difference between the p-value using clean data and data with outliers is obtained with the value $\alpha=0$ corresponding to the classical MLE. 
%\newpage

\subsection{First Word Data}
\begin{figure}[h!!]
\centering
	\includegraphics[scale=0.65]{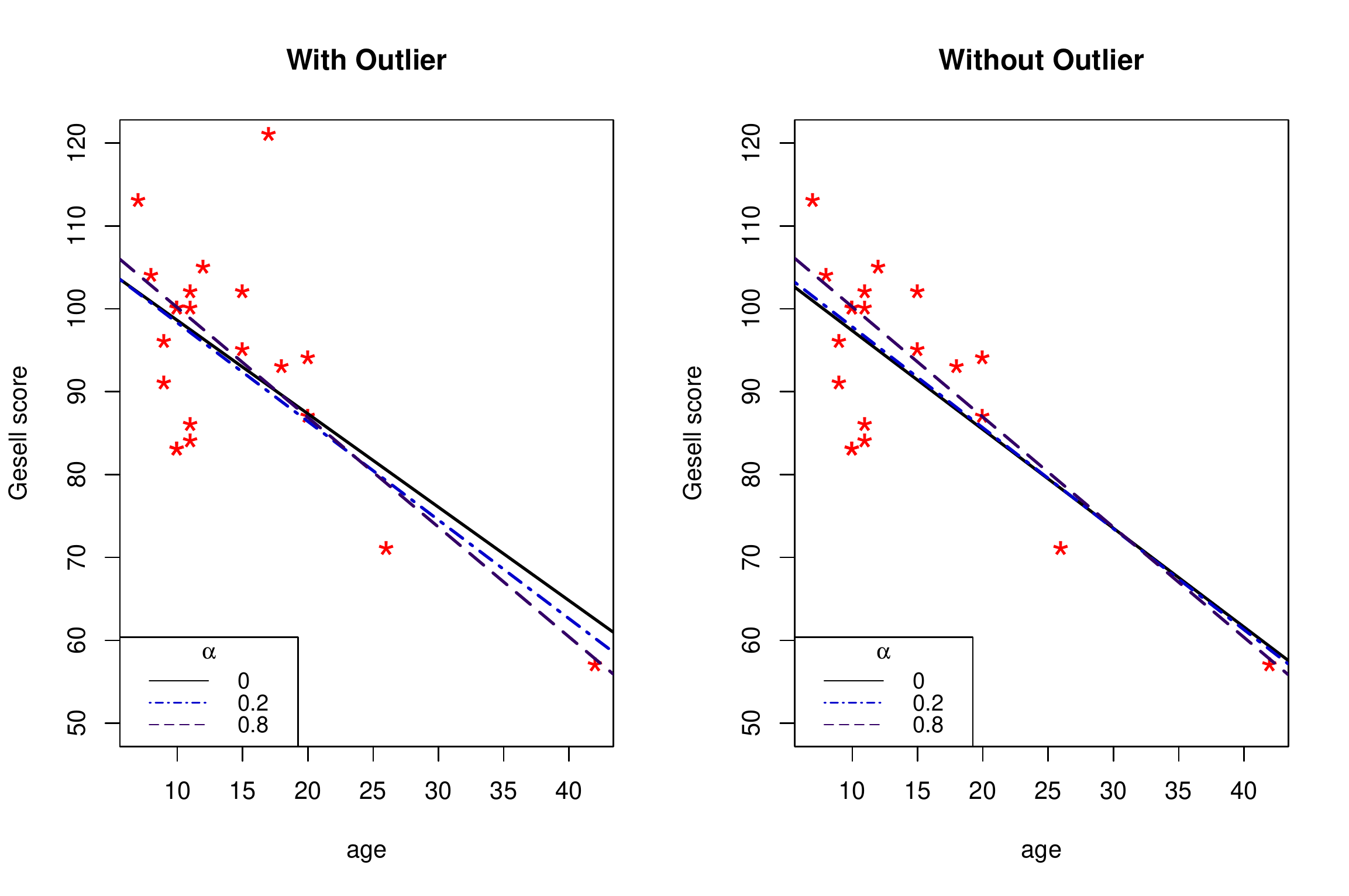}
	\caption{Plots of the data-points and fitted regression lines for the First Word Data using several minimum RP estimators before and after deleting the outliers.}
	\label{fig:words1}
\end{figure}

These data, originally presented in Mickey et al. \cite{Mickey}, consist on $n=21$ observations and relate the age in which children speak their first word to their Gesell adaptative score, a meassure of mental ability. By means of a sequential approach to detect outliers via stepwise regression,  \cite{Mickey} concluded that observation $18$ was an outlier.    While  estimates of the regression coefficients   obtained with the MLE  do not change excesively when omitting this outlier (Figure \ref{fig:words1}), we do observe a greater change in the error variance estimation (Table \ref{table:words1}). As expected, minimum RP estimates for $\alpha>0$ remain more robust. 

\begin{table}[h!]
	\center
	\caption{The parameter estimates of the linear regression model for the First Word Data using several minimum RP estimators.}
	\label{table:words1}
\begin{tabular}{r rrr r rrr}
\hline
& \multicolumn{3}{c}{With outliers} & &\multicolumn{3}{c}{Without outliers} \\
\cline{2-4}  \cline{6-8} 
$\alpha$ & $\sigma$  & $\beta_0$  & $\beta_1$  && $\sigma$  & $\beta_0$  & $\beta_1$  \\ \hline
0     & 10.4845 & 109.8730 & -1.1269 &  & 8.1976 & 109.2816 & -1.1916 \\
0.2   & 9.7860  & 110.2068 & -1.1897 &  & 8.5501 & 110.0225 & -1.2183 \\
0.4   & 9.2980  & 110.8118 & -1.2338 &  & 8.7780 & 110.8276 & -1.2451 \\
0.6   & 9.0319  & 111.7370 & -1.2710 &  & 8.8019 & 111.8168 & -1.2767 \\
0.8   & 8.3349  & 113.4011 & -1.3246 &  & 8.1972 & 113.5345 & -1.3292 \\
1     & 4.4187  & 116.6086 & -1.4065 &  & 4.4187 & 116.6086 & -1.4065 \\ \hline
\end{tabular}
\end{table}

As in the previous example, we consider the following tests
\begin{align}
\label{testword1} \operatorname{H}_0 &: \beta_0 = 112.56,\\ 
\label{testword2} \operatorname{H}_0 &: \beta_1 = -1.28,\\  
\label{testword3} \operatorname{H}_0 &: (\beta_0, \beta_1) = (112.56,-1.28) , 
\end{align}
where the values $112.56$ and $0.73$ are respectively the mean value of the estimated coefficients $\beta_0$ and $\beta_1$ with the different values of $\alpha$ and using the original data (with outliers) listed in Table \ref{table:words1}. Table \ref{table:pvalueword} shows the p-values obtained by using the corresponding Wald-type test statistics, $W_n(\widehat{\beta}_0)$, $W_n(\widehat{\beta}_1)$ and $W_n(\widehat{\boldsymbol{\beta}}).$

\begin{table}[h]
	\centering
	\label{table:pvalueword}
	\caption{p-value obtained for the tests (\ref{testword1})-(\ref{testword3}) using the corresponding Wald-type test statistics.}
	\begin{tabular}{r rrr r rrr}
		\hline
		& \multicolumn{3}{c}{With outliers} & &\multicolumn{3}{c}{Without outliers} \\
		\cline{2-4}  \cline{6-8} 
		$\alpha$ & $W_n(\widehat{\beta}_0)$ & $W_n(\widehat{\beta}_1)$ & $W_n(\widehat{\boldsymbol{\beta}})$ 
		&& $W_n(\widehat{\beta}_0)$ & $W_n(\widehat{\beta}_1)$ & $W_n(\widehat{\boldsymbol{\beta}})$ \\ 
		\hline
	0 & 0.072 & 0.098 & 0.071 && 0.013 & 0.293 & 0.001 \\ 
	0.2 & 0.110 & 0.331 & 0.065 && 0.065 & 0.485 & 0.007 \\ 
	0.4 & 0.243 & 0.636 & 0.098 && 0.234 & 0.719 & 0.061 \\ 
	0.6 & 0.596 & 0.949 & 0.318 && 0.628 & 0.997 & 0.311 \\ 
	0.8 & 0.590 & 0.620 & 0.588 && 0.529 & 0.583 & 0.529 \\ 
	1 & 0.001 & 0.078 & 0.000 && 0.001 & 0.078 & 0.000 \\ 
		\hline
	\end{tabular}
\end{table}
The results highlight again the gain in robustness.

\section{Concluding remarks}

In this paper we have presented the minimum RP estimators for the case of i.n.i.d.o.  Wald-type tests based on them are also developed. Classical MLE and Wald-test are obtained as a particular case of these new estimators and tests. In particular, we have studied the case of MLRM. Through the study of the influence functions and the development of an extensive simulation study we prove their robustness from a theoretical and practical point of view, respectively. Application to different models  is
a problem that will be of interest for further consideration.\\

\noindent \textbf{Acknowledgements:} This research is supported by the Spanish Grants no. PGC2018-095194-B-100, no. FPU19/01824  and no. FPU16/03104. 
%\clearpage
\appendix
\section{Proof of Results}

\subsection{Proof of Theorem \ref{Th1} \label{ap:con}}

The proof follows similar steps that the proof presented in \cite{Ghosh2013} for the minimum DPD estimators  for i.n.i.d.o and the proof presented in \cite{Leroy} for the MLE with i.n.i.d.o.

To prove the existence, with probability tending to $1$, of a consistent sequence of solutions of the system of equations (\ref{2.5}), we study the behaviour of the objective function in (\ref{2.4}), $H_{n}^{\alpha}(\boldsymbol{\theta}),$ on a neighbourhood of the true parameter value. We consider the sphere $Q_{a}$ with center at the true value of the parameter $\boldsymbol{\theta}^*$ and radius $a.$ We will show that for any sufficiently small $a$
\[
H_{n}^{\alpha}(\boldsymbol{\theta})<H_{n}^{\alpha}(\boldsymbol{\theta}%
_{0})\qquad
\]
with probability tending to $1$ at all points $\boldsymbol{\theta}$ on the surface of $Q_{a}.$ This inequality ensures  that the objective function $H_{n}^{\alpha}(\boldsymbol{\theta})$ has a local maximum in the interior of $Q_{a}$. Since $H_{n}^{\alpha}(\boldsymbol{\theta})$ is differentiable the system of equations (\ref{2.5}) must be satisfied at a local maximum. Therefore, for any $a>0$, the system of equations (\ref{2.5})
has a solution $\widehat{\boldsymbol{\theta}}_{n}(a)$ within $Q_{a}$ verifying
\[
\lim_{n\rightarrow\infty}P_{\boldsymbol{\theta}^*}\left(  \left\Vert
\widehat{\boldsymbol{\theta}}_{n}(a)-\boldsymbol{\theta}^*\right\Vert
\right)  =1.
\]
%To obtain the needed facts concerning the behaviour of $H_{n}^{\alpha}(\boldsymbol{\theta})$ on $Q_{a}$ for small $a$, 
We consider a Taylor series expansion of $H_{n}^{\alpha}(\boldsymbol{\theta})$ around $\boldsymbol{\theta
}_{0}=\left(  \theta_{0}^{1},...,\theta_{0}^{p}\right)  ,$%

\begin{align}
H_{n}^{\alpha}(\boldsymbol{\theta})-H_{n}^{\alpha}(\boldsymbol{\theta}^*) &
=%
%TCIMACRO{\tsum \limits_{j=1}^{p}}%
%BeginExpansion
{\textstyle\sum\limits_{j=1}^{p}}
%EndExpansion
\left(  \frac{\partial H_{n}^{\alpha}(\boldsymbol{\theta})}{\partial\theta
_{j}}\right)  _{\boldsymbol{\theta=\theta}_{0}}\left(  \theta_{j}-\theta
_{j}^{0}\right)  \label{2.6}\\
&  +\frac{1}{2}%
%TCIMACRO{\tsum \limits_{j=1}^{p}}%
%BeginExpansion
{\textstyle\sum\limits_{j=1}^{p}}
%EndExpansion%
%TCIMACRO{\tsum \limits_{k=1}^{p}}%
%BeginExpansion
{\textstyle\sum\limits_{k=1}^{p}}
%EndExpansion
\left(  \frac{\partial^{2}H_{n}^{\alpha}(\boldsymbol{\theta})}{\partial
\theta_{j}\partial\theta_{k}}\right)  _{\boldsymbol{\theta=\theta}_{0}}\left(
\theta_{j}-\theta_{j}^{0}\right)  \left(  \theta_{k}-\theta_{k}^{0}\right)
\nonumber\\
&  +\frac{1}{6}%
%TCIMACRO{\tsum \limits_{j=1}^{p}}%
%BeginExpansion
{\textstyle\sum\limits_{j=1}^{p}}
%EndExpansion%
%TCIMACRO{\tsum \limits_{k=1}^{p}}%
%BeginExpansion
{\textstyle\sum\limits_{k=1}^{p}}
%EndExpansion%
%TCIMACRO{\tsum \limits_{l=1}^{p}}%
%BeginExpansion
{\textstyle\sum\limits_{l=1}^{p}}
%EndExpansion
\left(  \frac{\partial^{3}H_{n}^{\alpha}(\boldsymbol{\theta})}{\partial
\theta_{j}\partial\theta_{k}\partial\theta_{l}}\right)  _{\boldsymbol{\theta
=\theta}^{\ast}}\left(  \theta_{j}-\theta_{j}^{0}\right)  \left(  \theta
_{k}-\theta_{k}^{0}\right)  \left(  \theta_{l}-\theta_{l}^{0}\right)
\nonumber\\
&  =%
%TCIMACRO{\tsum \limits_{j=1}^{p}}%
%BeginExpansion
{\textstyle\sum\limits_{j=1}^{p}}
%EndExpansion
\frac{1}{n}%
%TCIMACRO{\tsum \limits_{i=1}^{n}}%
%BeginExpansion
{\textstyle\sum\limits_{i=1}^{n}}
%EndExpansion
\left(  \frac{\partial V_{i}(Y_{i};\boldsymbol{\theta})}{\partial\theta_{j}%
}\right)  _{\boldsymbol{\theta=\theta}_{0}}\left(  \theta_{j}-\theta_{j}%
^{0}\right)  \nonumber\\
&  +\frac{1}{2}%
%TCIMACRO{\tsum \limits_{j=1}^{p}}%
%BeginExpansion
{\textstyle\sum\limits_{j=1}^{p}}
%EndExpansion%
%TCIMACRO{\tsum \limits_{k=1}^{p}}%
%BeginExpansion
{\textstyle\sum\limits_{k=1}^{p}}
%EndExpansion
\frac{1}{n}%
%TCIMACRO{\tsum \limits_{i=1}^{n}}%
%BeginExpansion
{\textstyle\sum\limits_{i=1}^{n}}
%EndExpansion
\left(  \frac{\partial^{2}V_{i}(Y_{i};\boldsymbol{\theta})}{\partial\theta
_{j}\partial\theta_{k}}\right)  _{\boldsymbol{\theta=\theta}_{0}}\left(
\theta_{j}-\theta_{j}^{0}\right)  \left(  \theta_{k}-\theta_{k}^{0}\right)
\nonumber\\
&  +\frac{1}{6}%
%TCIMACRO{\tsum \limits_{j=1}^{p}}%
%BeginExpansion
{\textstyle\sum\limits_{j=1}^{p}}
%EndExpansion%
%TCIMACRO{\tsum \limits_{k=1}^{p}}%
%BeginExpansion
{\textstyle\sum\limits_{k=1}^{p}}
%EndExpansion%
%TCIMACRO{\tsum \limits_{l=1}^{p}}%
%BeginExpansion
{\textstyle\sum\limits_{l=1}^{p}}
%EndExpansion
\frac{1}{n}%
%TCIMACRO{\tsum \limits_{i=1}^{n}}%
%BeginExpansion
{\textstyle\sum\limits_{i=1}^{n}}
%EndExpansion
\left(  \frac{\partial^{3}V_{i}(Y_{i};\boldsymbol{\theta})}{\partial\theta
_{j}\partial\theta_{k}\partial\theta_{l}}\right)  _{\boldsymbol{\theta=\theta
}^{\ast}}\left(  \theta_{j}-\theta_{j}^{0}\right)  \left(  \theta_{k}%
-\theta_{k}^{0}\right)  \left(  \theta_{l}-\theta_{l}^{0}\right)  \nonumber\\
&  =L_{1}+L_{2}+L_{3},\nonumber
\end{align}
where $\boldsymbol{\theta}^{\ast}$ belong to the interior of the ball centred on $\boldsymbol{\theta}^*$ and radius $a.$
We study separately right-hand terms $L_{1},L_{2}$ and $L_{3}$ in (\ref{2.6}). 

Using assumption \ref{itm:C6}, we have that
\[
A_{j}^{(n}=\left(  \frac{\partial H_{n}^{\alpha}(\boldsymbol{\theta}%
)}{\partial\theta_{j}}\right)  _{\boldsymbol{\theta=\theta}_{0}}=\frac{1}{n}%
%TCIMACRO{\tsum \limits_{i=1}^{n}}%
%BeginExpansion
{\textstyle\sum\limits_{i=1}^{n}}
%EndExpansion
\left(  \frac{\partial V_{i}(Y_{i};\boldsymbol{\theta})}{\partial\theta_{j}%
}\right)  _{\boldsymbol{\theta=\theta}_{0}}\overset{P}{\underset{}{\rightarrow
}}\frac{1}{n}%
%TCIMACRO{\tsum \limits_{i=1}^{n}}%
%BeginExpansion
{\textstyle\sum\limits_{i=1}^{n}}
%EndExpansion
E_{\boldsymbol{\theta}^*}\left[  \left(  \frac{\partial V_{i}%
(Y;\boldsymbol{\theta})}{\partial\theta_{j}}\right)  _{\boldsymbol{\theta
=\theta}_{0}}\right]  =0.
\]
We are going to establish the last equality,%
\[
\frac{\partial V_{i}(Y_{i};\boldsymbol{\theta})}{\partial\theta_{j}}=\frac
{1}{L_{\alpha}^{i}\left(  \boldsymbol{\theta}\right)  ^{2}}\left(  \alpha
f_{i}(Y,\boldsymbol{\theta})^{\alpha}u_{j}(Y,\boldsymbol{\theta})L_{\alpha
}^{i}\left(  \boldsymbol{\theta}\right)  -\frac{\partial L_{\alpha}^{i}\left(
\boldsymbol{\theta}\right)  }{\partial\theta_{j}}f_{i}(Y,\boldsymbol{\theta
})^{\alpha}\right)  .
\]
with $u_j(y,\boldsymbol{\theta})=\frac{\partial \log (f_i(Y,\boldsymbol{\theta}))}{\partial {\theta_j}}$. But
\begin{align*}
\frac{\partial L_{\alpha}^{i}\left(  \boldsymbol{\theta}\right)  }%
{\partial\theta_{j}}  &  =\frac{\alpha}{\alpha+1}\left(  \int f_{i}%
(y,\boldsymbol{\theta})^{\alpha+1}dy\right)  ^{\frac{\alpha}{\alpha+1}%
-1}\left(  \alpha+1\right)  \int f_{i}(y,\boldsymbol{\theta})^{\alpha+1}%
u_{j}(y,\boldsymbol{\theta})dy\\
&  =\alpha\left(  \int f_{i}(y,\boldsymbol{\theta})^{\alpha+1}dy\right)
^{\frac{\alpha}{\alpha+1}-1}\int f_{i}(y,\boldsymbol{\theta})^{\alpha+1}%
u_{j}(y,\boldsymbol{\theta})dy.
\end{align*}
Therefore,%
\begin{align*}
L_{\alpha}^{i}\left(  \boldsymbol{\theta}\right)  ^{2}\frac{\partial
V_{i}(Y;\boldsymbol{\theta})}{\partial\theta_{j}}  &  =\alpha f_{i}%
(Y,\boldsymbol{\theta})^{\alpha}u_{j}(Y,\boldsymbol{\theta})L_{\alpha}%
^{i}\left(  \boldsymbol{\theta}\right) \\
&  \left[  \alpha\left(  \int f_{i}(y,\boldsymbol{\theta})^{\alpha
+1}dy\right)  ^{\frac{\alpha}{\alpha+1}-1}\int f_{i}(y,\boldsymbol{\theta
})^{\alpha+1}u_{j}(y,\boldsymbol{\theta})dy\right]  f_{i}(Y,\boldsymbol{\theta
})^{\alpha}%
\end{align*}
and
\begin{small}
\begin{align*}
E_{\boldsymbol{\theta}^*}\left[  \left(  \frac{\partial V_{i}(Y_{i}%
;\boldsymbol{\theta})}{\partial\theta_{j}}\right)  _{\boldsymbol{\theta
=\theta}_{0}}\right]   &  =%
%TCIMACRO{\dint }%
%BeginExpansion
{\displaystyle\int}
%EndExpansion
\left(  \alpha f_{i}(y,\boldsymbol{\theta}^*)^{\alpha}u_{j}%
(y,\boldsymbol{\theta}^*)L_{\alpha}^{i}\left(  \boldsymbol{\theta}%
_{0}\right)  \right)  f_{i}(y,\boldsymbol{\theta}^*)dy\\
&  -%
%TCIMACRO{\dint }%
%BeginExpansion
{\displaystyle\int}
%EndExpansion
\left[  \alpha\left(  \int f_{i}(y,\boldsymbol{\theta}^*)^{\alpha
+1}dy\right)  ^{\frac{\alpha}{\alpha+1}-1}\int f_{i}(y,\boldsymbol{\theta}%
_{0})^{\alpha+1}u_{j}(y,\boldsymbol{\theta}^*)dy\right]  f_{i}%
(y,\boldsymbol{\theta}^*)^{\alpha+1}dy\\
&  =0.
\end{align*}
\end{small}
On the other hand, we denote
\[
B_{jk}^{(n}=\left(  \frac{\partial^{2}H_{n}^{\alpha}(\boldsymbol{\theta}%
)}{\partial\theta_{j}\partial\theta_{k}}\right)  _{\boldsymbol{\theta=\theta
}_{0}},
\]
and applying again condition \ref{itm:C6}, we obtain the convergence
\begin{small}
\[
B_{jk}^{(n}=\left(  \frac{\partial^{2}H_{n}^{\alpha}(\boldsymbol{\theta}%
)}{\partial\theta_{j}\partial\theta_{k}}\right)  _{\boldsymbol{\theta=\theta
}_{0}}=\frac{1}{n}%
%TCIMACRO{\tsum \limits_{i=1}^{n}}%
%BeginExpansion
{\textstyle\sum\limits_{i=1}^{n}}
%EndExpansion
\left(  \frac{\partial^{2}V_{i}(Y_{i};\boldsymbol{\theta})}{\partial\theta
_{j}\partial\theta_{k}}\right)  _{\boldsymbol{\theta=\theta}_{0}%
}\overset{P}{\underset{}{\rightarrow}}\frac{1}{n}%
%TCIMACRO{\tsum \limits_{i=1}^{n}}%
%BeginExpansion
{\textstyle\sum\limits_{i=1}^{n}}
%EndExpansion
E_{\boldsymbol{\theta}^*}\left[  \left(  \frac{\partial^{2}V_{i}%
(Y;\boldsymbol{\theta})}{\partial\theta_{j}\partial\theta_{k}}\right)
_{\boldsymbol{\theta=\theta}_{0}}\right]  =\left(  -\boldsymbol{\Psi}_{n}\right)  _{jk}.
\]
\end{small}

Finally, applying condition \ref{itm:C6} to the third derivative, we have%
\[
\left(  \frac{\partial^{3}H_{n}^{\alpha}(\boldsymbol{\theta})}{\partial
\theta_{j}\partial\theta_{k}\partial\theta_{l}}\right)  _{\boldsymbol{\theta
=\theta}^{\ast}}=\frac{1}{n}%
%TCIMACRO{\tsum \limits_{i=1}^{n}}%
%BeginExpansion
{\textstyle\sum\limits_{i=1}^{n}}
%EndExpansion
\left(  \frac{\partial^{3}V_{i}(Y_{i};\boldsymbol{\theta})}{\partial\theta
_{j}\partial\theta_{k}\partial\theta_{l}}\right)  _{\boldsymbol{\theta=\theta
}^{\ast}}\overset{P}{\underset{}{\rightarrow}}\frac{1}{n}%
%TCIMACRO{\tsum \limits_{i=1}^{n}}%
%BeginExpansion
{\textstyle\sum\limits_{i=1}^{n}}
%EndExpansion
E_{\boldsymbol{\theta}_{\ast}}\left[  \left(  \frac{\partial^{3}%
V_{i}(Y;\boldsymbol{\theta})}{\partial\theta_{j}\partial\theta_{k}%
\partial\theta_{l}}\right)  _{\boldsymbol{\theta=\theta}^{\ast}}\right]  .
\]
Assumption \ref{itm:C5} ensures the existence of $ M_{jkl}, j,k,l= 1,..,p$ s.t.
\[
\left\vert \left(  \frac{\partial^{3}V_{i}(Y_{i};\boldsymbol{\theta}%
)}{\partial\theta_{j}\partial\theta_{k}\partial\theta_{l}}\right)
_{\boldsymbol{\theta=\theta}^{\ast}}\right\vert \leq M_{jkl}^{\left(
i\right)  }\left(  y\right),
\]
and therefore there exists $\gamma_{jkl}^{(i)}(y)$ verifying%
\begin{equation}
0\leq\left\vert \gamma_{jkl}(y)\right\vert \leq1 \label{2.7}%
\end{equation}
in such a way that
\[
\left(  \frac{\partial^{3}V_{i}(Y_{i};\boldsymbol{\theta})}{\partial\theta
_{j}\partial\theta_{k}\partial\theta_{l}}\right)  _{\boldsymbol{\theta=\theta
}^{\ast}}=M_{jkl}^{\left(  i\right)  }\left(  Y_{i}\right)  \gamma_{jkl}^{(i)}%
(Y_{i}).
\]
and
\[
E_{\boldsymbol{\theta}_{\ast}}\left[  M_{jkl}^{\left(  i\right)  }\left(
Y_{i}\right)  \right]  =m_{jkl},
\]
with
\[
\left\vert \frac{1}{n}%
%TCIMACRO{\tsum \limits_{i=1}^{n}}%
%BeginExpansion
{\textstyle\sum\limits_{i=1}^{n}}
%EndExpansion
E_{\boldsymbol{\theta}_{\ast}}\left[  \left(  \frac{\partial^{3}%
V_{i}(Y;\boldsymbol{\theta})}{\partial\theta_{j}\partial\theta_{k}%
\partial\theta_{l}}\right)  _{\boldsymbol{\theta=\theta}^{\ast}}\right]
\right\vert <m_{jkl}.
\]
The previous convergence provide that for all $a$ and  for all $\varepsilon$ there exists $n_{0}$ such that for all
$n>n_{0}$ we have
\[%
\begin{array}
[c]{l}%
P\left(  \left\vert A_{j}^{(n}\right\vert >a^{2}\right)  <\frac{\varepsilon
}{p+p^{2}+p^{3}}\\
P\left(  \left\vert B_{jk}^{(n}-\left(  \Psi_{n}\right)  _{jk}\right\vert \geq
a\right)  <\frac{\varepsilon}{p+p^{2}+p^{3}}\\
P\left(  \left\vert \frac{1}{n}%
%TCIMACRO{\tsum \limits_{i=1}^{n}}%
%BeginExpansion
{\textstyle\sum\limits_{i=1}^{n}}
%EndExpansion
\left(  \frac{\partial^{3}V_{i}(Y_{i};\boldsymbol{\theta})}{\partial\theta
_{j}\partial\theta_{k}\partial\theta_{l}}\right)  _{\boldsymbol{\theta=\theta
}^{\ast}}\right\vert \geq2m_{jkl}\right)  \leq\frac{\varepsilon}%
{p+p^{2}+p^{3}}.
\end{array}
\]
We shall denote now by $S^{\ast}$ the event containing the $p+p^{2}+p^{3}$
inequalities,%
\[
\left\{  \left\vert A_{j}^{(n}\right\vert >a^{2};\text{ }\left\vert
B_{jk}^{(n}-\left(  -\Psi_{n}\right)  _{jk}\right\vert \geq a;\text{
}\left\vert \frac{1}{n}%
%TCIMACRO{\tsum \limits_{i=1}^{n}}%
%BeginExpansion
{\textstyle\sum\limits_{i=1}^{n}}
%EndExpansion
\left(  \frac{\partial^{3}V_{i}(Y_{i};\boldsymbol{\theta})}{\partial\theta
_{j}\partial\theta_{k}\partial\theta_{l}}\right)  _{\boldsymbol{\theta=\theta
}^{\bullet}}\right\vert \geq2m_{jkl}\right\}  .
\]
It is clear that $P\left(  S^{\ast}\right)  <\varepsilon$ and $P(\left(
S^{\ast}\right)  ^{C})\geq1-\varepsilon.$ In the following we denote
$S=\left(  S^{\ast}\right)  ^{C}.$ We finally study the sign of $H_{n}^{\alpha}(\boldsymbol{\theta})-H_{n}^{\alpha
}(\boldsymbol{\theta}^*)$ under the event $S$ and for $\boldsymbol{\theta\in}Q_{a}.$

Since $\boldsymbol{\theta\in}Q_{a}$ in $S$ holds
\begin{equation}
\left\vert L_{1}\right\vert =\left\vert \frac{1}{n}%
%TCIMACRO{\tsum \limits_{j=1}^{p}}%
%BeginExpansion
{\textstyle\sum\limits_{j=1}^{p}}
%EndExpansion%
%TCIMACRO{\tsum \limits_{i=1}^{n}}%
%BeginExpansion
{\textstyle\sum\limits_{i=1}^{n}}
%EndExpansion
\left(  \frac{\partial V_{i}(Y_{i};\boldsymbol{\theta})}{\partial\theta_{j}%
}\right)  _{\boldsymbol{\theta=\theta}_{0}}\left(  \theta_{j}-\theta_{j}%
^{0}\right)  \right\vert \leq p\text{ }a\text{ }a^{2}\label{L1}%
\end{equation}
and
\[
\left\vert \frac{1}{2}%
%TCIMACRO{\tsum \limits_{j=1}^{p}}%
%BeginExpansion
{\textstyle\sum\limits_{j=1}^{p}}
%EndExpansion%
%TCIMACRO{\tsum \limits_{k=1}^{p}}%
%BeginExpansion
{\textstyle\sum\limits_{k=1}^{p}}
%EndExpansion
\left\{  \left(  B_{jk}^{(n}-\left(  -\Psi_{n}\right)  _{jk}\right)  \right\}
\left(  \theta_{j}-\theta_{j}^{0}\right)  \left(  \theta_{k}-\theta_{k}%
^{0}\right)  \right\vert \leq\frac{1}{2}p^{2}\text{ }a^{2}\text{ }a.
\]
%In equation (\ref{L1}) we have used that in $S,$%
%\[
%\left\vert A_{j}^{(n}\right\vert =\left\vert \frac{1}{n}%
%%TCIMACRO{\tsum \limits_{i=1}^{n}}%
%%BeginExpansion
%{\textstyle\sum\limits_{i=1}^{n}}
%%EndExpansion
%\left(  \frac{\partial V_{i}(Y_{i};\boldsymbol{\theta})}{\partial\theta_{j}%
%}\right)  _{\boldsymbol{\theta=\theta}_{0}}\right\vert \leq a^{2}\text{ and
%}\left\vert \theta_{j}-\theta_{j}^{0}\right\vert <a.
%\]

We now consider the negative quadratic form%
\[
\boldsymbol{A}=-\frac{1}{2}%
%TCIMACRO{\tsum \limits_{j=1}^{p}}%
%BeginExpansion
{\textstyle\sum\limits_{j=1}^{p}}
%EndExpansion%
%TCIMACRO{\tsum \limits_{k=1}^{p}}%
%BeginExpansion
{\textstyle\sum\limits_{k=1}^{p}}
%EndExpansion
\left(  \theta_{j}-\theta_{j}^{0}\right)  \left(  \theta_{k}-\theta_{k}%
^{0}\right)  \frac{1}{n}%
%TCIMACRO{\tsum \limits_{i=1}^{n}}%
%BeginExpansion
{\textstyle\sum\limits_{i=1}^{n}}
%EndExpansion
E_{\boldsymbol{\theta}^*}\left[  \left(  \frac{\partial^{2}V_{i}%
(Y;\boldsymbol{\theta})}{\partial\theta_{j}\partial\theta_{k}}\right)
_{\boldsymbol{\theta=\theta}_{0}}\right]  .
\]
 An orthogonal transformation can reduce the quadratic form $\boldsymbol{A}$ to its diagonal form $ \boldsymbol{A}= 
%TCIMACRO{\tsum \limits_{i=1}^{p}}%
%BeginExpansion
{\textstyle\sum\limits_{i=1}^{p}}
%EndExpansion
\lambda_{i}\xi_{i}^{2}$ with $%
%TCIMACRO{\tsum \limits_{i=1}^{p}}%
%BeginExpansion
{\textstyle\sum\limits_{i=1}^{p}}
%EndExpansion
\xi_{i}^{2}=%
%TCIMACRO{\tsum \limits_{i=1}^{p}}%
%BeginExpansion
{\textstyle\sum\limits_{i=1}^{p}}
%EndExpansion
\left(  \theta_{i}-\theta_{i}^{p}\right)  ^{2}=a^{2}.$ Shorting the negatives eigenvalues
$\lambda_{i}$  we get
\[%
%TCIMACRO{\tsum \limits_{i=1}^{p}}%
%BeginExpansion
{\textstyle\sum\limits_{i=1}^{p}}
%EndExpansion
\lambda_{i}\xi_{i}^{2}\leq-\lambda_{0}a^{2}<0.
\]
A study of the sign of the function
$
\frac{1}{2}p^{2}a^{3}-\lambda_{0}a^{2}%
$
proves that we can find $c>0,a_{0}>0$ so that for $a<a_{0}$%
\begin{align*}
\left\vert L_{2}\right\vert =&\left\vert \frac{1}{2}
{\textstyle\sum\limits_{j=1}^{p}}{\textstyle\sum\limits_{k=1}^{p}}
\left\{  \left(  B_{jk}^{(n}-\left(  -\Psi_{n}\right)  _{jk}\right)  \right\}
\left(  \theta_{j}-\theta_{j}^{p}\right)  \left(  \theta_{k}-\theta_{k}%
^{p}\right) \right. \\
&\left.  +\frac{1}{2}{\textstyle\sum\limits_{j=1}^{p}}{\textstyle\sum\limits_{k=1}^{p}}
\left(  \theta_{j}-\theta_{j}^{p}\right)  \left(  \theta_{k}-\theta_{k}%
^{p}\right)  \left(  -\Psi_{n}\right)  _{jk}\right\vert <-ca^{2}.
\end{align*}

Lastly,%
\begin{align*}
\left\vert L_{3}\right\vert  &  =\left\vert \frac{1}{6}%
%TCIMACRO{\tsum \limits_{j=1}^{p}}%
%BeginExpansion
{\textstyle\sum\limits_{j=1}^{p}}
%EndExpansion%
%TCIMACRO{\tsum \limits_{k=1}^{p}}%
%BeginExpansion
{\textstyle\sum\limits_{k=1}^{p}}
%EndExpansion%
%TCIMACRO{\tsum \limits_{l=1}^{p}}%
%BeginExpansion
{\textstyle\sum\limits_{l=1}^{p}}
%EndExpansion
\frac{1}{n}%
%TCIMACRO{\tsum \limits_{i=1}^{n}}%
%BeginExpansion
{\textstyle\sum\limits_{i=1}^{n}}
%EndExpansion
\left(  \frac{\partial^{3}V_{i}(Y_{i};\boldsymbol{\theta})}{\partial\theta
_{j}\partial\theta_{k}\partial\theta_{l}}\right)  _{\boldsymbol{\theta=\theta
}^{\bullet}}\left(  \theta_{j}-\theta_{j}^{p}\right)  \left(  \theta
_{k}-\theta_{k}^{p}\right)  \left(  \theta_{l}-\theta_{l}^{p}\right)
\right\vert \\
&  <\frac{2}{6}%
%TCIMACRO{\tsum \limits_{j=1}^{p}}%
%BeginExpansion
{\textstyle\sum\limits_{j=1}^{p}}
%EndExpansion%
%TCIMACRO{\tsum \limits_{k=1}^{p}}%
%BeginExpansion
{\textstyle\sum\limits_{k=1}^{p}}
%EndExpansion%
%TCIMACRO{\tsum \limits_{l=1}^{p}}%
%BeginExpansion
{\textstyle\sum\limits_{l=1}^{p}}
%EndExpansion
m_{jkl}a^{3}=a^{3}b,
\end{align*}
being
\[
b=\frac{2}{6}%
%TCIMACRO{\tsum \limits_{j=1}^{p}}%
%BeginExpansion
{\textstyle\sum\limits_{j=1}^{p}}
%EndExpansion%
%TCIMACRO{\tsum \limits_{k=1}^{p}}%
%BeginExpansion
{\textstyle\sum\limits_{k=1}^{p}}
%EndExpansion%
%TCIMACRO{\tsum \limits_{l=1}^{p}}%
%BeginExpansion
{\textstyle\sum\limits_{l=1}^{p}}
%EndExpansion
m_{jkl}.
\]
Therefore,%
\[
H_{n}^{\alpha}(\boldsymbol{\theta})-H_{n}^{\alpha}(\boldsymbol{\theta}%
_{0})<pa^{3}-ca^{2}+ba^{3}%
\]
and $pa^{3}-ca^{2}+ba^{3}<0$ if and only if $a<\frac{c}{b+p}.$ Therefore
assuming $a<\frac{c}{b+p}$ we get that in the event $S$%
\[
H_{n}^{\alpha}(\boldsymbol{\theta})-H_{n}^{\alpha}(\boldsymbol{\theta}%
_{0})<0\text{ }\forall\boldsymbol{\theta}\in Q_{a}.
\]
Thus the event $C$ involving all $\boldsymbol{\theta}\in Q_{a}$ s.t.
$H_{n}^{\alpha}(\boldsymbol{\theta})-H_{n}^{\alpha}(\boldsymbol{\theta}^*<0$, is contained in $S$, i.e. $P\left(  C\right)  \geq P(S)>1-\varepsilon.$
Choosing $\ a$ lower than $\min\left(  a_{0},\frac{c}{b+p}\right)  ,$ we have
\[
\lim_{n\rightarrow\infty}P\left(  \forall\boldsymbol{\theta}\in Q_{a}%
/H_{n}^{\alpha}(\boldsymbol{\theta})-H_{n}^{\alpha}(\boldsymbol{\theta}%
_{0})<0\right)  =1.
\]
Thus, there exists $\widehat{\boldsymbol{\theta}}_{n}(a)$ $\in Q_{a}$,
i,e., $\left\Vert \widehat{\boldsymbol{\theta}}_{n}(a)-\boldsymbol{\theta}%
_{0}\right\Vert <a$ such that $H_{n}^{\alpha}(\boldsymbol{\theta})$ has a
local maximum in $\widehat{\boldsymbol{\theta}}_{n}(a),$ i.e.,%
\[
\forall\text{ }a\leq\min\left(  a_{0},\frac{c}{b+p}\right)  ,
\]
we obtain the required convergence
\[
\lim_{n\rightarrow\infty}P\left(  \left\Vert \widehat{\boldsymbol{\theta}}%
_{n}(a)-\boldsymbol{\theta}^*\right\Vert <a\right)  =1.
\]

\bigskip
\subsection{Proof of Theorem \ref{Th2} \label{ap:asy}}

We denote
\[
H_{n,j}^{\alpha}(\boldsymbol{\theta})=\frac{\partial H_{n}^{\alpha
}(\boldsymbol{\theta})}{\partial\theta_{j}}%
\]
with $H_{n}^{\alpha}(\boldsymbol{\theta})$ defined in (\ref{2.4}). A Taylor
expansion of $H_{n,j}^{\alpha}(\boldsymbol{\theta})$ around
$\boldsymbol{\theta}^*,$ gives,%
\begin{align*}
H_{n,j}^{\alpha}(\boldsymbol{\theta})=&H_{n,j}^{\alpha}(\boldsymbol{\theta}^*)+{\textstyle\sum\limits_{k=1}^{p}}\left(  \frac{\partial^{2}H_{n}^{\alpha}(\boldsymbol{\theta})}{\partial\theta_{j}\partial\theta_{k}}\right)  _{\boldsymbol{\theta=\theta}_{0}}\left(\theta_{k}-\theta_{k}^{0}\right) \\
& +\frac{1}{2}{\textstyle\sum\limits_{k=1}^{p}}{\textstyle\sum\limits_{l=1}^{p}}\left(  \frac{\partial^{3}H_{n}^{\alpha}(\boldsymbol{\theta})}{\partial\theta_{j}\partial\theta_{k}\partial\theta_{l}}\right)  _{\boldsymbol{\theta=\theta}^{\ast}}\left(  \theta_{k-}\theta_{k}^{0}\right)  \left(  \theta_{l}-\theta_{l}^{0}\right)
\end{align*}
with $\boldsymbol{\theta}^{\ast}$ in the segment conecting $\boldsymbol{\theta
}$ and $\boldsymbol{\theta}^*.$ It is clear that at the minimum RP estimator the function $H_{n,j}^{\alpha}$ vanishes, $H_{n,j}^{\alpha}(\widehat{\boldsymbol{\theta}}_{\alpha})=0.$
Therefore,
\begin{align*}
H_{n,j}^{\alpha}(\boldsymbol{\theta}^*)=&-{\textstyle\sum\limits_{k=1}^{p}}\left(  \frac{\partial^{2}H_{n}^{\alpha}(\boldsymbol{\theta})}{\partial
\theta_{j}\partial\theta_{k}}\right)  _{\boldsymbol{\theta=\theta}_{0}}\left(
\widehat{\theta}_{\alpha,k}-\theta_{k}^{0}\right)\\
&  -\frac{1}{2}{\textstyle\sum\limits_{k=1}^{p}}{\textstyle\sum\limits_{l=1}^{p}}\left(  \frac{\partial^{3}H_{n}^{\alpha}(\boldsymbol{\theta})}{\partial
\theta_{j}\partial\theta_{k}\partial\theta_{l}}\right)  _{\boldsymbol{\theta
=\theta}^{\ast}}\left(  \widehat{\theta}_{\alpha,k}-\theta_{k}^{0}\right)
\left(  \widehat{\theta}_{\alpha,l}-\theta_{l}^{0}\right)  .
\end{align*}
Using that
\[
H_{n,j}^{\alpha}(\boldsymbol{\theta}^*)=\frac{1}{n}%
%TCIMACRO{\tsum \limits_{i=1}^{n}}%
%BeginExpansion
{\textstyle\sum\limits_{i=1}^{n}}
%EndExpansion
\left(  \frac{\partial V_{i}(Y,\boldsymbol{\theta)}}{\partial\theta_{j}%
}\right)  _{\boldsymbol{\theta=\theta}_{0}},
\]
it holds
\begin{align*}
	\frac{1}{\sqrt{n}}{\textstyle\sum\limits_{i=1}^{n}}
	\left(  \frac{\partial V_{i}(Y,\boldsymbol{\theta)}}{\partial\theta_{j}}\right)  _{\boldsymbol{\theta=\theta}_{0}}=&\sqrt{n}{\sum_{k=1}^{p}}
	\left(\widehat{\theta}_{\alpha,k}-\theta_{k}^{0}\right)  
	\bigg\{ -\left(	\frac{\partial^{2}H_{n}^{\alpha}(\boldsymbol{\theta})}{\partial\theta_{j}\partial\theta_{k}}\right)_{\boldsymbol{\theta=\theta}_{0}}\\
	&-\frac{1}{2} \sum_{k=1}^{p}\sum_{l=1}^{p}
	\left( \frac{\partial^{3}H_{n}^{\alpha}(\boldsymbol{\theta})}{\partial \theta_{j}\partial\theta_{k}\partial\theta_{l}}\right)  _{\boldsymbol{\theta=\theta}^{\ast}}\left(  \widehat{\theta}_{\alpha,l}-\theta_{l}^{0}\right)
	\bigg\}.
\end{align*}
If we denote,
\begin{align*}
Z_{kn} &  =\sqrt{n}%
%TCIMACRO{\tsum \limits_{k=1}^{p}}%
%BeginExpansion
{\textstyle\sum\limits_{k=1}^{p}}
%EndExpansion
\left(  \widehat{\theta}_{\alpha,k}-\theta_{k}^{0}\right),  \\
A_{jkn} &  =-\left(  \frac{\partial^{2}H_{n}^{\alpha}(\boldsymbol{\theta}%
)}{\partial\theta_{j}\partial\theta_{k}}\right)_{\boldsymbol{\theta=\theta
}_{0}}-\frac{1}{2}%
%TCIMACRO{\tsum \limits_{k=1}^{p}}%
%BeginExpansion
{\textstyle\sum\limits_{k=1}^{p}}
%EndExpansion%
%TCIMACRO{\tsum \limits_{l=1}^{p}}%
%BeginExpansion
{\textstyle\sum\limits_{l=1}^{p}}
%EndExpansion
\left(  \frac{\partial^{3}H_{n}^{\alpha}(\boldsymbol{\theta})}{\partial
\theta_{j}\partial\theta_{k}\partial\theta_{l}}\right)  _{\boldsymbol{\theta
=\theta}^{\ast}}\left(  \widehat{\theta}_{\alpha,l}-\theta_{l}^{0}\right),  \\
T_{jn} &  =\frac{1}{\sqrt{n}}%
%TCIMACRO{\tsum \limits_{i=1}^{n}}%
%BeginExpansion
{\textstyle\sum\limits_{i=1}^{n}}
%EndExpansion
\left(  \frac{\partial V_{i}(Y,\boldsymbol{\theta)}}{\partial
\boldsymbol{\theta}}\right)  _{\boldsymbol{\theta=\theta}_{0}},
\end{align*}
we can write
\[
T_{jn}=%
%TCIMACRO{\tsum \limits_{k=1}^{p}}%
%BeginExpansion
{\textstyle\sum\limits_{k=1}^{p}}
%EndExpansion
A_{jkn}Z_{kn}.
\]

Finally, we define the following vectors
\begin{align*}
\boldsymbol{Z}_{n} &  =\left(  Z_{1n},...,Z_{pn}\right)  ^{T},\\
\boldsymbol{T}_{n} &  =\left(  T_{1n},...,T_{pn}\right)  ^{T},\\
\boldsymbol{A}_{n} &  =(A_{jkn})_{j=1,...,p;k=1,...,p}.
\end{align*}
It is clear that
\begin{align*}
\boldsymbol{T}_{n} & =\boldsymbol{A}_{n}\boldsymbol{Z}_{n}\\
 &=\left(  \frac{1}{\sqrt{n}}%
%TCIMACRO{\tsum \limits_{i=1}^{n}}%
%BeginExpansion
{\textstyle\sum\limits_{i=1}^{n}}
%EndExpansion
\left(  \frac{\partial V_{i}(Y,\boldsymbol{\theta)}}{\partial\theta_{1}%
}\right)  _{\boldsymbol{\theta=\theta}_{0}},...,\frac{1}{\sqrt{n}}%
%TCIMACRO{\tsum \limits_{i=1}^{n}}%
%BeginExpansion
{\textstyle\sum\limits_{i=1}^{n}}
%EndExpansion
\left(  \frac{\partial V_{i}(Y,\boldsymbol{\theta)}}{\partial\theta_{p}%
}\right)  _{\boldsymbol{\theta=\theta}_{0}}\right)  ^{T}\\
&  =\frac{1}{\sqrt{n}}%
%TCIMACRO{\tsum \limits_{i=1}^{n}}%
%BeginExpansion
{\textstyle\sum\limits_{i=1}^{n}}
%EndExpansion
\left(  \frac{\partial V_{i}(Y,\boldsymbol{\theta)}}{\partial
\boldsymbol{\theta}}\right)  _{\boldsymbol{\theta=\theta}_{0}},
\end{align*}
and it is a simple exercise to verify that $V_{i}(Y,\boldsymbol{\theta)}$, $i=1,\dots,n$, are
independent with
\[
E_{\boldsymbol{\theta}^*}\left[  \left(  \frac{\partial V_{i}%
(Y,\boldsymbol{\theta)}}{\partial\boldsymbol{\theta}}\right)
_{\boldsymbol{\theta=\theta}_{0}}\right]  =0
\]
and
\[
Var_{\boldsymbol{\theta}^*}\left[  \left(  \frac{\partial V_{i}%
(Y,\boldsymbol{\theta)}}{\partial\boldsymbol{\theta}}\right)
_{\boldsymbol{\theta=\theta}_{0}}\right]  <\infty.
\]
By Assumption \ref{itm:C7} and applying the multivariate extension of Lindeberg-Levy central limit theorem we get
\[
\sqrt{n}\boldsymbol{\Omega}_{n}^{-\frac{1}{2}}\boldsymbol{T}_{n}%
\underset{n\rightarrow\infty}{\overset{L}{\rightarrow}}N(\boldsymbol{0}%
_{p},\boldsymbol{I}_{p})
\]
or equivalently%
\[
\sqrt{n}\boldsymbol{\Omega}_{n}^{-\frac{1}{2}}\boldsymbol{A}_{n}%
\boldsymbol{Z}_{n}\underset{n\rightarrow\infty}{\overset{L}{\rightarrow}%
}N(\boldsymbol{0}_{p},\boldsymbol{I}_{p}).
\]
By assumption \ref{itm:C5},
\[
\left(  \frac{\partial^{3}H_{n}^{\alpha}(\boldsymbol{\theta})}{\partial
\theta_{j}\partial\theta_{k}\partial\theta_{l}}\right)  _{\boldsymbol{\theta
=\theta}^{\ast}}%
\]
is bounded with probability tending to one. Therefore based on the consistency
of $\widehat{\boldsymbol{\theta}}_{\alpha}$ we have that the second term of
$A_{jkn}$ converges to zero in probability.
Moreover
\[
-\left(  \frac{\partial^{2}H_{n}^{\alpha}(\boldsymbol{\theta})}{\partial
\theta_{j}\partial\theta_{k}}\right)  _{\boldsymbol{\theta=\theta}_{0}%
}\underset{n\rightarrow\infty}{\overset{P}{\rightarrow}}\left(  \Psi
_{n}\right)  _{jk}%
\]
and hence%
\[
\boldsymbol{\Omega}_{n}^{-\frac{1}{2}}\left(  \boldsymbol{A}_{n}%
-\boldsymbol{\Psi}_{n}\right)  \boldsymbol{Z}_{n}\underset{n\rightarrow
\infty}{\overset{P}{\rightarrow}}\boldsymbol{0}_{p}.
\]
Therefore,%
\[
\boldsymbol{\Omega}_{n}^{-\frac{1}{2}}\boldsymbol{\Psi}_{n}\boldsymbol{Z}%
_{n}\underset{n\rightarrow\infty}{\overset{L}{\rightarrow}}N(\boldsymbol{0}%
_{p},\boldsymbol{I}_{p})
\]
and finally%
\[
\boldsymbol{\Omega}_{n}^{-\frac{1}{2}}\boldsymbol{\Psi}_{n}\left(
\widehat{\boldsymbol{\theta}}_{\alpha}-\boldsymbol{\theta}^*\right)
\underset{n\rightarrow\infty}{\overset{L}{\rightarrow}}N(\boldsymbol{0}%
_{p},\boldsymbol{I}_{p}).
\]

\bigskip

\subsection{Proof of Theorem \ref{th:asym_test}}
We have by (\ref{eq:asym_est}) that
\begin{align*}
\sqrt{n}(\widehat{\boldsymbol{\theta}}_{\alpha}-\boldsymbol{\theta}^*)\underset{n\rightarrow\infty}{\overset{L}{\rightarrow}}N(\boldsymbol{0}_{p},\boldsymbol{\Sigma}_{\alpha}(\boldsymbol{\theta}^*)),
\end{align*}
where $\boldsymbol{\Sigma}_{\alpha}(\boldsymbol{\theta}^*)=\lim_{n\rightarrow\infty} \boldsymbol{\Psi}_n(\boldsymbol{\theta}^*) \boldsymbol{\Omega}_n^{-1}(\boldsymbol{\theta}^*)\boldsymbol{\Psi}_n(\boldsymbol{\theta}^*)$. Therefore,
\begin{align*}
\sqrt{n}(\boldsymbol{M}^T\widehat{\boldsymbol{\theta}}_{\alpha}-\boldsymbol{m})\underset{n\rightarrow\infty}{\overset{L}{\rightarrow}}N(\boldsymbol{0}_{p},\boldsymbol{M}^T\boldsymbol{\Sigma}_{\alpha}(\boldsymbol{\theta}^*)\boldsymbol{M}).
\end{align*}
As $rank(\boldsymbol{M})=p$, we have that
\begin{align*}
n(\boldsymbol{M}^T\widehat{\boldsymbol{\theta}}_{\alpha}-\boldsymbol{m})^T(\boldsymbol{M}^T\boldsymbol{\Sigma}_{\alpha}(\boldsymbol{\theta}^*)\boldsymbol{M})^{-1}(\boldsymbol{M}^T\widehat{\boldsymbol{\theta}}_{\alpha}-\boldsymbol{m})
\end{align*}
converges in law to a chi-square distribution with $p$ degrees of freedom. But under $H_0$, $\boldsymbol{\Sigma}_{\alpha}(\boldsymbol{\theta}^0)=\boldsymbol{\Sigma}_{\alpha}(\boldsymbol{\theta}^*)$, and thus $\boldsymbol{W}_n^0(\boldsymbol{\theta}^0)$ converges  to a chi-square distribution with $p$ degrees of freedom.

\subsection{Proof of Theorem \ref{th:power1}}
A first-order Taylor expansion of $\ell(\boldsymbol{\theta})$ around $\boldsymbol{\theta}^*$ at $\widehat{\boldsymbol{\theta}}_{\alpha}$ is given by
\begin{align*}
\ell\widehat{\boldsymbol{\theta}}_{\alpha}-\ell(\boldsymbol{\theta}^*)=\left.\frac{\partial \ell(\boldsymbol{\theta})}{\partial \boldsymbol{\theta}^T} \right|_{\boldsymbol{\theta}=\boldsymbol{\theta}^*} (\widehat{\boldsymbol{\theta}}_{\alpha}-\boldsymbol{\theta}^*)+o_p(n^{-1/2}).
\end{align*}
Then the asymptotic distribution of the random variable $\sqrt{n}(\widehat{\boldsymbol{\theta}}_{\alpha}-\boldsymbol{\theta}^*)$ matches the asymptotic distribution of the random variable $\left.\frac{\partial \ell(\boldsymbol{\theta})}{\partial \boldsymbol{\theta}^T} \right|_{\boldsymbol{\theta}=\boldsymbol{\theta}^*} \sqrt{n}(\widehat{\boldsymbol{\theta}}_{\alpha}-\boldsymbol{\theta}^*$ and the result follows.

\subsection{Proof of Theorem \ref{th:asym_test_composite}}
We have by (\ref{eq:asym_est}) that
\begin{align*}
\sqrt{n}(\widehat{\boldsymbol{\theta}}_{\alpha}-\boldsymbol{\theta}^*)\underset{n\rightarrow\infty}{\overset{L}{\rightarrow}}N(\boldsymbol{0}_{p},\boldsymbol{\Sigma}_{\alpha}(\boldsymbol{\theta}^*)),
\end{align*}
where $\boldsymbol{\Sigma}_{\alpha}(\boldsymbol{\theta}^*)=\lim_{n\rightarrow\infty} \boldsymbol{\Psi}_n(\boldsymbol{\theta}^*) \boldsymbol{\Omega}_n^{-1}(\boldsymbol{\theta}^*)\boldsymbol{\Psi}_n(\boldsymbol{\theta}^*)$. Therefore,
\begin{align*}
\sqrt{n}(\boldsymbol{M}^T\widehat{\boldsymbol{\theta}}_{\alpha}-\boldsymbol{m})\underset{n\rightarrow\infty}{\overset{L}{\rightarrow}}N(\boldsymbol{0}_{p},\boldsymbol{M}^T\boldsymbol{\Sigma}_{\alpha}(\boldsymbol{\theta}^*)\boldsymbol{M}).
\end{align*}
As $rank(\boldsymbol{M})=r$, we have that
\begin{align*}
n(\boldsymbol{M}^T\widehat{\boldsymbol{\theta}}_{\alpha}-\boldsymbol{m})^T(\boldsymbol{M}^T\boldsymbol{\Sigma}_{\alpha}(\boldsymbol{\theta}^*)\boldsymbol{M})^{-1}(\boldsymbol{M}^T\widehat{\boldsymbol{\theta}}_{\alpha}-\boldsymbol{m})
\end{align*}
converges in law to a chi-square distribution with $r$ degrees of freedom. But $\boldsymbol{\Sigma}_{\alpha}(\widehat{\boldsymbol{\theta}}_{\alpha})$ is a consistent estimator of $\boldsymbol{\Sigma}_{\alpha}(\boldsymbol{\theta}^*)$, and thus $\boldsymbol{W}_n^0(\widehat{\boldsymbol{\theta}}_{\alpha})$ converges  to a chi-square distribution with $r$ degrees of freedom.

\subsection{Proof of Theorem \ref{thm:contiguous} \label{app:contiguous}}

A Taylor series expansion of $\boldsymbol{M}^T\boldsymbol{\theta} - \boldsymbol{m}$ around $\boldsymbol{\theta}_n$ yields
\begin{align*}
	\boldsymbol{M}^T\widehat{\boldsymbol{\theta}}_\alpha - \boldsymbol{m} & =  \boldsymbol{M}^T\boldsymbol{\theta}_n - \boldsymbol{m} + \boldsymbol{M}^T(\widehat{\boldsymbol{\theta}}_\alpha - \boldsymbol{\theta}_n) + o(||\widehat{\boldsymbol{\theta}}_\alpha-\boldsymbol{\theta}_n||_1 )\\
	&= n^{-1/2}\boldsymbol{M}^T\boldsymbol{d} - \boldsymbol{m} + \boldsymbol{M}^T(\widehat{\boldsymbol{\theta}}_\alpha - \boldsymbol{\theta}_n) + o(||\widehat{\boldsymbol{\theta}}_\alpha-\boldsymbol{\theta}_n||_1).
\end{align*}
Using Theorem \ref{thm:asymptotic}, $\sqrt{n}(\widehat{\boldsymbol{\theta}}_\alpha - \boldsymbol{\theta}_n) \underset{n\rightarrow\infty}{\overset{L}{\rightarrow}} N(\boldsymbol{0},\boldsymbol{\Sigma}_{\alpha})$ and $\sqrt{n}o(||\widehat{\boldsymbol{\theta}}_\alpha-\boldsymbol{\theta}_n||_1 ) = o_p(1),$ we get the asymptotic convergence
$$\sqrt{n}\left(\boldsymbol{M}^T\widehat{\boldsymbol{\theta}}_\alpha \right) \underset{n\rightarrow\infty}{\overset{L}{\rightarrow}} N(\boldsymbol{M}^T\boldsymbol{d}, \boldsymbol{M}^T\boldsymbol{\Sigma}_{\alpha}\boldsymbol{M}). $$
We now consider the random variable $\boldsymbol{Z} = \sqrt{n}\boldsymbol{M}^T\widehat{\boldsymbol{\theta}}_\alpha\left(\boldsymbol{M}^T\boldsymbol{\Sigma}_{\alpha}\boldsymbol{M}\right)^{-1/2}$ satisfying
$$\boldsymbol{Z} \underset{n\rightarrow\infty}{\overset{L}{\rightarrow}} N(\left(\boldsymbol{M}^T\boldsymbol{\Sigma}_{\alpha}\boldsymbol{M}\right)^{-1/2} \boldsymbol{M}^T\boldsymbol{d}, \boldsymbol{I}_{r\times r}).$$
 Hence, the asymptotic distribution of the quadratic form $\boldsymbol{W} = \boldsymbol{Z}^T\boldsymbol{Z}$ is given  by  a non-central chi-square distribution with $r$ degrees of freedom and non-centrality parameter 
 $$\delta = \boldsymbol{d}^T \boldsymbol{M}[\boldsymbol{M}^T\boldsymbol{\Sigma}_{\alpha}(\widehat{\boldsymbol{\theta}}_{\alpha})\boldsymbol{M}]^{-1} \boldsymbol{M}^T \boldsymbol{d}.$$

\subsection{Proof of Theorem \ref{th:IF} \label{app:IF}}

The IF of the functional $\boldsymbol{T}_{\alpha}(G_1,\dots,G_n)$ with contamination in the $i_0$-th direction will be obtained replacing $\boldsymbol{\theta}_{\varepsilon}^{i_0}$ and $g_{i_0,\varepsilon}$ in $\boldsymbol{\theta}$ and $g_{i_0}$ respectively in the equality (\ref{eq:IF2}), differentiating with respect to $\varepsilon$ and evaluating the corresponding equality in $\varepsilon=0$. 

In (\ref{eq:IF2}) we replace $\boldsymbol{\theta}$ by $\boldsymbol{\theta}_{\varepsilon}^{i_0}$ and $g_{i_0}(y)$ by
$$
g_{i_0,\varepsilon}=(1-\varepsilon)g_{i_0}(y)+\varepsilon \Delta_{t_{i_0}}(y),
$$
and for $i \neq i_0$ we consider the original $g_i(y)$. We get

\begin{align} \label{eq:IF3}
\frac{1}{n} {\textstyle\sum\limits_{i=1}^{n}} & \dfrac{\int f_i(y,\boldsymbol{\theta}_{\varepsilon}^{i_0})^{1+\alpha}\boldsymbol{u}_i(y,\boldsymbol{\theta}_{\varepsilon}^{i_0})dy}{\int f_i(y,\boldsymbol{\theta}_{\varepsilon}^{i_0})^{1+\alpha}dy}- \frac{1}{n} {\textstyle\sum\limits_{i\neq i_0}^{n}}  \dfrac{\int f_i(y,\boldsymbol{\theta}_{\varepsilon}^{i_0})^{\alpha}g_i(y)\boldsymbol{u}_i(y,\boldsymbol{\theta}_{\varepsilon}^{i_0})dy}{\int f_i(y,\boldsymbol{\theta}_{\varepsilon}^{i_0})^{\alpha}g_i(y)dy} \notag \\
&-   \dfrac{\int f_{i_0}(y,\boldsymbol{\theta}_{\varepsilon}^{i_0})^{\alpha}g_{i_0,\varepsilon}(y)\boldsymbol{u}_{i_0}(y,\boldsymbol{\theta}_{\varepsilon}^{i_0})dy}{\int f_{i_0}(y,\boldsymbol{\theta}_{\varepsilon}^{i_0})^{\alpha}g_{i_0,\varepsilon}(y)dy}=\boldsymbol{0}_p.
\end{align}
Now, we denote

\begin{align*}
\boldsymbol{\zeta}_{i,\alpha}(\boldsymbol{\theta}_{\varepsilon}^{i_0})&=\frac{\int f_i(y,\boldsymbol{\theta}_{\varepsilon}^{i_0})^{1+\alpha}\boldsymbol{u}_i(y,\boldsymbol{\theta}_{\varepsilon}^{i_0})dy}{\int f_i(y,\boldsymbol{\theta}_{\varepsilon}^{i_0})^{1+\alpha}dy},\\
\boldsymbol{\zeta}^*_{i,\alpha}(\boldsymbol{\theta}_{\varepsilon}^{i_0})&=\frac{\int f_i(y,\boldsymbol{\theta}_{\varepsilon}^{i_0})^{\alpha}g_i(y)\boldsymbol{u}_i(y,\boldsymbol{\theta}_{\varepsilon}^{i_0})dy}{\int f_i(y,\boldsymbol{\theta}_{\varepsilon}^{i_0})^{\alpha}g_i(y)dy},\\
\boldsymbol{\zeta}^{**}_{i,\alpha}(\boldsymbol{\theta}_{\varepsilon}^{i_0})&=\frac{\int f_{i_0}(y,\boldsymbol{\theta}_{\varepsilon}^{i_0})^{\alpha}g_{i_0,\varepsilon}(y)\boldsymbol{u}_{i_0}(y,\boldsymbol{\theta}_{\varepsilon}^{i_0})dy}{\int f_{i_0}(y,\boldsymbol{\theta}_{\varepsilon}^{i_0})^{\alpha}g_{i_0,\varepsilon}(y)dy}.
\end{align*}
Therefore, (\ref{eq:IF3}) can be written as,
\begin{align}\label{eq:IF4}
\frac{1}{n}\sum\limits_{i=1}^{n} \boldsymbol{\zeta}_{i,\alpha}(\boldsymbol{\theta}_{\varepsilon}^{i_0})-\frac{1}{n} {\textstyle\sum\limits_{i\neq i_0}^{n}}\boldsymbol{\zeta}^*_{i,\alpha}(\boldsymbol{\theta}_{\varepsilon}^{i_0})- \boldsymbol{\zeta}^{**}_{i_0,\alpha}(\boldsymbol{\theta}_{\varepsilon}^{i_0})=\boldsymbol{0}.
\end{align}
Now, we have
\begin{align*}
\frac{\partial \boldsymbol{\zeta}_{i,\alpha}(\boldsymbol{\theta}_{\varepsilon}^{i_0})}{\partial \varepsilon}=&\left(\int f_i(y,\boldsymbol{\theta}_{\varepsilon}^{i_0})^{1+\alpha} dy  \right)^{-2}  \left\{ \left[ (1+\alpha)\int f_i(y,\boldsymbol{\theta}_{\varepsilon}^{i_0})^{\alpha}\frac{\partial f_i(y,\boldsymbol{\theta}_{\varepsilon}^{i_0})}{\partial \boldsymbol{\theta}_{\varepsilon}^{i_0}}\frac{\partial \boldsymbol{\theta}_{\varepsilon}^{i_0}}{\partial \varepsilon} \boldsymbol{u}_i(y, \boldsymbol{\theta}_{\varepsilon}^{i_0})dy \right. \right. \\
& \left.  \left. \quad  + \int f_i(y,\boldsymbol{\theta}_{\varepsilon}^{i_0})^{1+\alpha}\frac{\partial \boldsymbol{u}_i(y,\boldsymbol{\theta}_{\varepsilon}^{i_0})}{\partial \boldsymbol{\theta}_{\varepsilon}^{i_0})}\frac{\partial \boldsymbol{\theta}_{\varepsilon}^{i_0}}{\partial \varepsilon} dy \right] \int f_i(y,\boldsymbol{\theta}_{\varepsilon}^{i_0})^{1+\alpha} dy \right.\\
&\left. - \left[ (1+\alpha)\int f_i(y,\boldsymbol{\theta}_{\varepsilon}^{i_0})^{\alpha}\frac{\partial f_i(y,\boldsymbol{\theta}_{\varepsilon}^{i_0})}{\partial \boldsymbol{\theta}_{\varepsilon}^{i_0}}\frac{\partial \boldsymbol{\theta}_{\varepsilon}^{i_0}}{\partial \varepsilon} dy \right] \int f_i(y,\boldsymbol{\theta}_{\varepsilon}^{i_0})^{\alpha}\ dy\right\},
\end{align*}
and 
\begin{align*}
& \left.\frac{\partial \boldsymbol{\zeta}_{i,\alpha}(\boldsymbol{\theta}_{\varepsilon}^{i_0})}{\partial \varepsilon}\right|_{\varepsilon=0}=  \dfrac{IF(t_{i_0},\boldsymbol{T}_{\alpha},G_1,\dots,G_n)}{\left(\int f_i(y,\boldsymbol{\theta}_{\varepsilon}^{i_0})^{1+\alpha} dy  \right)^{2} }\\
& \times \left\{ \left[(1+\alpha)  \int f_{i}(y,\boldsymbol{\theta})^{\alpha+1} \boldsymbol{u}^T_{i}(y,\boldsymbol{\theta})\boldsymbol{u}_{i}(y,\boldsymbol{\theta})dy + \int f_{i}(y,\boldsymbol{\theta})^{\alpha+1} \frac{\partial \boldsymbol{u}_{i}(y,\boldsymbol{\theta})}{\partial \boldsymbol{\theta}}dy\right]\int f_{i}(y,\boldsymbol{\theta})^{\alpha+1} dy \right.\\
& \left. \quad  - (1+\alpha) \left(\int f_{i}(y,\boldsymbol{\theta})^{\alpha+1} \boldsymbol{u}_{i}(y,\boldsymbol{\theta}) dy \right)\left(\int f_{i}(y,\boldsymbol{\theta})^{\alpha+1} \boldsymbol{u}_{i}(y,\boldsymbol{\theta}) dy \right)^T\right\}\\
& = \dfrac{IF(t_{i_0},\boldsymbol{T}_{\alpha},G_1,\dots,G_n)}{\left(\int f_i(y,\boldsymbol{\theta}_{\varepsilon}^{i_0})^{1+\alpha} dy  \right)^{2} } \boldsymbol{A}_{i,\alpha}(\boldsymbol{\theta}),
\end{align*}
with
\begin{align*}
\boldsymbol{A}_{i,\alpha}(\boldsymbol{\theta})=&\left[(1+\alpha)  \int f_{i}(y,\boldsymbol{\theta})^{\alpha+1} \boldsymbol{u}^T_{i}(y,\boldsymbol{\theta})\boldsymbol{u}_{i}(y,\boldsymbol{\theta})dy + \int f_{i}(y,\boldsymbol{\theta})^{\alpha+1} \frac{\partial \boldsymbol{u}_{i}(y,\boldsymbol{\theta})}{\partial \boldsymbol{\theta}}dy\right]\int f_{i}(y,\boldsymbol{\theta})^{\alpha+1} dy.
\end{align*}
Therefore, 
\begin{align*}
\left. \frac{\partial }{\partial \varepsilon}\sum\limits_{i=1}^n\boldsymbol{\zeta}_{i,\alpha}(\boldsymbol{\theta}_{\varepsilon}^{i_0})\right|_{\varepsilon=0} &= \left. \sum\limits_{i=1}^n\frac{\partial \boldsymbol{\zeta}_{i,\alpha}(\boldsymbol{\theta}_{\varepsilon}^{i_0})}{\partial \varepsilon}\right|_{\varepsilon=0}\\
&=  IF(t_{i_0},\boldsymbol{T}_{\alpha},G_1,\dots,G_n)\sum\limits_{i=1}^n\frac{\boldsymbol{A}_{i,\alpha}(\boldsymbol{\theta})}{\left(\int f_i(y,\boldsymbol{\theta}_{\varepsilon}^{i_0})^{1+\alpha} dy  \right)^{2}}.
\end{align*}
Now,

\begin{align*}
\frac{\partial \boldsymbol{\zeta}^{*}_{i,\alpha}(\boldsymbol{\theta}_{\varepsilon}^{i_0})}{\partial \varepsilon}=&\left(\int f_i(y,\boldsymbol{\theta}_{\varepsilon}^{i_0})^{\alpha}g_i(y) dy  \right)^{-2}  \left\{ \left[ \alpha\int f_i(y,\boldsymbol{\theta}_{\varepsilon}^{i_0})^{\alpha-1}g_i(y)\frac{\partial f_i(y,\boldsymbol{\theta}_{\varepsilon}^{i_0})}{\partial \boldsymbol{\theta}_{\varepsilon}^{i_0}}\frac{\partial \boldsymbol{\theta}_{\varepsilon}^{i_0}}{\partial \varepsilon} \boldsymbol{u}_i(y, \boldsymbol{\theta}_{\varepsilon}^{i_0})dy \right. \right. \\
& \left.  \left. \quad  + \int f_i(y,\boldsymbol{\theta}_{\varepsilon}^{i_0})^{\alpha}g_i(y)\frac{\partial \boldsymbol{u}_i(y,\boldsymbol{\theta}_{\varepsilon}^{i_0})}{\partial \boldsymbol{\theta}_{\varepsilon}^{i_0})}\frac{\partial \boldsymbol{\theta}_{\varepsilon}^{i_0}}{\partial \varepsilon} dy \right] \int f_i(y,\boldsymbol{\theta}_{\varepsilon}^{i_0})^{\alpha}g_i(y) dy \right.\\
&\left. - \left[ \alpha\int f_i(y,\boldsymbol{\theta}_{\varepsilon}^{i_0})^{\alpha-1}g_i(y)\frac{\partial f_i(y,\boldsymbol{\theta}_{\varepsilon}^{i_0})}{\partial \boldsymbol{\theta}_{\varepsilon}^{i_0}}\frac{\partial \boldsymbol{\theta}_{\varepsilon}^{i_0}}{\partial \varepsilon} dy \right] \int f_i(y,\boldsymbol{\theta}_{\varepsilon}^{i_0})^{\alpha}dy\right\}.
\end{align*}
This is, 

\begin{small}
\begin{align*}
& \left.\frac{\partial \boldsymbol{\zeta}^{*}_{i,\alpha}(\boldsymbol{\theta}_{\varepsilon}^{i_0})}{\partial \varepsilon}\right|_{\varepsilon=0}=  \dfrac{IF(t_{i_0},\boldsymbol{T}_{\alpha},G_1,\dots,G_n)}{\left(\int f_i(y,\boldsymbol{\theta}_{\varepsilon}^{i_0})^{\alpha}g_i(y) dy  \right)^{2} }\\
& \times \left\{ \left[\alpha  \int f_{i}(y,\boldsymbol{\theta})^{\alpha}g_i(y) \boldsymbol{u}^T_{i}(y,\boldsymbol{\theta})\boldsymbol{u}_{i}(y,\boldsymbol{\theta})dy + \int f_{i}(y,\boldsymbol{\theta})^{\alpha}g_i(y) \frac{\partial \boldsymbol{u}_{i}(y,\boldsymbol{\theta})}{\partial \boldsymbol{\theta}}dy\right]\int f_{i}(y,\boldsymbol{\theta})^{\alpha}g_i(y) dy \right.\\
& \left.- \alpha \left(\int f_{i}(y,\boldsymbol{\theta})^{\alpha} g_i(y)\boldsymbol{u}_{i}(y,\boldsymbol{\theta}) dy \right)\left(\int f_{i}(y,\boldsymbol{\theta})^{\alpha}g_i(y) \boldsymbol{u}_{i}(y,\boldsymbol{\theta}) dy \right)^T \right\}\\
& = \dfrac{IF(t_{i_0},\boldsymbol{T}_{\alpha},G_1,\dots,G_n)}{\left(\int f_i(y,\boldsymbol{\theta}_{\varepsilon}^{i_0})^{\alpha}g_i(y) dy  \right)^{2} } \boldsymbol{A}^{*}_{i,\alpha}(\boldsymbol{\theta}),
\end{align*}
\end{small}
\noindent with
\begin{small}
\begin{align*}
\boldsymbol{A}^*_{i,\alpha}(\boldsymbol{\theta})=&\left[\alpha  \int f_{i}(y,\boldsymbol{\theta})^{\alpha}g_i(y) \boldsymbol{u}^T_{i}(y,\boldsymbol{\theta})\boldsymbol{u}_{i}(y,\boldsymbol{\theta})dy + \int f_{i}(y,\boldsymbol{\theta})^{\alpha}g_i(y) \frac{\partial \boldsymbol{u}_{i}(y,\boldsymbol{\theta})}{\partial \boldsymbol{\theta}}dy\right]\int f_{i}(y,\boldsymbol{\theta})^{\alpha}g_i(y) dy \\
&- \alpha \left(\int f_{i}(y,\boldsymbol{\theta})^{\alpha} g_i(y)\boldsymbol{u}_{i}(y,\boldsymbol{\theta}) dy \right)\left(\int f_{i}(y,\boldsymbol{\theta})^{\alpha}g_i(y) \boldsymbol{u}_{i}(y,\boldsymbol{\theta}) dy \right)^T.
\end{align*}
\end{small}
Therefore, 
\begin{align*}
\left. \frac{\partial }{\partial \varepsilon}\sum\limits_{i=1}^n\boldsymbol{\zeta}^*_{i,\alpha}(\boldsymbol{\theta}_{\varepsilon}^{i_0})\right|_{\varepsilon=0} &= \left. \sum\limits_{i=1}^n\frac{\partial \boldsymbol{\zeta}^*_{i,\alpha}(\boldsymbol{\theta}_{\varepsilon}^{i_0})}{\partial \varepsilon}\right|_{\varepsilon=0}\\
&=  IF(t_{i_0},\boldsymbol{T}_{\alpha},G_1,\dots,G_n)\sum\limits_{i=1}^n\frac{\boldsymbol{A}^*_{i,\alpha}(\boldsymbol{\theta})}{\left(\int f_i(y,\boldsymbol{\theta}_{\varepsilon}^{i_0})^{\alpha}g_i(y) dy  \right)^{2}}.
\end{align*}
In a similar manner,
\begin{align*}
\left. \frac{\partial }{\partial \varepsilon}\boldsymbol{\zeta}^{**}_{i_0,\alpha}(\boldsymbol{\theta}_{\varepsilon}^{i_0})\right|_{\varepsilon=0} &= \left. \frac{\partial \boldsymbol{\zeta}^{**}_{i_0,\alpha}(\boldsymbol{\theta}_{\varepsilon}^{i_0})}{\partial \varepsilon}\right|_{\varepsilon=0}= \dfrac{\boldsymbol{\ell}_{i_0,\alpha}(\boldsymbol{\theta})}{\left( \int f_{i_0}(y,\boldsymbol{\theta})^{\alpha} g_{i_0}(y)dy\right)^2}, 
\end{align*}
with
\begin{align*}
\boldsymbol{\ell}_{i_0,\alpha}(\boldsymbol{\theta})=f_{i_0}(y,\boldsymbol{\theta}) \int f_{i_0}(y,\boldsymbol{\theta})^{\alpha+1} \boldsymbol{u}_{i_0}(y,\boldsymbol{\theta})dy -f_{i_0}(y,\boldsymbol{\theta})\boldsymbol{u}_{i_0}(y,\boldsymbol{\theta}) \int f_{i_0}(y,\boldsymbol{\theta})^{\alpha+1} dy.
\end{align*}

Therefore, equality (\ref{eq:IF4}) can be written as
\begin{align*}
& IF(t_{i_0},\boldsymbol{T}_{\alpha},G_1,\dots,G_n)\left\{\frac{1}{n}\sum_{i=1}^n\left[\dfrac{\boldsymbol{A}_{i,\alpha}(\boldsymbol{\theta})}{\left( \int f_{i}(y,\boldsymbol{\theta})^{1+\alpha} dy\right)^2}-\dfrac{\boldsymbol{A}^*_{i,\alpha}(\boldsymbol{\theta})}{\left( \int f_{i}(y,\boldsymbol{\theta})^{\alpha} g_{i}(y)dy\right)^2}\right] \right\}\\
&\quad +\dfrac{\boldsymbol{\ell}_{i_0,\alpha}(\boldsymbol{\theta})}{\left( \int f_{i_0}(y,\boldsymbol{\theta})^{\alpha} g_{i_0}(y)dy\right)^2}=\boldsymbol{0}.
\end{align*}
Finally,
\begin{align*}
IF(t_{i_0}, \boldsymbol{T}_{\alpha},G_1,\dots,G_n)=(\boldsymbol{M}_{n,\alpha}(\boldsymbol{\theta}))^{-1}\dfrac{-\boldsymbol{\ell}_{i_0,\alpha}(\boldsymbol{\theta})}{\left( \int f_{i_0}(y,\boldsymbol{\theta})^{\alpha} g_{i_0}(y)dy\right)^2},
\end{align*}
where
\begin{align*}
\boldsymbol{M}_{n,\alpha}(\boldsymbol{\theta})=\frac{1}{n}{\sum\limits_{i=1}^{n}}\left[\dfrac{\boldsymbol{A}_{i,\alpha}(\boldsymbol{\theta})}{\left( \int f_{i}(y,\boldsymbol{\theta})^{1+\alpha} dy\right)^2}-\dfrac{\boldsymbol{A}^*_{i,\alpha}(\boldsymbol{\theta})}{\left( \int f_{i}(y,\boldsymbol{\theta})^{\alpha} g_{i}(y)dy\right)^2}\right].
\end{align*}

\subsection{Proof of Lemma \ref{lemma:MLRM}}

This proof is very similar to that of \cite{Ghosh2013} (Lemma 6.1).

\subsection{Proof of Theorem \ref{thm:asymptotic}}

The consistence follows directly from Lemma \ref{lemma:MLRM} and Theorem \ref{Th1}, while the asymptotic distribution is obtained applying Lemma \ref{lemma:MLRM}  and Theorem \ref{Th2} to the MLRM.

%\clearpage

\end{document}